
%
%
%
%
%
%
\input amstex
\documentstyle{amsppt}
\loadbold
\def\cstar{$C^*$-algebra}
\def\esg{$E_0$-semigroup}

\def\<{\left<}										
\def\>{\right>}
\define\sumplus{\sideset \and^\oplus \to\sum}
\magnification=\magstep 1

\topmatter
\title Path Spaces, continuous tensor products, and
$E_0$-semigroups
\endtitle

\author William Arveson
\endauthor

\affil Department of Mathematics\\
University of California\\Berkeley CA 94720, USA
\endaffil

\date 14 November 1994
\enddate
\thanks This research was supported in part by
NSF grant DMS92-43893
\endthanks
\abstract
We classify all continuous
tensor product systems of Hilbert spaces which are
``infinitely divisible" in the sense that they have
an associated logarithmic structure.  These results
are applied to the theory of \esg s to deduce that every
\esg\ possessing sufficiently many ``decomposable" operators
must be cocycle conjugate to a $CCR$ flow.

A {\it path space} is an abstraction of the set
of paths in a topological space, on which there is given
an associative rule of concatenation.  A {\it metric path space}
is a pair $(P,g)$ consisting of a path space $P$ and a function
$g:P^2\to \Bbb C$ which behaves as if it were
the logarithm of a multiplicative
inner product.
The logarithmic structures associated with
infinitely divisible product systems are such objects.
The preceding results are based on a classification
of metric path spaces.
\endabstract


\toc
\specialhead{} Introduction
\endspecialhead
\specialhead{} Part I. Path Spaces
\endspecialhead
\subhead 1.  Definitions and examples
\endsubhead
\subhead 2.  Additive forms and multiplicative forms
\endsubhead
\subhead 3.  Exactness of cocycles
\endsubhead
\subhead 4.  Classification of additive forms
\endsubhead
\specialhead{} Part II. Continuous Tensor Products
\endspecialhead
\subhead 5. Continuity of the modulus
\endsubhead
\subhead 6.  Inner products of decomposable vectors
\endsubhead
\subhead 7.  Continuity and normalization
\endsubhead
\subhead 8.  Continuous logarithms
\endsubhead
\subhead 9.  Infinite divisibility
\endsubhead
\subhead 10.  Existence of measurable propagators
\endsubhead
\specialhead{} Part III. Applications
\endspecialhead
\subhead 11. Decomposable continuous tensor products
\endsubhead
\subhead 12.  Decomposable $E_0$-semigroups
\endsubhead
\endtoc

\endtopmatter
\vfill\eject

\document

\subheading{Introduction}
Let $\alpha = \{\alpha_t: t>0\}$ be an \esg.  That is, each
$\alpha_t$ is a normal unit-preserving $*$-endomorphism of
$\Cal B(H)$ such that $\alpha_s\circ\alpha_t = \alpha_{s+t}$,
and which is continuous in the sense that for all $A\in\Cal B(H)$,
$\xi,\eta\in H$, $\<\alpha_t(A)\xi,\eta\>$ is
continuous in $t$.

A {\it unit} for $\alpha$ is a strongly continuous semigroup
$U = \{U(t): t\geq 0\}$ of bounded operators on $H$ satisfying
$$
\alpha_t(A)U(t) = U(t) A,
$$
for every $t\geq 0$, $A\in\Cal B(H)$.  $\alpha$ is called
{\it completely spatial} if there is a  $t>0$ such that $H$ is
spanned by the ranges of all operator products of the form
$$
U_1(t_1)U_2(t_2)\dots U_n(t_n)
$$
where $U_1, U_2,\dots,U_n$ are units, $t_1, t_2,\dots,t_n$ are
nonnegative reals summing to $t$, and $n$ is an arbitrary positive
integer.  Completely spatial \esg s are those which have
``sufficiently many" units.

In \cite{2, section 7}, completely spatial \esg s are shown to be completely
classified up to cocycle conjugacy by their numerical index.  This
is established at the level of continuous tensor product systems.  In
somewhat more detail, a {\it product system} is a measurable family
of separable Hilbert spaces $E = \{E_t: t>0\}$ which is endowed with
an associative multiplication that ``acts like tensoring" in the sense
that for any choices $x, x^\prime\in E_s$, $y,y^\prime\in E_t$, the
products $xy$, $x^\prime y^\prime$ both belong to $E_{s+t}$ and we have
$$
\align
\<xy,x^\prime y^\prime\> &=\<x,x^\prime\> \<y,y^\prime\>,\qquad {\text and} \\
E_{s+t} &= \overline{span} E_s E_t.
\endalign
$$
In addition, there is a third key axiom which plays the role of local
triviality for Hermitian vector bundles.  The precise definition of
product systems can be found in \cite{2, Definition 1.4}.
The intuition is that
$E_t$ resembles a continuous tensor product
$$
E_t = \bigotimes_{0<s<t} H_s,\qquad H_s = H \tag{I.1}
$$
of copies of a single Hilbert space $H$.  However, this
heuristic picture of continuous tensor products cannot be pushed too far.
While the formula (I.1) can be made precise for certain
standard examples, it is also known that there are many product
systems for which the ``germ" $H$ fails to exist.

Every \esg\ $\alpha$ determines a product system $E_\alpha$.  The
$t^{\text th}$ Hilbert space $E_\alpha(t)$ is the linear space of operators
$$
E_\alpha(t) = \{T\in\Cal B(H): \alpha_t(A)T = TA, \forall A\in\Cal B(H)\},
$$
with inner product defined by
$$
\<S,T\>{\bold 1} = T^*S,\qquad S,T\in E_\alpha(t),
$$
and with multiplication given by ordinary operator multiplication.  It is
known that two \esg s are cocycle conjugate iff their product systems
$E_\alpha$, $E_\beta$ are isomorphic \cite{2, Corollary of Theorem 3.18};
moreover, every abstract
product system is associated with some \esg\ \cite{5, Corollary 5.17}.

These remarks show that, up to cocycle conjugacy, the
theory of \esg s is equivalent to the
theory of continuous tensor product systems.  Consequently,
a central component of our approach to \esg s has been to develop
the theory of product systems.  The classification of \esg s
described above was obtained by classifying product systems which possess
a property that corresponds to complete spatiality of their product systems
(this was called {\it divisibility\/} in \cite{2}).

In this paper, we extend that result so as to include certain product
systems which do not appear to contain any units {\it a priori}, but which
do contain sufficiently many vectors that resemble ``elementary tensors".
Such product systems are called {\it decomposable}.  The corresponding
property of \esg s is described as follows.  Fix $t>0$.  An operator
$T \in E_\alpha(t)$ is called {\it decomposable} if, for every
$0<s<t$ there are operators $A\in E_\alpha(s)$, $B\in E_\alpha(t-s)$
such that $T = AB$.  Let $\Cal D(t)$ denote the set of all decomposable
operators in $E_\alpha(t)$.  It is easy to see that if $H$ is spanned
by vectors of the form
$$
\Cal D(t_0)H = \{T\xi: T\in\Cal D(t_0),\xi\in H\}
$$
for some particular $t_0 > 0$, then we have
$$
H = [\Cal D(t)H] \tag{I.2}
$$
for every $t > 0$.  $\alpha$ is called {\it decomposable} if (I.2) is
valid for some $t>0$.  Notice that if $U_1, U_2, \dots, U_n$ are units
for $\alpha$ and $t_1, t_2,\dots,t_n$ are positive numbers summing to
$t$, then every operator product of the form
$$
U_1(t_1)U_2(t_2)\dots U_n(t_n)
$$
is a decomposable operator in $E_\alpha(t)$ because each $U_j$ is a
semigroup.  Thus, every completely spatial \esg\ is decomposable.

On the other hand, decomposable \esg s are not required to contain any
units {\it a priori}.  The results of this paper imply that decomposable
\esg s are completely spatial, and are therefore classified to cocycle
conjugacy by their numerical index.  In particular, decomposable \esg s
necessarily have plenty of units.

Decomposability translates into an important property of abstract
product systems, and we want to discuss this property and its
role in structural issues.  Let $E = \{E(t): t>0\}$ be an abstract
product system and fix $t>0$.  A nonzero vector $x\in E(t)$ is called
{\it decomposable} if for every $s\in(0,t)$ there are vectors
$a\in E(s)$, $b\in E(t-s)$ such that $x = ab$.  Let $D(t)$ be the
set of all decomposable vectors in $E(t)$.  $D(t)$ can be empty.
But if it is large enough that it spans $E(t)$ for some particular
$t$, then it is easy to show that
$$
E(t) = \overline{span} D(t) \tag{I.3}
$$
for every positive $t$.  A product system satisfying (I.3) for every
positive $t$ is called {\it decomposable}.  It is also easy to show
that if $E$ is isomorphic to the product system $E_\alpha$  of some
\esg\ $\alpha$, then $E$ is decomposable iff $\alpha$ is (Proposition 12.2).

If we think of the multiplication in a product system $E$ as representing
the tensor product operation, then the vectors in $D(t)$ are ``elementary
tensors".  In the heuristic picture in which one thinks of $E(t)$ as a
continuous tensor product
$$
E(t) = \bigotimes_{0 < s < t} H_s, \qquad H_s = H
$$
of copies of a Hilbert space $H$, the decomposable vectors have a
corresponding heuristic representation
$$
x = \operatornamewithlimits\otimes_{0 < s < t} x_s \tag{1.4}
$$
where $x_s$ is a nonzero vector in $H$ for every $s\in (0,t)$.  Since
finite tensor products of Hilbert spaces of the form
$$
\underbrace{H\otimes H\otimes\dots\otimes H}_{n \text{ times}}
$$
are always spanned by vectors of the form
$$
x = x_1\otimes x_2\otimes\dots \otimes x_n
$$
with $x_i\in H$, decomposability appears as a natural property of
product systems.  However, notice that Theorem 11.1 below implies
that only the simplest product systems are decomposable.

Product systems occur naturally in probability theory.  For example,
any random distribution over the interval $(0,\infty)$ having
stationary independent increments gives rise to a product system
(see \cite{2, pp 14--16} for an example based on Poisson processes).  Some
of these examples (such those arising from Poisson processes) do not
appear to have ``enough" units to be divisible in the sense of \cite{2}.
Nevertheless, they do.  They do because in constructing these examples,
one essentially constructs the sets $D(t)$, $t>0$, and then defines
the product system $E$ by setting
$$
E(t) = \overline{span} D(t), \qquad t>0.
$$
Such a product system is clearly decomposable and thus, by Theorem 11.1
below, is simply one of the standard ones.

Our approach to decomposable product systems is
based on the concept of a {\it metric path space}.  Roughly
speaking, a metric path space is a structure which behaves as
if it were the ``logarithm" of a product system.  Path spaces are
very general objects, and we give a variety of examples.
The most natural constructions of product systems all begin with a
metric path space of some kind...an example of this being
the construction starting with a Poisson process as above.
We give an essentially complete
classification of metric path spaces in part I.  In part II,
we show that every decomposable product system is associated
with an essentially unique metric path space, and in part III
these results are applied to the classification of product
systems and  \esg s.

We still do not know a ``natural" construction of a product system
without any units, or even one with some units but not one of the
standard ones.  Such product systems are known to exist by work of
Powers \cite{25, 27} together with the results of \cite{2}.  Powers'
constructions are very indirect.  We believe that there should be a
natural way of constructing such product systems, and we offer
that as a basic unsolved problem.  The fact that we do not yet know
how to solve it shows how poorly understood continuous tensor
products are today.

It is appropriate to discuss the relation of our results
to the work of Araki and Woods on continuous tensor products \cite{1}.
Araki and Woods were interested in generalizing results of von Neumann
on infinite tensor products of Hilbert spaces
\cite{21,22} to the case of continuous tensor products.  Their
precise formulation of a continuous tensor product of Hilbert spaces
is expressed in terms of a nonatomic Boolean algebra of type I factors.
Our formulation of continuous tensor products
is in terms of the product systems defined above.
While both concepts seek to make precise the intuitive
idea of continuous tensor products of Hilbert spaces,
they are formulated in substantially different and
mathematically inequivalent ways.

For example,
the Hilbert spaces of \cite{1} are parameterized by the elements of
a $\sigma$-algebra of sets, whereas the fiber spaces of product systems
are parameterized by positive real numbers.  If one starts with an
abstract product system $E$ then there is no underlying Hilbert space
in view.  Even if $E$ did act naturally on some Hilbert space (as it
does when it is defined in terms of an \esg), one still cannot use it to
define a Boolean algebra of type I factors.  It is true that one can
associate a type I factor with subintervals of $(0,\infty)$, but
one does not know how to extend that set function to arbitrary Borel sets in
$(0,\infty)$...certainly not so that one can be sure
the extended set function will take values
in the set of type I factors.

Thus, product systems are more general than the structures of
Araki-Woods.  But even that assertion is not exactly true.  We also
do not know how to start with a complete Boolean algebra of type I
factors over the interval $(0,\infty)$ and use that structure to define a
product system.  In order to do that, one needs a bit more symmetry
with respect to positive translations.  Still, I am comfortable
with the philosophical position that product
systems are generalizations of complete Boolean
algebras of type I factors over the interval $(0, \infty)$.

In more pragmatic terms, we have been unable to apply the work of
Araki-Woods to the theory of product systems or to the theory of
\esg s, even though there are apparent similarities in
certain results.  For example, Theorem 11.1 below is similar
to Theorem 6.1 of \cite{1}.  Notice that the hypothesis
of the latter, namely
$$
\inf_B d(\Psi; B) > 0,
$$
corresponds to our hypothesis of decomposability because of
Lemma 6.3 of \cite{1}.  In fact, in establishing the classification
results of \cite{2} as well as those of the present paper,
we have had to start from scratch...basing our proofs on different principles
than the proofs of corresponding results of \cite{1}.
Nevertheless, we certainly acknowledge the considerable
influence of \cite{1} in guiding
our thinking about continuous tensor products.

\vfill
\eject

\heading{Part I: Path Spaces}
\endheading

\subheading{1.  Definitions and examples}

A {\it path space} is a family of sets over the
positive real axis
$$
\pi: P\to (0, \infty)
$$
which is endowed with an associative multiplication
$(x,y)\in P\times P \mapsto xy\in P$ satisfying the following
two axioms.  Letting $P(t) = \pi^{-1}(t)$ be the fiber
over $t$ we require that for every $s,t>0$,
$$
P(s)\cdot P(t) = P(s+t)\tag{1.1.1}
$$
and that for all $x_1$, $x_2\in P(s)$, $y_1$, $y_2\in P(t)$
$$
x_1y_1=x_2y_2 \implies x_1=x_2 \text{ and } y_1=y_2.\tag{1.1.2}
$$

Together, the conditions (1.1) assert that the projection
$\pi$ should obey $\pi(xy) = \pi(x)+\pi(y)$, and that for fixed
$0<s<t$, each element $z\in P(t)$ factors uniquely into a
product $z=xy$, with $x\in P(s)$ and $y\in P(t-s)$.

Path spaces are generalizations of continuous cartesian
products, whereby one starts with a basic set $X$ and
defines fiber spaces $P(t)$ by
$$
P(t) = X^{(0,t]}.\tag{1.2.1}
$$
$P$ is defined as the total set
$\pi: P=\{(t,f): t>0, f\in P(t)\}\to (0,\infty)$
with its natural projection $\pi(t,f) = t$.
The multiplication in $P$ is defined by
$$
(s,f)\cdot (t,g) = (s+t,h)
$$
where $h$ is the concatenation
$$
h(x) =
\cases
f(x),& 0<x\leq s\\
g(x-s),& s<x\leq s+t.
\endcases\tag{1.2.2}
$$

\remark{Remarks}
The terminology {\it path space} suggests that
one should think of the elements of $P$ as representing
paths in some space.  Notice however that
if $X$ is a topological
space in the preceding example (1.2), the
functions in any particular $P(t)$ are not necessarily
continuous, nor even measurable.
In fact, the examples (1.2) are pathological.
We give some more typical examples
below.  But we emphasize that in general,
path spaces are not required
to support an additional topological structure,
nor even a Borel structure.  We will find that,
nevertheless, the intrinsic structure of path
spaces allows one to define a notion of
measurability (see Definition 2.2) and even continuity
(see section 8) for the functions defined on them.

Let $\pi:P\to (0,\infty)$ be a path space and
let $x,y\in P$.  We will say that $x$ is a left
divisor of $y$ is there is an element $z\in P$ with
$y=xz$.  The element $z$ is necessarily unique and,
if $x\in P(s)$ and $y\in P(t)$ then we must have
$0<s<t$ and $z\in P(t-s)$.  Fix $T>0$.  A section
$$
t\in (0,T] \mapsto x(t)\in P(t)
$$
is said to be {\it left-coherent} if $x(s)$ is
a left divisor of $x(t)$ whenever $s<t$.  Notice
that such a section is uniquely determined by
its last element $x(T)\in P(T)$ because of the
unique factorization property.  Thus there is
a bijective correspondence between elements
of $P(T)$ and left-coherent sections
$t\in (0,T]\mapsto x(t)\in P(t)$.

There is another description of left
coherent sections $t\in (0,T]\mapsto P(t)$ that
we will find useful.  Given such a section
we can define a family of elements
$\{x(s,t)\in P(t-s): 0\leq s<t\leq T\}$ by
$$
x(t) = x(s)x(s,t),\qquad \text{for }0<s<t\leq T,
$$
and where $x(0,t)$ is defined as $x(t)$.  This family
satisfies the consistency equation
$$
0\leq r<s<t\leq T \implies x(r,t) = x(r,s)x(s,t).  \tag{1.3}
$$
There is an obvious modification of these considerations
for fully defined sections $t\in (0,\infty)\mapsto x(t)\in P(t)$,
and in this case we obtain a family $\{x(s,t): 0\leq s<t<\infty\}$
satisfying (1.3) for unrestricted $0\leq r<s<t$.

\proclaim{Definition 1.4}
Let $I$ be an interval of the form $I=[0,T]$ with
$0<T<\infty$ or $I = [0,\infty)$.  A family
of elements $\{x(s,t)\in P(t-s): s,t\in I, s<t\}$ is called a
propagator on $I$ if it satisfies the equation (1.3) for
every $r,s,t\in I$ such that $r<s<t$.
\endproclaim

Operator-valued functions satisfying equation (1.3) are
naturally associated with solutions of time-dependent
linear differential equations of the form
$$
\frac{d}{dt}x(t) = x(t) a(t),\qquad 0\leq t \leq T,
$$
where $a(\cdot)$ is a given operator-valued function with
invertible values \cite{28, p. 282}.  For our
purposes, propagators are associated with left-coherent
sections as above.  Notice that the correspondence between
propagators and left-coherent sections is also bijective
since one can recover the section $\{x(t)\}$ from the
propagator $\{x(s,t)\}$ by setting
$$
x(t) = x(0,t).
$$

\example{Example 1.5}Let $V$ be a finite dimensional
vector space or a Banach space.  For every $t>0$, let
$P(t)$ denote the space of all continuous functions
$$
f:[0,t]\to V
$$
satisfying $f(0)=0$.  For $f\in P(s)$, $g\in P(t)$ we
define a concatenation $f*g \in P(s+t)$ by
$$
f*g(\lambda) =
\cases
f(\lambda),& 0\leq\lambda < s \\
f(s) + g(\lambda-s),& s\leq\lambda\leq s+t.
\endcases
$$
Notice that $f*g$ executes the path $f$ first, and then
it executes $g$.  If we assemble these spaces into a family
$$
P = \{(t,f): t>0, f\in P(t)\}
$$
with projection $\pi(t,f) = t$ and multiplication
$$
(s,f)\cdot(t,g) = (s+t, f*g),
$$
then we obtain a path space $\pi:P\to (0,\infty)$.
\endexample

\example{Example 1.6} One can define variations of example
1.5 with higher degrees of smoothness.  For example, let
$P(t)$ be the space of all continuously differentiable
functions $f: [0,t]\to V$ satisfying $f(0) = f^\prime(0) = 0$.
Then with the concatenation rule in which for $f\in P(s)$
and $g\in P(t)$, $f*g$ is defined by
$$
f*g(\lambda) =
\cases
f(\lambda),& 0\leq \lambda<s\\
f(s) + (\lambda-s)f^\prime(s) +g(\lambda-s),& s<\lambda\leq s+t
\endcases
$$
we obtain a path space structure on
$P = \{(t,f): t>0, f\in P(t)\}$ by way of
$\pi(t,f) = t$, $(s,f)\cdot(t,g) = (s+t,f*g)$.
\endexample

\example{Example 1.7}  As we will see in section 4, the most
important examples of path spaces arise as follows.

Let $\Cal C$ be a separable Hilbert space, and consider
the Hilbert space $L^2((0,\infty); \Cal C)$ of all square
integrable vector valued functions $f:(0,\infty)\to \Cal C$
with inner product
$$
\<f,g\> = \int_0^\infty \<f(x),g(x)\>\, dx.
$$
For every $t>0$, let $\Cal P_{\Cal C}(t)$ denote the subspace of
$L^2((0,\infty);\Cal C)$ consisting of all functions
$f$ satisfying $f(x) = 0$ almost everywhere for
$x\geq t$.  For $f\in \Cal P_{\Cal C}(s)$,
$g\in \Cal P_{\Cal C}(t)$ we define
the concatenation $f\boxplus g\in \Cal P_{\Cal C}(s+t)$ by
$$
f\boxplus g(\lambda) =
\cases
f(\lambda),& 0<\lambda<s\\
g(\lambda-s),& s\leq\lambda\leq s+t.
\endcases
$$
Let $\pi:\Cal P_{\Cal C}\to (0,\infty)$ be the
total space defined by
$$
\Cal P_{\Cal C} = \{(t,f): t>0, f\in \Cal P_{\Cal C}(t)\},\tag{1.8}
$$
where $\pi(t,f)=t$.  With multiplication
$(s,f)\cdot(t,g) = (s+t,f\boxplus g)$,
$\pi:\Cal P_{\Cal C}\to (0,\infty)$ becomes a path space.
\endexample

Needless to say, there are obvious variations of these
examples which can be formulated in different function spaces.

\remark{Remarks}  The preceding examples
involve generalized paths in vector
spaces.  It is less obvious how one might concatenate paths in
multiply connected spaces so as to obtain a path space structure.
For example, suppose we are given a
closed set $K\subseteq \Bbb R^2$ which represents an
``obstruction" and one defines $P(t)$ to be the set of all
continuous functions
$$
f: [0,t]\to \Bbb R^2\setminus K
$$
satisfying $f(0)=0$ (we assume, of course, that $0\notin K$).
One may not concatenate elements according to the
rule of example 1.5 since that rule of concatenation
will not necessarily provide a path $f*g$ that avoids $K$, even if
$f$ and $g$ both avoid $K$.  Nevertheless, there are
many ways to define such concatenation operations.  We offer
the following digression for the interested reader, though
we will make no further reference to this class of examples in
this paper.
\endremark

\example{Example 1.9}  Let $K$ be a closed subset of
$\Bbb R^n$, $n\geq 2$, having smooth boundary and which
does not contain $0$.  Let
$$
V:\Bbb R^n\setminus K\to \Bbb R
$$
be any $C^1$ function which tends to $+\infty$ near the
boundary of $K$ in the sense that  if $x_n$ is any sequence
for which the distance from $x_n$ to $K$ tends to $0$,
then $V(x_n)\to +\infty$.  For example, one may take
$$
V(x) = \sup_{y\in K}\frac{1}{|x-y|},\qquad x\in \Bbb R^n\setminus K.
$$
For such a $V$, the vector field
$$
F(x) = -\nabla V(x)
$$
becomes infinitely repulsive as $x$ approaches $K$.  If the
vector $v=\nabla V(0)$ is not zero then we can
replace $V$ with
$$
\tilde V(x) = V(x) - \sum_{k=1}^n v_kx_k
$$
in order to achieve the normalization $\nabla V(0) = 0$.

Let $P(t)$ be the set of all continuously differentiable
functions
$$
f:[0,t]\to \Bbb R^n\setminus K
$$
satisfying $f(0) = f^\prime(0) = 0$.  We can use the vector
field $F$ to construct a concatenation operation on
the family $P = \{(t,f): t>0, f\in P(t)\}$ as follows.
Starting with any continuous function $\phi:[0,t]\to\Bbb R^n$
with $\phi(0)=0$, there is a corresponding element
$f\in P(t)$ defined as the solution of the initial
value problem
$$
\align
f^\prime(s) &= \phi(s) -\nabla V(f(s)),\\
f(0) &= 0.
\endalign
$$
Notice that $f^\prime(0)=0$ because of the normalization
of $V$.  There is no problem with the existence of global
solutions of this differential equation because of the
nature of $V$.  Indeed, it is quite easy to show that
this differential equation has a unique solution over the
interval $0\leq s<\infty$, and we obtain $f$ by
restricting the solution to the interval $[0,t]$.
Conversely, every $f\in P(t)$ arises in this
way from the ``driving" function
$\phi(s) = f^\prime(s) + \nabla V(f(s))$,
$0\leq s\leq t$.  Thus $f\leftrightarrow \phi$ defines
a bijection of families of sets over $(0,\infty)$.

Thus we can transfer the concatenation of example 1.5 directly
to define a concatenation
$f,g\in P(s)\times P(t) \mapsto f*g\in P(s+t)$.  Explicitly,
given $f\in P(s)$, $g\in P(t)$, then $h = f*g$ is the element
of $P(s+t)$ defined by the initial value problem
$$
\align
h^\prime(\lambda) =& \phi(\lambda) - \nabla V(h(\lambda))\\
h(0) =&0
\endalign
$$
on the interval $0\leq \lambda\leq s+t$, where
$\phi:[0,s+t]\to \Bbb R^n$ is the driving function
$$
\phi(\lambda) =
\cases
f^\prime(\lambda),& 0\leq \lambda<s\\
f^\prime(s) + g^\prime(\lambda-s)+\nabla V(f(s))+\nabla V(g(\lambda-s)),&
s\leq \lambda\leq s+t.
\endcases
$$
Notice that $f*g$ agrees with $f$ on
the interval $0\leq \lambda\leq s$, but fails to agree
with $f(s) + g(\lambda-s)$ on the interval
$s\leq \lambda\leq s+t$.
\endexample

\remark{Remark}There are many potentials that one could use, and thus
there are infinitely many path space structures whose
fiber spaces are of the form
$$
P(t) = \{f:[0,t]\to\Bbb R^n\setminus K : f\in C^1,\quad f(0)=f^\prime(0) = 0\},
$$
$t>0$.  It is clear from the construction that
all of these path spaces are
isomorphic to the example 1.5.  On the other hand, notice
that if one is presented with just the path space structure of one
of these examples, it will not be possible to write down the relevant
isomorphism if one does not know the correct potential
$V$ (at least up to a constant).
\endremark

\subheading{2.  Additive forms and multiplicative forms}

Given a set $X$ and a positive definite function
$B:X\times X\to \Bbb C$, one can construct a Hilbert space
$H(X,B)$.  If $B$ and $C$ are two positive definite
functions
$$
\align
B:&X\times X\to \Bbb C\\
C:&Y\times Y\to \Bbb C
\endalign
$$
then the product $B\times C$ defines a positive definite
function on the cartesian product $(X\times Y)\times(X\times Y)$
by way of
$$
B\times C((x_1,y_1),(x_2,y_2)) = B(x_1,x_2)C(y_1,y_2),
$$
and one has a natural isomorphism of Hilbert spaces
$$
H(X\times Y,B\times C) = H(X,B)\otimes H(Y,C).
$$
Thus, if one is given a path space $\pi:P\to (0,\infty)$ and a
family of positive definite functions
$$
B_t:P(t)\times P(t) \to \Bbb C,\qquad t>0
$$
which is multiplicative in an appropriate sense, then
one would expect to obtain a continuous tensor product of
Hilbert spaces...more precisely, a {\it product system}
in the sense of \cite{2}.

The purpose of this section is to formulate these issues
in the case where $B_t$ is {\it infinitely divisible} in the
sense that there is a family of conditionally positive
definite functions
$$
g_t: P(t)\times P(t) \to \Bbb C\qquad t>0
$$
such that $B_t = e^{g_t}$.  Of course, multiplicative properties
of $B$ follow from appropriate additive properties of $g$.  A
classification of all product systems that can arise from
this construction appears in Corollary 4.33.

We begin by reviewing a few basic facts and some terminology
relating to conditionally positive definite
functions \cite{18}.  Let $X$ be a set and let
$g: X\times X\to \Bbb C$ be a function.  $g$ is called
conditionally positive definite if it is self-adjoint in
the sense that $g(y,x) = \overline{g(x,y)}$ for all
$x,y$, and is such that for every $x_1,x_2,\dots,x_n\in X$ and
$\lambda_1,\lambda_2,\dots,\lambda_n\in \Bbb C$ with
$\lambda_1+\lambda_2+\dots+\lambda_n=0$, we have
$$
\sum_{i,j=1}^n \lambda_i \overline\lambda_jg(x_i,x_j)\geq 0.
$$
If $p:X\times X\to \Bbb C$ is a positive definite function
and $\psi:X\to \Bbb C$ is arbitrary then
$$
g(x,y) = p(x,y) + \psi(x) + \overline{\psi(y)} \tag{2.1}
$$
is conditionally positive definite.  For $g$ and $p$ fixed, the
function $\psi$ satisfying (2.1) is not unique, but
it is unique up to a perturbation of the form
$$
\psi^\prime(x) = \psi(x) + ic
$$
where $c$ is a real constant.  Conversely, every conditionally
positive definite function $g$ can be decomposed into a sum
of the form (2.1).  For fixed $g$, the positive definite
function appearing in (2.1) is not unique.
If $g:X\times X\to\Bbb C$ is conditionally positive definite
then
$$
B(x,y) = e^{g(x,y)}
$$
defines a positive definite function $B: X\times X\to \Bbb C$.
The converse is false: there are self-adjoint functions
$g:X\times X\to \Bbb C$ which are not conditionally positive
definite whose exponentials $B = e^g$ are positive definite.
However, if there is a sequence $\lambda_n$ of positive numbers
which tends to zero as $n\to \infty$ and
$$
B_n(x,y) = e^{\lambda_n g(x,y)}\tag{2.2}
$$
is positive definite for every $n=1,2,\dots$, then $g$ is
conditionally positive definite.  Notice that if
$$
B_n(x,y) = e^{\frac{1}{n}g(x,y)}
$$
is positive definite
for every $n=1,2,\dots$ then
$$
B_1(x,y) = e^{g(x,y)} = B_n(x,y)^n
$$
has a positive definite $n^{\text th}$ root for every
$n=1,2,\dots$; such a positive definite function $B_1$
is called {\it infinitely divisible}.

Let $P$ be a path space, which will be fixed throughout
this section.  $P^2$ will denote the fiber product
$$
P^2 = \{(t,x,y): x,y\in P(t)\},
$$
with projection $\pi(t,x,y) = t$ and fiber spaces
$P^2(t) = P(t)\times P(t)$.  A function $g:P^2\to \Bbb C$
is called conditionally positive definite if for every
$t>0$, the restriction of $g$ to $P(t)\times P(t)$
is conditionally positive definite.

Given such a function $g: P^2\to \Bbb C$, we may construct a
Hilbert space $H(t)$ for every $t>0$.  Briefly, letting
$\Bbb C_0P(t)$ denote the set of all functions
$f:P(t)\to \Bbb C$ such that $f(x) = 0$ for all but finitely
many $x\in P(t)$ and which satisfy $\sum_x f(x) = 0$,
then $\Bbb C_0P(t)$ is a complex vector space and we may
define a positive semidefinite sesquilinear form $\<\cdot,\cdot\>$
on $\Bbb C_0P(t)$ by
$$
\<f_1,f_2\> = \sum_{x,y\in P(t)}f_1(x)\overline{f_2(y)} g(x,y).
$$
After passing to the quotient of $\Bbb C_0P(t)$ by the
subspace $\{f\in \Bbb C_0P(t): \<f,f\> = 0\}$ and completing
the resulting inner product space, we obtain a Hilbert space
$H(t)$.  We will say that $g$ is {\it separable} if $H(t)$ is
a separable Hilbert space for every $t>0$.

In spite of the fact that each fiber $P(t)$ is a lifeless
set without additional structure, there is a useful notion
of measurability for conditionally positive definite functions
$g:P^2\to \Bbb C$.  Fix $0<s<t$.  For every element
$y\in P(t)$ we may consider its associated propagator
$\{y(\lambda,\mu): 0\leq\lambda<\mu\leq t\}$.  Notice that for
every $x\in P(s)$ and every $\lambda$ in the interval
$[0,t-s]$ we can form the complex number $g(x,y(\lambda,\lambda+s)$.

\proclaim{Definition 2.2}$g$ is called measurable if
for every $0<s<t<\infty$, every pair of elements
$x_1,x_2\in P(s)$ and every $y\in P(t)$,
$$
\lambda\in (0,t-s)\mapsto g(x_1,y(\lambda,\lambda+s))
- g(x_2,y(\lambda,\lambda+s))
$$
defines a complex-valued Borel function.
\endproclaim
\noindent
Indeed, we will see that for the conditionally positive
definite functions that are of primary interest, the functions
appearing in Definition 2.2 are actually {\it continuous} (see
Theorem 4.3).  In view of the fact that we have imposed
no structure on path spaces beyond that which follows from their
rule of multiplication, this property appears noteworthy.

Finally, we introduce an appropriate notion of
additivity for conditionally positive definite
functions $g: P^2\to \Bbb C$.

\proclaim{Definition 2.3}
$g$ is called additive if, there is a function defined on
the full cartesian product
$$
\psi: P\times P\to \Bbb C
$$
such that for all $s,t>0$, $x_1,x_2\in P(s)$, $y_1,y_2\in P(t)$,
$$
g(x_1y_1,x_2y_2)-g(x_1,x_2)-g(y_1,y_2) =
\psi(x_1,y_1)+\overline{\psi(x_2,y_2)}.\tag{2.4}
$$
\endproclaim

\remark{Remarks} Notice that the domain of $\psi$, namely
$P\times P$, is larger than the domain of $g$, namely
$P^2$.  In fact, $P^2$ is the diagonal of $P\times P$:
$$
P^2 = \{(p,q)\in P\times P: \pi(p)=\pi(q)\},
$$
$\pi: P\to (o,\infty)$ being the natural projection.

The function $\psi$ of Definition 2.3 is called the
{\it defect} of $g$.  The defect of $g$ is not uniquely
determined by equation (2.4), but if $\psi_1$ and
$\psi_2$ both satisfy (2.4) then it is easy to see that
there must be a function
$c:(0,\infty)\times (0,\infty)\to\Bbb R$ such that
$\psi_2(x,y)=\psi_1(x,y)+ic(s,t)$ for every $x\in P(s)$,
$y\in P(t)$, $s,t>0$.
$g$ is called {\it exact} if there
is a function $\rho:P\to \Bbb C$ and a real-valued
function $c:(0,\infty)\times (0,\infty)\to \Bbb R$ such
that for every $s,t>0$, $x\in P(s)$, $y\in P(t)$,
$$
\psi(x,y) = \rho(xy) - \rho(x) -\rho(y) + ic(s,t).  \tag{2.5}
$$
Notice that when $g$ is exact we can replace it with
$$
g_0(x,y) = g(x,y) - \rho(x) -\overline{\rho(y)}
$$
to obtain a new conditionally positive definite
function $g_0:P^2\to \Bbb C$ which obeys the addition
formula (2.4) with zero defect.
\endremark

\proclaim{Definition 2.6}
Let P be a path space.  An additive form is a
function $g: P^2\to \Bbb C$ which restricts to a conditionally
positive definite function on $P(t)\times P(t)$ for every
$t>0$, and which is separable, measurable, and additive.

A metric path space is a pair $(P,g)$ consisting of
a path space $P$ and an additive form $g: P^2\to \Bbb C$.
\endproclaim

There are many natural examples of additive forms on path spaces.
We give two simple ones here that are important for probability
theory, and a third example
which will be central to what follows.

For every $t>0$ let $PC[0,t]$ denote the space of all
piecewise continuous real-valued functions $f: [0,t]\to \Bbb R$
and let $PC$ be the path space
$$
\align
PC&=\{(t,f): t>0,\quad f\in PC[0,t]\,\}\\
\pi(t,f)&=t
\endalign
$$
with concatenation defined by
$(s,f)(t,g) = (s+t,f*g)$ where
$$
f*g(\lambda) =
\cases
f(\lambda),& 0\leq\lambda<s\\
g(\lambda-s),& s\leq\lambda\leq s+t.
\endcases
$$
Let $c$ be a positive constant.
\example{Example 2.7: Gaussian forms on $PC$}
For $x_1,x_2\in PC[0,t]$, put
$$
g(x_1,x_2) = -c\int_0^t[x_1(\lambda)-x_2(\lambda)]^2\, d\lambda.
$$
\endexample

\example{Example 2.8: Poisson forms on $PC$}
Let $h$ be a second positive constant and
for $x_1$, $x_2\in PC[0,t]$ put
$$
g(x_1,x_2) = c\int_0^t [e^{ih(x_1(\lambda) - x_2(\lambda))} -1]\,d\lambda.
$$
\endexample

The forms $g$ defined in examples 2.7 and 2.8 are
essentially the covariance
functions associated with random processes of the type indicated
by their name.  Note that in both cases the processes have
stationary increments, and in fact 2.7 is the covariance function of
Brownian motion.  Notice too that the detailed structure of the
path space used in these examples is not critical.  For instance,
if one replaces
$PC[0,t]$ with the corresponding Skorohod space $D[0,t]$ and
imitates what was done above, the new examples will
share the essential features as those of 2.7 and 2.8.

More generally, with any continuous conditionally positive
definite function of two real
variables $\gamma:\Bbb R^2\to \Bbb C$
we can associate an additive form $g$ on $PC$
by way of
$$
g(x_1,x_2) = \int_0^t\gamma(x_1(\lambda),x_2(\lambda))\,d\lambda,
$$
for $x_1,x_2\in PC[0,t]$.

Notice too that in all of the preceding examples the defect
of $g$ is zero.  This will {\it not} be the case for additive
forms that are associated with decomposable product systems as in
in chapter II, and we will have to deal with additive forms
having nonzero defects.

\example{2.9: The standard examples}
Let $\Cal C$ be a separable Hilbert space and consider the
path space $\Cal P_{\Cal C}$ of example 1.7.  In this case
the additive form $g: \Cal P_{\Cal C}^2\to \Bbb C$ is simply
the inner product inherited from $L^2((0,\infty);\Cal C)$,
$$
g(f_1,f_2) = \int_0^t \<f_1(\lambda),f_2(\lambda)\>\,d\lambda,
$$
for $f_1$, $f_2\in \Cal P_{\Cal C}(t)$.

If we replace $\Cal C$ with a Hilbert space $\Cal C^\prime$
having the same dimension $n$ as $\Cal C$ then we obtain a new
path space $\Cal P_{\Cal C^\prime}$ and a new form
$$
g^\prime : \Cal P_{\Cal C^\prime}^2 \to \Bbb C.
$$
However, any unitary operator $W: \Cal C\to \Cal C^\prime$ induces
an obvious isomorphism of path space structures
$$
\tilde W: \Cal P_{\Cal C}\to \Cal P_{\Cal C^\prime}
$$
by way of $\tilde W f(\lambda)= Wf(\lambda)$ for $\lambda\in(0,t]$,
$t>0$.  It follows that
$$
g^\prime(\tilde Wf_1,\tilde Wf_2) = g(f_1,f_2)
$$
for $f_1$, $f_2\in \Cal P_{\Cal C}(t)$, $t>0$.  We conclude that,
up to isomorphism, the examples of 2.9 depend only on the
dimension $n$.  This sequence of metric path spaces will be
denoted $(\Cal P_n, g_n)$, $n=1,2,\dots,\infty$.  It is
convenient to include the degenerate pair $(\Cal P_0,g_0)$ where
$\Cal P_0$ is the trivial path space $\Cal P_0 = (0,\infty)\times \{0\}$
with multiplication $(s,0)(t,0) = (s+t,0)$
and additive form $g_0(x,y) = 0$ for all $x,y$.  Thus
$(\Cal P_n,g_n)$ is defined for every $n=0,1,2,\dots,\infty$.
\endexample

Suppose now that we are given a metric path space $(P,g)$, and
assume for the moment that $g$ has {\it defect zero}.  Then for
every $t>0$ we can define a positive definite function on $P(t)$
by
$$
x,y\in P(t) \mapsto e^{g(x,y)}.
$$
Let $E(t)$ be the Hilvert space obtained from this positive definite
function.  In more detail, there is a function $F_t: P(t)\to E(t)$ with the
property that $E(t)$ is spanned by the range of $F_t$ and
$$
\<F_t(x),F_t(y)\> = e^{g(x,y)}
$$
for every $x,y\in P(t)$.  It can be seen that the separability
hypothesis on $g$ implies that $E(t)$ is separable (in fact, the
separability hypothesis implies that $E(t)$ can be identified
with a subspace of the symmetric Fock space over a separable
one-particle space as at the end of section 4).

Let us examine the consequences of the formula (2.4) with
$\psi=0$.  Fixing $x_1,x_2\in P(s)$ and $y_1,y_2\in P(t)$,
we have
$$
e^{g(x_1y_1,x_2y_2)} = e^{g(x_1,x_2)} e^{g(y_1,y_2)}
$$
and hence
$$
\<F_{s+t}(x_1y_1),F_{s+t}(x_2y_2)\> =
\<F_{s}(x_1),F_{s}(x_2)\>\<F_{t}(y_1),F_{t}(y_2)\>.
$$
It follows that there is a unique bilinear map
$(\xi,\eta)\in E(s)\times E(t)\mapsto \xi\cdot\eta\in E(s+t)$
which satisfies $F_s(x)\cdot F_t(y) = F_{s+t}(xy)$ for all
$x\in P(s)$, $y\in P(t)$ and this extended mapping acts like
tensoring in that
$$
\<\xi_1\cdot\eta_1,\xi_2\cdot\eta_2\> =
\<\xi_1,\xi_2\>\<\eta_1,\eta_2\>
$$
for every $\xi_1,\xi_2\in E(s)$, $\eta_1,\eta_2\in E(t)$.

Thus we can define an associative operation in the total family
of Hilbert spaces $p:E\to (0,\infty)$ defined by
$$
\align
E &= \{(t,\xi): t>0,\quad \xi\in E(t)\, \}\\
p(t,\xi)&= t
\endalign
$$
by way of $(s,\xi)(t,\eta) = (s+t,\xi\cdot\eta)$.  This structure
$p:E\to (0,\infty)$ has the main features of a product system
\cite{2}.  However, since the total family $E$ carries no
natural Borel structure (because we are not given a Borel structure
on the total space of $P$), the measurability axioms for
product systems are meaningless here.  For this reason, we
will refer to $p: E\to (0,\infty)$ as the {\it product structure}
associated with the metric path space $(P,g)$.

The above construction required that $g$ have defect zero.  If
$g$ has nonzero defect but is exact, then one can modify this
construction so as to obtain a product structure in this case
as well (see section 4).

In general, Theorem 4.3 below implies that every additive
form on a path space is exact.  Moreover, we will find that
$(P,g)$ must be ``essentially"
isomorphic to one of the standard examples $(\Cal P_n,g_n)$,
$n=0, 1, 2,\dots, \infty$ in such a way that the
product structure associated with $(P,g)$
is either the trivial one having
one-dimensional spaces $E(t)$ or is isomorphic to one of the
standard product systems of \cite{2}.
Thus, {\it every\/} metric path space
gives rise to a product structure
that is completely understood.

\subheading{3.  Exactness of cocycles}
Let $\Cal C$ be a separable Hilbert space and let $L$
be the space of all measurable functions
$f:(0,\infty)\to \Cal C$ which are locally square integrable
in the sense that
$$
\int_0^T \|f(x)\|^2\, dx < \infty
$$
for every $T>0$.  The topology of $L$ is defined by the
sequence of seminorms
$$
\|f\|_n = (\int_0^n \|f(x)\|^2\, dx)^{1/2}
$$
$n=1,2,\dots$ and
$$
d(f,g) = \sum_1^\infty 2^{-n}\frac{\|f-g\|_n}{1+\|f-g\|_n}
$$
is a translation invariant metric on $L$ with respect to which
it becomes a separable Fr\`echet space.

$L$ is the dual of the inductive limit of Hilbert spaces
$$
L_0 = \bigcup_{T>0} L^2((0,t);\Cal C).
$$
$L_0$ is identified with the submanifold of $L^2((0,\infty);\Cal C)$
consisting of all functions having compact support; a net
$f_\alpha\in L_0$ converges to $f\in L_0$ iff there is a
$T>0$ such that $f_\alpha$ is supported in $(0, T)$ for
sufficiently large $\alpha$, and
$$
\lim_{\alpha}\int_0^T \|f_\alpha(x)-f(x)\|^2\, dx = 0.
$$
$L$ is isometrically anti-isomorphic to the dual of $L_0$
by way of the sesquilinear pairing
$$
f,g\in L_0\times L\mapsto \<f,g\> = \int_0^\infty \<f(x),g(x)\>\, dx.
$$
A function $t\in (0,\infty)\mapsto \phi_t\in L$ is called
{\it measurable} if it is a Borel function relative to the weak$^*$
topology on $L$, i.e.,
$$
t\in (0,\infty)\mapsto \<f,\phi_t\>
$$
should be a complex-valued Borel function for every $f\in L_0$.

\proclaim{Definition 3.1}
An additive cocycle is a measurable function
$t\in (0, \infty)\mapsto \phi_t\in L$ satisfying
$$
\phi_{s+t}(x) = \phi_s(x) + \phi_t(x+s) \qquad {\text a.e. }(dx)
$$
for every $s,t>0$.
\endproclaim

The purpose of this section is to prove the following
characterization of additive cocycles.

\proclaim{Theorem 3.2}
Let $\{\phi_t: t>0\}$ be an additive cocycle.  Then there
is a function $f\in L$ such that
$$
\phi_t(x) = f(x+t) - f(x) \qquad {\text a.e. }(dx)
$$
for every $t>0$.
\endproclaim

\remark{Remarks}
Theorem 3.2 bears a resemblance to known results about
multiplicative cocycles associated with transitive actions
of topological groups on topological spaces
\cite{20}.  However, our setting here differs
in several key aspects.  Rather than a group action
we have a transitive action of the additive semigroup
of positive reals.  Moreover, the elements of
$L$ must satisfy nontrivial integrability conditions.

Notice too that we make no assumption about the
boundedness of $\{\phi_t: t>0\}$.  That is to say,
we do {\it not} assume that
$$
\sup_{t>0} |\<f,\phi_t\>| \leq M_f < \infty
$$
for every $f\in L_0$.  Indeed, it is easy to
give examples of cocycles for which this
condition is not satisfied.  The fact that additive
cocycles can be unbounded means that softer techniques
that require some form of boundedness (specifically,
techniques involving the use of Banach limits on the
additive semigroup $\Bbb R^+$) are not available.
For these reasons we have taken some care to give
a complete proof of Theorem 3.2.
\endremark

\remark{Remark on measurability}
Every weak$^*$ measurable function $t\mapsto \phi_t\in L$
is a Borel mapping of metric spaces.
This means that for every $g\in L$,
the function
$$
t\in (0,\infty)\mapsto d(\phi_t,g)
$$
is Borel measurable.  To see that, notice that for each
$n=1,2,\dots$,
$$
\|\phi_t-g\|_n = \sup_{k\geq 1}|\<f_k, \phi_t-g\>|,
$$
the supremum extended over  a countable set of functions
$f_1,f_2,\dots$
which is dense in the unit ball of $L^2((0,n);\Cal C)$.  Hence
$t\mapsto \|\phi_t-g\|_n$ is a Borel function for each $n$.
It follows that
$$
d(\phi_t,g) =
\sum_{n=1}^\infty 2^{-n}\frac{\|\phi_t-g\|_n}{1+\|\phi_t-g\|_n}
$$
is a Borel function of $t$.
\endremark

The following lemma asserts that additive cocycles are
metrically continuous.  I am indebted to Calvin Moore for pointing
out a pertinent reference \cite{12, Th\'eor\`eme 4, p. 23} for
a key step in its proof.

\proclaim{Lemma 3.3}
Let $\{\phi_t: t>o\}$ be an additive cocycle.  Then
$t\mapsto \phi_t$ moves continuously in the
metric topology of $L$ and we have
$$
\lim_{t\to 0+}\phi_t = 0.
$$
\endproclaim

\demo{proof}
Let $\{T_t: t\geq 0\}$ be the natural translation semigroup
which acts in $L$ by
$$
T_tf(x) = f(x+t).
$$
For every $n=1,2,\dots$ and every $f\in L$ we have
$$
\lim_{t\to 0}\|T_tf - f\|_n^2 =
\lim_{t\to 0}\int_0^n \|f(x+t) - f(x)\|^2\, dx = 0,
$$
and hence $d(T_tf,f) \to 0$ as $t\to 0$.  Thus the
semigroup $\{T_t: t\geq 0\}$ is continuous in the metric
topology of $L$.

By the preceding remark, the function
$$
t\in (0,\infty) \mapsto \phi_t\in L
$$
is a Borel function taking values in a separable complete
metric space.  Thus there is a subset $N\subseteq (0,\infty)$
of the first category in $(0,\infty)$ such that the restriction
of this function to $(0,\infty)\setminus N$ is metrically
continuous \cite{19, p. 306}.

Let $t_0\in (0,\infty)$ and choose $t_n\in (0,\infty)$ such
that $t_n\to t_0$ as $n\to \infty$.  Since
$$
M = \cup_{n=0}^\infty N-t_n
$$
is a first category subset of $(0,\infty)$ we may find
$s\in (0,\infty)\setminus M$.  Then $s+t_n$ belongs to
$(0,\infty)\setminus N$ for every $n=0,1,2,\dots$, and it
follows that
$$
\lim_{n\to \infty}d(\phi_{s+t_n},\phi_{s+t_0}) = 0.
$$
Writing
$$
\phi_{t_0} + T_{t_0}\phi_s = \phi_{t_0+s} = \lim_{n\to \infty}\phi_{t_n+s}
=\lim_{n\to \infty}(\phi_{t_n} + T_{t_n}\phi_s)
$$
and noting that by the metric continuity of $\{T_t\}$ we have
$$
\lim_{n\to \infty}T_{t_n}\phi_s = T_{t_0}\phi_s,
$$
we conclude that
$$
\phi_{t_0}=\lim_{n\to \infty}\phi_{t_n}.
$$

For continuity at $t=0+$, notice that for every $s,t_0>0$ we have
$$
\phi_{s+t_0} = \phi_s + T_s\phi_{t_0}.
$$
Letting $s\to 0+$ and noting that
$\phi_{s+t_0} \to \phi_{s+t_0}$ and $T_s\phi_{t_0}\to\phi_{t_0}$,
we obtain
$$
\lim_{s\to 0+}\phi_s = 0,
$$
as required \qed
\enddemo

\demo{proof of Theorem 3.2}
By the Fubini theorem we may find a Borel function
$$
v:(t,x)\in (0,\infty)\times(0,\infty)\mapsto v(t,x)\in\Cal C
$$ such that
$v(t,x) = \phi_t(x)$
almost everywhere with respect to the product measure $dt\,dx$
on $(0,\infty)\times(0,\infty)$.  Because of the cocycle condition
$$
\phi_{s+t}(x) = \phi_s(x) + \phi_t(x+s)\qquad \text{a.e. }(dx)
$$
it follows that $v$ must satisfy
$$
v(s+t,x) = v(s,x) + v(t,x+s)
$$
for almost every triple $(s,t,x) \in (0,\infty)^3$ with respect
to $ds\,dt\,dx$.  By the Fubini theorem there is a Borel set
$N\subseteq (0,\infty)$ of measure zero such that for every
$x\in (0,\infty)\setminus N$ we have
$$
v(s+t,x) = v(s,x) + v(t,x+s)
$$
almost everywhere $ds\,dt$.  Choose a sequence
$x_1>x_2>\dots$ in $(0,\infty)\setminus N$ which decreases to $0$.

We claim that there is a sequence of Borel functions
$f_n:(x_n,\infty)\to \Cal C$ satisfying
$$
f_n(\xi+t) - f_n(\xi) = v(t,\xi) \tag{3.5.1}
$$
almost everywhere on $(x_n,\infty)\times (0,\infty)$ with respect
to $d\xi\, dt$, and which
is coherent in that
$$
f_{n+1}\restriction_{(x_n,\infty)} =
f_n\qquad \text{almost everywhere}.  \tag{3.5.2}
$$
Proceeding inductively, we define $f_1$ by
$$
f_1(\xi) = v(\xi-x_1,x_1),\qquad \xi>x_1.
$$
To verify (3.5.1), note that for $\lambda,t>0$ and
$\xi = \lambda+x_1$ we have
$$
f_1(\lambda+x_1+t) - f_1(\lambda+x_1) = v(\lambda+t,x_1) - v(\lambda,x_1)
$$
which by (3.4) is almost everywhere ($d\lambda\, dt$) equal to
$v(t,x_1+\lambda)$, hence we have (3.5.1).

Assuming now that $f_k$ has been defined for $1\leq k\leq n$,
define $g:(x_{n+1},\infty)\to \Cal C$ by
$$
g(\xi) = v(\xi -x_{n+1}, x_{n+1}).
$$
If we replace $x_1$ with $x_{n+1}$ in the argument of the preceding
paragraph we obtain the conclusion
$$
g(\xi+t) - g(\xi) = v(t,\xi)
$$
almost everywhere ($d\xi\, dt$) on $(x_{n+1},\infty)\times (0,\infty)$.
Restricting $\xi$ to the interval $(x_n,\infty)\subseteq(x_{n+1},\infty)$
then gives
$$
g(\xi+t)-g(\xi) = v(t,\xi) = f_n(\xi+t) - f_n(\xi)
$$
almost everywhere ($d\xi\, dt$) on $(x_n,\infty)\times (0,\infty)$.
Thus the function $h:(x_n,\infty)\to \Cal C$ defined by
$h(x) = g(x) - f_n(x)$ is Borel-measurable and translation invariant
in the sense that for almost every $t\geq 0$ we have
$$
h(\xi+t) = h(\xi)\qquad \text{a.e. }(d\xi).
$$
It follows that there is a vector $c\in \Cal C$ such that
$$
h(\xi) = c\qquad \text{a.e. }(d\xi).
$$
If we set $f_{n+1} = g(x) - c$ then we obtain both
required conditions (3.5.1) and (3.5.2).

Because of the coherence property (3.5.2), there is a Borel
function $f:(0,\infty)\to \Cal C$ satisfying
$$
f(\xi) = f_n(\xi) \qquad \text{a.e. }(d\xi)
$$
on $(x_n,\infty)$, for every $n=1,2,\dots$; and because
of (3.5.1) we have
$$
f(\xi+t) - f(\xi) = v(t,\xi)
$$
for almost every pair $(\xi,t)\in (0,\infty)\times(0,\infty))$
with respect to the product measure $d\xi\, dx$.

It follows that
$$
f(\xi+t) -f(\xi) = \phi_t(\xi)
$$
almost everywhere $d\xi\,dt$, and by another application of
Fubini's theorem we may conclude that there is a Borel set
$N\subseteq (0,\infty)$ such that for all
$t\in (0,\infty)\setminus N$ we have
$$
\phi_t(x) = f(x+t) - f(x) \qquad \text{a.e. }(dx). \tag{3.6}
$$

We show next that the exceptional set $N$ of (3.6) can be eliminated.
To that end, consider the vector space $\Cal F$
consisting of all Borel functions
$$
F:(0,\infty)\to \Cal C,
$$
where we make the traditional identification of two functions
that agree almost everywhere.  We endow this space with the
topology of local convergence in measure.  More precisely,
a net $f_\alpha\in \Cal F$ converges to $F\in \Cal F$ iff for
every $T>0$, the restrictions $F_\alpha\restriction_{(0,T)}$
converge in Lebesgue measure to $F\restriction_{(0,T)}$.
$\Cal F$ is metrizable as a separable complete metric space
by way of
$$
d(F,G) = \sum_{n=1}^\infty 2^{-n} d_n(F,G)
$$
where for $n=1,2,\dots$,
$$
d_n(F,G) = \int_0^n \frac{\|F(x) - G(x)\|}{1+\|F(x) - G(x)\|}\,dx .
$$
The translation semigroup $\{T_t: t\geq 0\}$ defined by
$$
T_tF(x) = F(x+t)
$$
acts continuously on the Fr\`echet space $\Cal F$.

It follows from these remarks that the right side of (3.6) is
continuous in $t$, provided that we consider $T_tf-f$ as an
element of $\Cal F$.  Since the inclusion map
$L\subseteq \Cal F$ is continuous, it follows from Lemma 3.3
that the left side of (3.6) defines a continuous function
$$
t\in (0,\infty)\mapsto \phi_t\in \Cal F.
$$
Equation (3.6) implies that these two continuous functions
agree on the complement of a null set, and hence they agree
for all $t>0$.

Thus, we have a Borel function $f:(0,\infty)\to L$ with the
property that for every $t>0$,
$$
\phi_t(x) = f(x+t) - f(x) \qquad \text{a.e. }(dx).  \tag{3.7}
$$
It remains to show that the condition $\phi_t\in L$ for
every $t>0$ implies that $f$ itself belongs to $L$.  We will
deduce that from the following result, the proof of which is
based on an argument shown to me by Henry Helson, who kindly
consented to its inclusion in this paper.

\proclaim{Lemma 3.8}
Let $f$ be a nonnegative Borel function defined on $(0,\infty)$
satisfying
$$
\int_0^T|f(x+t)-f(x)|^2\, dx < \infty
$$
for every $T>0$ and every $t>0$.  Then for every $T>0$ we have
$$
\int_0^T |f(x)|^2\, dx < \infty.
$$
\endproclaim
\demo{proof}
Fix $T>0$.  It suffices to show that for every nonnegative
funtion $g\in L^2(0,T)$ we have
$$
\int_0^T f(x)g(x)\, dx < \infty.
$$

To that end, find a function $u\in L^2(0,1)$ such that
$u(x)>0$ for every $x$ and
$$
\int_0^1f(x)u(x)\,dx < \infty,
$$
and define a function $\Phi: (0,\infty)\to \Bbb R^+$ by
$$
\Phi(t) = \int_0^1 |f(x+t)-f(x)|u(x)\, dx.
$$
We claim that $\Phi$ is continuous and tends to $0$ as
$t\to 0+$.  Indeed, we may apply Lemma 3.3 to the additive
cocycle
$$
\phi_t(x) = f(x+t) - f(x)
$$
in the case where $\Cal C$ is the complex numbers to conclude
that $t\in (0,\infty)\mapsto \phi_t\in L$ is metrically continuous
and tends to $0$ as $t\to 0+$.  Since $F\in L\mapsto |F|\in L$
is clearly continuous, the same is true of the modulus
$t\in (0,\infty)\mapsto |\phi_t|\in L$, and the claim follows.

Notice next that for every positive $g\in L^2$,
$$
\int_0^1\int_0^T f(x+t)g(t)u(x)\, dt\,dx < \infty.  \tag{3.10}
$$
Indeed, since
$$
f(x+t)g(t)u(x) \leq |\phi_t(x)|g(t)u(x) +f(x)g(t)u(x),
$$
and since
$$
\int_0^1\int_0^T |\phi_t(x)|g(t)u(x)\, dt\, dx
=\int_0^T\Phi(t)g(t)\,dt < \infty
$$
and
$$
\int_0^1\int_0^Tf(x)g(t)u(x)\, dt\, dx =
\int_0^Tf(x)u(x)\,dx \int_0^1 g(t)\,dt < \infty,
$$
(3.10) follows.

Since $u>0$, (3.10) implies that
$$
\int_0^T f(x+t)g(t)\, dt <\infty
$$
for almost every $x\in (0,1)$.  Nor for every $x>0$ we have
$$
\int_0^T f(t)g(t) \, dt \leq
\int _0^T f(x+t)g(t)\, dt + \int_0^T |f(x+t) - f(x)|g(t)\, dt.
$$
The first integral on the right is finite for certain values of
$x\in (0,1)$ and the second one is finite for all $x > 0$.  Thus
(3.9) follows.  \qed
\enddemo

To complete the proof of theorem 3.2, we must show that any
measurable function $f:(0,\infty)\to \Cal C$ for which the
differences $\phi_t(x) = f(x+t)-f(x)$ belong to $L$ must
itself belong to $L$.  Equivalently, if $f:(0,\infty)\to \Cal C$
is a Borel function for which
$$
\int_0^T\|f(x+t)-f(x)\|^2\, dx < \infty
$$
for every $T>0$ and every $t>0$, then
$$
\int_0^T\ \|f(x)\|^2\, dx < \infty \tag{3.11}
$$
for every $T>0$.

To see that, fix such an $f$ and consider $F(x) = \|f(x)\|$.
We have
$$
|F(x+t) - F(x)| \leq \|f(x+t)-f(x)\|,
$$
hence $F$ satisfies the hypotheses of Lemma 3.8.  It follows
from Lemma 3.8 that $F$ is locally in $L^2$.  \qed
\enddemo

\subheading{4.  Classification of additive forms}
Before stating the main result on classification
of metric path spaces we introduce the concept of
a strongly spanning set.  Let $H$ be a Hilbert space
and let $S$ be a subset of $H$.  We will write
$e^H$ for the symmetric Fock space over the one-particle
space $H$,
$$
e^H = \sumplus_{n= 0}^\infty H^{(n)},
$$
$H^{(n)}$ denoting the symmetric tensor product
of $n$ copies of $H$ when $n\geq 1$, and where
$H^{(0)}$ is defined as $\Bbb C$.
Consider the exponential map $\exp:H\to e^H$,
defined by
$$
\exp(\xi)=\sum_{n=0}^\infty\frac{1}{\sqrt{n!}}\,\xi^n,
$$
$\xi^n$ denoting $\xi^{\otimes n}$ if $n\geq 1$ and
$\xi^0 = 1\in \Bbb C$.  It is well known
that $e^H$ is spanned by
$$
\exp(H) = \{\exp(\xi): \xi\in H\}.
$$

\proclaim{Definition 4.1}
A subset $S\subseteq H$ is said to strongly span $H$
if $e^H$ is spanned by the set of vectors
$\exp(S) = \{exp(\xi): \xi\in S\}$.
\endproclaim

\remark{Remarks}  When we use the term `span', we of course
mean {\it closed} linear span.
Every strongly spanning set must span $H$, but the
converse if false, as the following remarks show.
In general, every vector $\zeta\in e^H$ give rise
to a holomorphic function $f_\zeta:H\to \Bbb C$ by way of
$$
f_\zeta(\xi)=\<\exp(\xi),\zeta\>.
$$
Indeed, if we let $\zeta_n$ be the projection of $\zeta$
onto $H^{(n)}$, then we have a representation of
$f_\zeta$ as a power series
$$
f_\zeta(\xi) = \sum_0^\infty \frac{1}{\sqrt{n!}}\<\xi^n,\zeta_n\>,
$$
and because $\sum_n\|\zeta_n\|^2 = \|\zeta\|^2 < \infty$,
this power series converges absolutely and uniformly
over the ball in $H$ of radius $R$, for every $R>0$.
Thus, $\Cal F=\{f_\zeta: \zeta\in e^H\}$ is a complex
vector space of entire functions defined on $H$.  $S$ is
a strongly spanning set iff
$$
F\in \Cal F,\quad F(S) = \{0\}\implies F=0.
$$

If $H$ is finite dimensional then every holomorphic
polynomial $F:H\to \Bbb C$ belongs to $\Cal F$.  In
particular, if $H$ is two dimensional and $\{e_1,e_2\}$
is an orthonormal basis for $H$ then
$$
S = \{\lambda(e_1+e_2): \lambda\in\Bbb C\}\cup
\{\lambda(e_1-e_2): \lambda\in\Bbb C\}
$$
clearly spans $H$ because it contains $e_1+e_2$ and
$e_1-e_2$.  On the other hand, $S$ is not strongly
spanning because
$$
F(\xi) = \<\xi,e_1\>^2 - \<\xi,e_2\>^2
$$
is a nonzero polynomial which vanishes on $S$.  More
generally, any spanning subset $S$ of a finite dimensional
$H$ which is contained in the zero set of a nontrivial
polynomial will fail to be strongly spanning.  This will
be the case whenever $S$ is an algebraic set, or an algebraic
variety in $H$.
\endremark

The following is our main classification of metric path
spaces.

\proclaim{Theorem 4.3}
Let $P$ be a path space and let $g:P^2\to\Bbb C$ be a
separable measurable additive form.  Then there is a separable
Hilbert space $\Cal C$, a complex-valued function
$\rho:P\to \Bbb C$ and a mapping of fiber spaces
$$
\log:P\to \Cal P_{\Cal C}
$$
such that $\log(xy) = \log(x)\boxplus\log(y)$ for every
$x,y\in P$, and such that for every $t>0$, $x_1,x_2\in P(t)$,
$$
\align
\log(P(t))& \text{ strongly spans }\Cal P_{\Cal C}(t)\tag{4.3.1}\\
g(x_1,x_2)& =
\<\log(x_1),\log(x_2)\> + \rho(x_1) + \overline{\rho(x_2)}\tag{4.3.2}
\endalign
$$
\endproclaim

\remark{Remarks}
The assertion that $\log$ is a mapping of fiber spaces means
$\log(P(t))\subseteq\Cal P_{\Cal C}(t)$ for every $t>0$.  Thus,
$\log$ defines a homomorphism of the path space structure of
$P$ into that of $\Cal P_{\Cal C}$.  Property 4.3.1 asserts that,
even though $\log(P(t))$ may not be dense in $\Cal P_{\Cal C}(t)$,
it is a rich enough subset so that
$$
\overline{\text{span}}\exp(P(t)) = e^{\Cal P_{\Cal C}(t)}.
$$

Finally, notice that property 4.3.2 implies that $g$ is an
{\it exact} form in the sense of (2.5).  Indeed, if we take
$x_1,x_2\in P(s)$ and $y_1,y_2\in P(t)$ then by definition
of the concatenation operation $\boxplus$ in $\Cal P_{\Cal C}$
we have
$$
\<\log(x_1)\boxplus\log(y_1),\log(x_2)\boxplus\log(y_2)\> =
\<\log(x_1),\log(x_2\> + \<\log(y_1),\log(y_2\>.
$$
Using the fact that $\log(x_ky_k) = \log(x_k)\boxplus\log(y_k)$
for $k=1,2$ we find that
$$
g(x_1y_1,x_2y_2)-g(x_1,x_2)-g(y_1,y_2) =
\psi(x_1,y_1) + \overline{\psi(x_2,y_2)}
$$
where $\psi(x,y) = \rho(xy)-\rho(x)-\rho(y)$ and $\rho$ is
the function given by 4.3.2.
\endremark

The proof of Theorem 4.3 will occupy most of the remainder of
this section, and will proceed along the following lines.
We first use $g$ to associate a Hilbert space $H_t$ with
$P(t)$ for every $t>0$.  We then show that for $s<t$, $H_s$
embeds naturally in $H_t$ so that we can form an inductive
limit of Hilbert spaces
$$
H_\infty = \lim_{\rightarrow}H_t.
$$
We introduce a strongly continuous semigroup of isometries
acting in $H_\infty$ which will turn out to be {\it pure}.
This implies that $H_\infty$ can be coordinatized in such a way
that it becomes an $L^2$ space of vector valued functions
$$
H_\infty = L^2((0,\infty); \Cal C)
$$
and the semigroup of isometries is the natural shift semigroup.
Finally, we use the results of section 3 to solve a
cohomological problem.  Once that is accomplished we can
define the required ``logarithm" $\log:P\to \Cal P_\Cal C$
and verify its properties.

\example{Definition of $H_t$}
Fix $t>0$.  Let $\Bbb C_0P(t)$ denote the complex vector
space of all finitely nonzero functions $f:P(t)\to \Bbb C$
satisfying the condition $\sum_xf(x) = 0$, and let $\<\cdot,\cdot\>$
be the sesquilinear form defined on $\Bbb C_0P(t)$ by
$$
\<f,g\> = \sum_{x,y\in P(t)}f(x)\overline{g(y)}g(x,y).
$$
$\<\cdot,\cdot\>$ is positive semidefinite, and after passing to
the quotient of $\Bbb C_0P(t)$ by the subspace of null functions
$K_t=\{f: \<f,f\> = 0\}$ we obtain an inner product space, whose
completion is denoted $H_t$.

Now $\Bbb C_0P(t)$ is spanned by the set of all differences
$\{\delta_x-\delta_y: x,y\in P(t)\}$, $\delta_z$ denoting the
unit function
$$
\delta_z(u) =
\cases
1,& \text{if }u=z\\
0,& \text{otherwise}.
\endcases
$$
Hence $H_t$ is spanned by the set $\{[x]-[y]: x,y\in P(t)\}$,
where $[x]-[y]$ denotes the element $\delta_x-\delta_y+K_t$.
The inner product in $H_t$ is characterized by
$$
\<[x_1]-[y_1],[x_2]-[y_2]\> =
g(x_1,x_2)-g(x_1,y_2)-g(y_1,x_2) + g(y_1,y_2), \tag{4.4}
$$
for $x_1,x_2,y_1,y_2\in P(t)$.  Notice that although we
have written $[x]-[y]$ as if it were a difference, it is
not actually the difference it appears to be
since $[x]$ and $[y]$ do not belong
to $H_t$.  It is in fact a two-variable function which
satisfies a certain cocycle identity. But the notation is
convenient provided one is careful never to treat $[x]$ and
$[y]$ as if they were elements of $H_t$.

$H_t$ is separable because
of our separability hypothesis on $g$ (See the discussion preceding
Definition 2.2).  When it is necessary to distinguish between
the various inner products we will write $\<\cdot,\cdot\>_t$ for
the inner product on $H_t$.
\endexample

\example{Embedding $H_s$ in $H_t$ for $s<t$}
Fix $s,t$ with $0<s<t$ and choose an element $e\in P(t-s)$.
We want to show that there is an isometric
linear map of $H_s$ into $H_t$ which carries differences
of the form $[x_1]-[x_2]$ with $x_i\in P(s)$
into $[x_1e]-[x_2e]$, and moreover
that this isometry does not depend on the
particular choice of $e\in P(t-s)$. To that end,
we claim that for all $x_i,y_i\in P(s)$, and
$z_i\in P(t-s)$ $i=1,2$, we have
$$
\<[x_1e]-[x_2e],[y_1z_1]-[y_2z_2]\>_t =
\<[x_1]-[x_2],[y_1]-[y_2]\>_s\tag{4.5}
$$
Indeed, because of (4.4) the left side is
$$
g(x_1e,y_1z_1)-g(x_1e,y_2z_2)-g(x_2e,y_1z_1) + g(x_2e,y_2z_2)\tag{4.6}
$$
But by (2.4) we have,
$$
g(x_ie,y_jz_j) = g(x_i,y_j) + g(e,z_j) + \psi(x_i,e)+\overline{\psi(y_j,z_j)}.
$$
It follows that for $j=1,2$,
$$
g(x_1e,y_jz_j) - g(x_2e,y_jz_j) =
g(x_1,y_j) - g(x_2,y_j) + \psi(x_1,e) - \psi(x_2,e).
$$
When we subtract this expression for $j=2$ from
the expression for $j=1$ the terms involving $\psi$ cancel
and we are left with
$$
g(x_1,y_1)-g(x_1,y_2) -g(x_2,y_1) + g(x_2,y_2)
$$
which is the right side of (4.5).

By taking $z_1=z_2=e$ in (4.5) we see that
$$
\<[x_1e]-[x_2e],[y_1e]-[y_2e]\>_t =
\<[x_1]-[x_2],[y_1]-[y_2]\>_s
$$
and hence there is a unique linear isometry
$U(t,s): H_s\to H_t$ such that
$$
U(t,s) :[x_1]-[x_2]\mapsto [x_1e]-[x_2e]
$$
for every $x_1,x_2\in P(s)$.  Moreover, since $H_t$
is spanned by elements of the form $[y_1z_1]-[y_2z_2]$
for $y_i\in P(s)$, $z_i\in P(t-s)$ (here we use the
fact that $P(t) = P(s)P(t-s)$), it also follows from
(4.5) that $U(t,s)$ is independent of the particular
choice of $e\in P(t-s)$; in more concrete terms, if
$e$ and $f$ are two elements of $P(t-s)$,
then for all $x_1,x_2\in P(s)$ we have
$$
[x_1e]-[x_2e] = [x_1f]-[x_2f].
$$

It follows from the latter that we have the consistency
relation
$$
U(t,s)U(s,r) = U(t,r)
$$
for all $0<r<s<t$.  Indeed, if we choose $e_1\in P(r-s)$ and
$e_2\in P(t-s)$ then for $x_1,x_2\in P(r)$ we have
$$
U(t,s)U(s,r)([x_1]-[x_2])= U(t,s)([x_1e_1]-[x_2e_1])=
[x_1e_1e_2]-[x_2e_1e_2],
$$
and the right side is $U(t,r)([x_1]-[x_2])$ simply because
$e_1e_2$ is an element of $P(t-s)$.
\endexample

Thus we can form the inductive limit of inner product
spaces
$$
\lim_{\rightarrow}H_t.
$$
Explicitly, this consists of all functions
$t\in (0,\infty)\mapsto \xi_t\in H_t$ having the property
that there is a $T = T_\xi>0$ such that for all $t>T$,
$$
\xi_t=U(t,T)\xi_T.
$$
The inner product in the inductive limit is defined
by
$$
\<\xi,\eta\>=\lim_{t\to\infty}\<\xi_t,\eta_t\>_t.
$$
$H_\infty$ is defined as the completion of
$$\lim_{\rightarrow}H_t,$$
and it is a separable Hilbert space.

Choose $x_1,x_2\in P(t)$.  By a slight abuse of notation we
will write $[x_1]-[x_2]$ for the element of $H_\infty$ defined
by the function
$$
\xi_\lambda =
\cases
0,& \text{for }\lambda\leq t\\
U(\lambda,t)([x_1]-[x_2]),&\text{for } \lambda>t.
\endcases
$$
Notice that {\it by definition} of $H_\infty$, we will
have
$$
[x_1e]-[x_2e] = [x_1]-[x_2]
$$
for every $x_1,x_2\in P(t)$, $t>0$ and for an arbitrary
element $e$ of $P$.  $H_\infty$ is spanned by the set of
formal differences
$$
\{[x_1]-[x_2]: x_i\in P(t), t>0\}.
$$
Finally, note that the inner product in $H_\infty$ is
defined by its values on these formal differences as
follows.  Choose $s\neq t$, $x_1,x_2\in P(s)$ and
$y_1,y_2\in P(t)$.  In order to evaluate the inner
product $\<[x_1]-[x_2],[y_1]-[y_2]\>$ we may suppose
that $s<t$.  Choose any $e\in P(t-s)$.  Then since
$[x_1]-[x_2] = [x_1e]-[x_2e]$ we have
$$
\align
\<[x_1]-[x_2],[y_1]-[y_2]\> &= \<[x_1e]-[x_2e],[y_1]-[y_2]\>\\
&=g(x_1e,y_1) - g(x_1e,y_2) - g(x_2e,y_1) +g(x_2e,y_2).
\endalign
$$

\example{The subspaces $N_t\subseteq H_\infty$}
For every $t>0$ we define a subspace $N_t$ of $H_\infty$
as follows
$$
N_t=\overline{\text{span}}\{[x_1]-[x_2]:x_i\in P(t)\}.
$$
Choose $0<s<t$ and $x_1,x_2\in P(s)$.
The preceding remarks imply that the element of $H_\infty$
represented by the difference $[x_1]-[x_2]$
can be identified with a difference $[y_1]-[y_2]$ of elements
from $P(t)$ by taking $y_i=x_ie$ for some $e\in P(t-s)$.  It
follows that the spaces $N_t$ are increasing
$$
0\leq s\leq t\implies N_s\subseteq N_t \tag{4.7.1}.
$$
Moreover, since the images of the spaces $H_t,\ t>0$ in
$H_\infty$ span $H_\infty$ we also have
$$
\overline{\cup_{t>0}N_t} = H_\infty.   \tag{4.7.2}
$$
\endexample

\example{The semigroup $\{U_t:t\geq 0\}$}
We now introduce a semigroup of isometries $\{U_t: t\geq 0\}$
acting on $H_\infty$.  Fix $t>0$ and choose an element
$e\in P(t)$.

There is a formula analogous to (4.5) in which
the order of multiplication is reversed.
That is, for $0<s<t$, $x_i,y_i\in P(s)$ and $e,z_i\in P(t-s)$
we claim that
$$
\<[ex_1]-[ex_2],[z_1y_1]-[z_2y_2]\>_t =
\<[x_1]-[x_2],[y_1]-[y_2]\>_s  \tag{4.8}
$$
Notice that the inner product on the left (resp. right)
is taken in the Hilbert space $H_t$ (resp. $H_s$).  The
proof of (4.8) is the same as the proof of (4.5).  The identity
(4.8) implies that if we choose an element $f\in P(t)$ for
some $t>0$ then for any $s>0$ and any
pair of elements $x_1,x_2\in P(s)$, the element of $H_\infty$
defined by $[fx_1]-[fx_2]$ does not depend on the particular
choice of $f\in P(t)$ in that for every $g\in P(t)$ we have
$$
[fx_1]-[fx_2] = [gx_1]-[gx_2]. \tag{4.9}
$$
Note that in (4.9) the vectors on
both sides belong to $H_\infty$.
Moreover, if $y_1,y_2$ is a second pair of elements
in $P(s)$ then (4.8) also implies
$$
\<[fx_1]-[fx_2],[fy_1]-[fy_2]\> = \<[x_1]-[x_2],[y_1]-[y_2]\>.
$$

Now fix $t>0$ and choose $f\in P(t)$.  It follows that there
is a unique isometry $U_{t,s}:N_s\to N_{s+t}$ which satisfies
$$
U_{t,s}([x_1]-[x_2]) = [fx_1]-[fx_2].
$$
Because of (4.9), $U_{t,s}$ does not depend on the choice
of $f$.  Notice too that $U_{t,s}$ does not depend on $s$.
Indeed, if $0<s_1<s_2$ and $x_1,x_2\in P(s_1)$ then
for any element $z\in P(s_2-s_1)$ we have
$$
U_{t,s_1}([x_1]-[x_2]) =[fx_1]-[fx_2] = [(fx_1)z]-[(fx_2)z],
$$
while
$$
U_{t,s_2}([x_1]-[x_2])=U_{t,s_2}([x_1z]-[x_2z])
=[f(x_1z)]-[f(x_2z)].
$$
The right hand sides of these two formulas agree because
of the associativity of the multiplication in $P$.

By (4.7.2) there is a unique isometry $U_t:H_\infty\to H_\infty$
satisfying
$$
U_t([x_1]-[x_2]) = [fx_1]-[fx_2]
$$
for all $x_1,x_2\in P(s)$, $s>0$.
For $t=0$ we set $U_0={\bold 1}$.  Finally, note that
$\{U_t: t\geq 0\}$ is a semigroup.  Indeed, given $s,t>0$
we choose $f\in P(s)$ and $g\in P(t)$, and note that for
every $x_1,x_2\in P(r)$ we have
$$
U_sU_t([x_1]-[x_2]) = U_s([gx_1]-[gx_2]) = [fgx_1]-[fgx_2].
$$
The right side must be $U_{s+t}([x_1]-[x_2])$ because
the product $fg$ belongs to $P(s+t)$.  This shows that
$U_sU_t = U_{s+t}$ on each $N_r$, and by (4.7.2) it follows
that $U_sU_t = U_{s+t}$.
\endexample

\example{Strong continuity}
Since $H_\infty$ is a separable Hilbert space, the strong
continuity of $\{U_t: t\geq 0\}$ will follow is we prove
that the operator function
$$
t\in(0,\infty)\mapsto U_t\in \Cal B(H_\infty)
$$
is weakly measurable in the sense that $\<U_\lambda\xi,\eta\>$
defines a Borel function on $0<\lambda<\infty$ for every
$\xi,\eta\in H_\infty$ \cite{2, Proposition 2.5 (ii)}.
In turn, because of the semigroup property it suffices to verify
this for $\lambda$ restricted to the interval $0<\lambda\leq 1$.
To that end, we claim that for any pair of vectors
$\xi_1,\xi_2$ in the spanning set
$$
\bigcup_{t>0}\{[x_1]-[x_2]: x_i\in P(t)\},
$$
the function
$$
\lambda\in (0,1]\mapsto \<U_\lambda\xi_1, \xi_2\> \tag{4.10}
$$
is Borel measurable.

To see this suppose that $\xi=[x_1]-[x_2]$ with $x_i\in P(s)$
and $\xi_2=[y_1]-[y_2]$ with $y_i\in P(t)$.  By replacing
$y_1,y_2$ with $y_1v, y_2v$ for an appropriate $v\in P$ we may assume
that $t$ is as large as we please, and in particular we may
assume that $t>s+1$.  Choosing elements $f\in P(t-s-\lambda)$
and $e\in P(\lambda)$, we have
$$
U_\lambda([x_1]-[x_2]) = U_\lambda([x_1f]-[x_2f])
=[ex_1f]-[ex_2f],
$$
and thus
$$
\<U_\lambda\xi_1,\xi_2\>=\<[ex_1f]-[ex_2f],[y_1]-[y_2]\>
=\alpha_{11}-\alpha_{12}-\alpha_{21}+\alpha_{22},\tag{4.11}
$$
where
$$
\alpha_{ij} = g(ex_if,y_j).
$$
In order to calculate the terms $\alpha_{ij}$ we make use
of the propagators of $y_2$ and $y_2$ to obtain the
factorizations
$$
y_j=Y_j(0,\lambda)y_j(\lambda,\lambda+t)
=y_j(0,\lambda)y_j(\lambda,\lambda+s)y_j(\lambda+s,t).
$$
Thus
$$
\align
\alpha_{ij}=& g(ex_if,y_j(0,\lambda)y_j(\lambda,\lambda+t)) \\
=&g(e,y_j(0,\lambda))+g(x_if,y_j(\lambda,\lambda+t))
+\psi(e,x_if)+\overline{\psi(y_j(0,\lambda),y_j(\lambda,\lambda+t))}\\
=&g(e,y_j(0,\lambda))+g(x_i,y_j(\lambda,\lambda+s))
+g(f,y_j(\lambda+s,t))\ +\\
&\psi(x_i,f)+\overline{\psi(y_j(\lambda,\lambda+s),y_j(\lambda+s,t))}
+\psi(e,x_if)+\overline{\psi(y_j(0,\lambda),y_j(\lambda,\lambda+t))}.
\endalign
$$
Noting that $\alpha_{ij}$ has the form
$$
\alpha_{ij} = g(x_i,y_j(\lambda,\lambda+s))+u_i+v_j
$$
for appropriate complex numbers $u_1,u_2,v_1,v_2$ (which
depend on $\lambda$), it follows that
the $u$'s and $v$'s cancel out of the right side of (4.11) and
we are left with
$$
\align
\alpha_{11}-\alpha_{12}-\alpha_{21}+\alpha_{22} =
&g(x_1,y_1(\lambda,\lambda+s))-g(x_2,y_1(\lambda,\lambda+s))+\\
&g(x_2,y_2(\lambda,\lambda+s))-g(x_1,y_2(\lambda,\lambda+s)).
\endalign
$$
Since $g$ is a measurable form, each of the two functions
$$
\lambda\in(0,1]\mapsto
g(x_1,y_j(\lambda,\lambda+s))-g(x_2,y_j(\lambda,\lambda+s))
$$
$j=1,2$ is a Borel function, and thus the right side of the
previous formula is a difference of Borel functions.
\endexample

\example{Purity of $\{U_t: t\geq0\}$}
We claim next that the semigroup $\{U_t: t\geq0\}$ is
{\it pure} in the sense that
$$
\bigcap_{t>0}U_tH_\infty = \{0\}.  \tag{4.12}
$$
This is a consequence of (4.7.1), (4.7.2) and the following.
\proclaim{Lemma 4.13}
For every $t>0$, $H_\infty$ decomposes into a direct sum
$$
H_\infty = N_t\oplus U_tH_\infty.
$$
\endproclaim
\demo{proof}
Fix $t>0$.  We show first that $N_t$ is orthogonal to
$U_tH_\infty$.  For that, it suffices to show that for every
$r>0$ and for $x_1,x_2\in P(t)$, $y_1,y_2\in P(r)$ we have
$$
\<[x_1]-[x_2],U_t([y_1]-[y_2])\> = 0.  \tag{4.14}
$$
Choose elements $e\in P(t)$, $f\in P(r)$.  Then we have
$$
\align
[x_1]-[x_2]&= [x_1f]-[x_2f],\qquad\text{and}\\
U_t([y_1]-[y_2])&=[ey_1]-[ey_2].
\endalign
$$
Thus the left side of (4.14) has the form
$$
\<[x_1f]-[x_2f],[ey_1]-[ey_2]\> =
\alpha_{11}-\alpha_{12}-\alpha_{21}+\alpha_{22},
$$
where
$$
\alpha_{ij} = g(x_if,ey_j).
$$
Using the definition of additive forms (2.3) we have
$$
\alpha_{ij} =g(x_i,e)+g(f,y_j) +\psi(x_i,f)+\overline{\psi(e,y_j)}
=u_i+v_j
$$
where
$$
\align
u_i=&g(x_i,e)+\psi(x_i,f)\\
v_j=&g(f,y_j)+\overline{\psi(e,y_j)}.
\endalign
$$
It follows that all of the $u$'s and $v$'s cancel and we
are left with the required formula
$$
\alpha_{11}-\alpha_{12}-\alpha_{21}+\alpha_{22}=0.
$$

To show that $N_t\cup U_tH_\infty$ spans $H_\infty$ it is enough
to show that for every $r>t$ and every pair $x_1,x_2\in P(r)$,
we have
$$
[x_1]-[x_2]\in N_t+U_tH_\infty.
$$
If we factor $x_i=a_ib_i$ where $a_i=x_i(0,t)\in P(t)$ and
$b_i=x_i(t,r)\in P(r-t)$, then we have
$$
\align
[x_1]-[x_2] = [a_1b_1]-[a_2b_2] =
&([a_1b_1]-[a_2b_1])+([a_2b_1]-[a_2b_2]) \\
=&([a_1]-[a_2])+U_t([b_1]-[b_2])
\endalign
$$
and the right side clearly belongs to $N_t+U_tH_\infty$ \qed
\enddemo
\endexample

\example{$\{U_t: t\geq0\}$ as a shift}
A familiar theorem asserts that every strongly continuous
pure semigroup of isometries is unitarily equivalent to a
direct sum of copies of the semigroup of simple unilateral
shifts acting on $L^2(0,\infty)$.  From this it follows that
we can replace $H_\infty$ with the Hilbert space
$L^2((0,\infty);\Cal C)$ of all square integrable vector
valued measurable functions $\xi: (0,\infty)\to \Cal C$
with inner product
$$
\<\xi,\eta\>=\int_0^\infty \<\xi(x),\eta(x)\>\, dx
$$
in such a way that $\{U_t: t\geq 0\}$ becomes the semigroup
$$
U_t\xi(x) =
\cases
\xi(x-t),& x>t\\
0,& 0<x\leq t.
\endcases
$$
After making this identification, we find that the
range of $U_t$ consists of all functions
$\xi\in L^2((0,\infty); \Cal C)$ which vanish almost
everywhere in the interval $0<x\leq t$.  From Lemma 4.13
we conclude that for every $t>0$,
$$
\overline{\text{span}}\{[x_1]-[x_2]:x\in P(t)\}=L^2((0,t);\Cal C),
\tag{4.14}
$$
the right side denoting the subspace of all functions
$\xi\in L^2((0,\infty);\Cal C)$ which vanish almost everywhere
outside the interval $0<x\leq t$.  Since $H_\infty$ is
a separable Hilbert space, it follows that $\Cal C$ must
be separable as well.

Finally, we remind the reader that the rules for left
and right multiplication in these new ``coordinates" are
the same as they were in $H_\infty$: if $x_1,x_2\in P(s)$
and $u$ is any element of $P(t)$, then
$$
\align
[x_1u]-[x_2u] &= [x_1]-[x_2],\qquad\text{and}\\
[ux_1]-[ux_2] &= U_t([x_1]-[x_2]).
\endalign
$$

Equation (4.14) identifies the space
$N_t=\overline{\text{span}}\{[x_1]-[x_2]:x_i\in P(t)\}$
with the space $\Cal P_{\Cal C}(t)$ of example 2.9 for every
$t>0$.  It remains to define the ``logarithm" mapping
$\log: P\to \Cal P_{\Cal C}$ with the properties
asserted in Theorem 4.3.
\endexample

\example{Definition of the logarithm}
In order to define the logarithm we must first show that certain
$2$-cocycles are trivial.  These cocycles are associated with
globally defined left-coherent sections
$$
t\in (0,\infty)\mapsto e_t\in P(t)
$$
whose existence is established in the following.

\proclaim{Lemma 4.15}
Choose any element $e\in P(1)$.  Then there is a
left-coherent family of elements $\{e_t: t>0\}$
with the property that $e_1 = e$.
\endproclaim
\demo{proof}
For $0<t\leq 1$ we set $e_t=e(0,t)$ where
$\{e(s,t):0\leq s<t\leq 1\}$ is the propagator associated
with $e$ as in section 1.  For $n<t\leq n+1$ we set
$$
e_t=e_1^ne(0,t-n).
$$
It is clear that $\{e_t: t>0\}$ has the required properties\qed
\enddemo

Choose such a section $\{e_t:t>0\}$, which will be fixed
throughout the remainder of this section.  Define a function
$\Gamma:(0,\infty)\times(0,\infty)\to L^2((0,\infty);\Cal C)$
by
$$
\Gamma(s,t) = [e_se_t]-[e_{s+t}].
$$
We will see presently that $\Gamma$ is an additive
$2$-cocycle in the sense that for all $r,s,t>0$ we have
$$
\Gamma(r+s,t)-\Gamma(r,s+t) -U_r\Gamma(s,t)+\Gamma(r,s) = 0,
$$
see (4.20).  The following asserts that $\Gamma$ is exact
in the sense that we require.

\proclaim{Theorem 4.16}
There is a measurable function
$t\in (0,\infty)\mapsto \phi_t\in L^2((0,\infty);\Cal C)$
such that
$$
\align
\phi_t(x)=&0 \qquad\text{a.e. outside } 0<x\leq t,\tag{4.16.1}\\
\Gamma(s,t) =& \phi_{s+t}-\phi_s -U_s\phi_t,\qquad s,t>0. \tag{4.16.2}
\endalign
$$
\endproclaim

\remark{Remark 4.17}
The measurability assertion of Theorem 4.16 simply means that for every
$\xi\in L^2((0,\infty);\Cal C)$,
$$
t\in (0,\infty)\mapsto \<\phi_t,\xi\>
$$
is a complex-valued Borel function.  Because $L^2((0,\infty);\Cal C)$
is separable, this is equivalent to measurability of
$t\mapsto \phi_t$ relative to the metric topology of
$L^2((0,\infty);\Cal C)$.
Notice too that (4.16.1) asserts that $\phi_t\in \Cal P_{\Cal C}(t)$
for every $t>0$.

Assume, for the moment, that Theorem 4.16 has been proved.
We can then define a fiber map
$$
\log :P\to \Cal P_{\Cal C}
$$
as follows.  For $z\in P(t)$, $t>0$ we put
$$
\log(z)=[z]-[e_t]-\phi_t.
$$
Notice that $\log(z)\in \Cal P_{\Cal C}(t)$ because
both $[z]-[e_t]$ and $\phi_t$ belong to $\Cal P_{\Cal C}(t)$.
For $x,y\in P$ we claim:
$$
\log(xy)=\log(x)\boxplus\log(y).  \tag{4.18}
$$
This is to say that, if $x\in P(s)$ and $y\in P(t)$ then
$$
\log(xy) = \log(x) + U_s(\log(y));
$$
equivalently,
$$
[xy]-[e_{s+t}] -\phi_{s+t}=[x]-[e_s]-\phi_s+U_s([y]-[e_t]-\phi_t).
$$
To see that this is the case, note that
$$
\align
[xy]-[e_{s+t}]=&([xy]-[e_sy]) +([e_sy]-[e_se_t])+([e_se_t]-[e_{s+t}])\\
=&[x]-[e_s] +U_s([y]-[e_t])+\Gamma(s,t).
\endalign
$$
Using (4.16.2) to substitute for $\Gamma(s,t)$ and subtracting
$\phi_{s+t}$ from both sides, we obtain (4.18).
\endremark

\demo{proof of Theorem 4.16}
The argument will proceed as follows.  We first find a family
$\{u_t: t>0\}$ of Borel functions
$$
u_t:(0,\infty)\to \Cal C
$$
which satisfy
$$
\Gamma(s,t) = u_{s+t}-u_s-U_s(u_t)
$$
almost everywhere, for every $s,t>0$.  The family
$\{u_t: t>0\}$ will {\it not} satisfy (4.16.1), but
these functions will be locally in $L^2$ in the sense that
$$
\int_0^T \|u_t(x)\|^2\, dx < \infty
$$
for every $T>0$.  We will then use the results of
section 3 to find a locally $L^2$ function $w:(0,\infty)\to \Cal C$
with the property that for every $t>0$
$$
u_t(x+t) = w(x+t) - w(x)
$$
almost everywhere $(dx)$ on the interval $(0,\infty)$.  Such
a function $w$ can be subtracted from $u_t$
$$
\phi_t(x) =
\cases
u_t(x)-w(x),& 0<x\leq t\\
0,& x>t
\endcases
$$
so that the modification $\{\phi_t: t>0\}$ satisfies
both (4.16.1) and (4.16.2).

\proclaim{Lemma 4.19}
For each $s,t>0$, $\Gamma(s,t)(x)$ vanishes almost everywhere
($dx$) outside the interval $0<x<s+t$.  $\Gamma$ is a $2$-cocycle
in the sense that for every $r,s,t>0$ we have
$$
\Gamma(r+s,t) - \Gamma(r,s+t) -
U_r\Gamma(s,t) + \Gamma(r,s) = 0.\tag{4.20}
$$
\endproclaim
\demo{proof}
Since both $e_se_t$ and $e_{s+t}$ belong to $P(s+t)$,
$\Gamma(s,t) = [e_se_t]-[e_{s+t}]$ belongs to
$$
\overline{\text{span}}\{[x_1]-[x_2]:x_i\in P(s+t)\}=L^2((0,s+t);\Cal C)\}.
$$
Moreover, since $e_{s+t} = e_se(s,s+t)$, we see that
$$
\Gamma(s,t) = [e_se_t]-[e_se(s,s+t)] =U_s([e_t]-[e(s,s+t)])
$$
belongs to the range of $U_s$ and hence $\Gamma(s,t)$ must
vanish almost everywhere on the interval $0<x\leq s$.

To prove (4.20), notice that
$$
U_r\Gamma(s,t) = U_r([e_se_t]-[e_{s+t}]) = [e_re_se_t]-[e_re_{s+t}],
$$
and hence by the definition of $\Gamma$ we have
$$
\align
\Gamma(r+s&,t) - \Gamma(r,s+t) - U_r\Gamma(s,t) \\
&=[e_{r+s}e_t]-[e_{r+s+t}]-[e_re_{s+t}]+[e_{r+s+t}]
-[e_re_se_t]+[e_re_{s+t}] \\
&=[e_{r+s}e_t]-[e_re_se_t] = [e_{r+s}]-[e_re_s]=-\Gamma(r,s),
\endalign
$$
as required\qed
\enddemo

\proclaim{Lemma 4.21}
For fixed $t>0$ and $\xi\in L^2((0,\infty);\Cal C)$,
the function
$$
s\mapsto \<\Gamma(s,t),\xi\>
$$
is Borel-measurable.
\endproclaim
\demo{proof}
Fix $t$.  Noting that
$$
\align
\Gamma(s,t) &= [e_se_t]-[e_{s+t}] =[e_se_t]-[e_se(s,s+t)]\\
&=U_s([e_t]-[e(s,s+t)]),
\endalign
$$
we have
$$
\<\Gamma(s,t),\xi\>=\<[e_t]-[e(s,s+t)],U_s^*\xi\>,
$$
for every $s,t>0$.  Since $s\mapsto U_s\xi$ is (metrically)
continuous, it suffices to show that the function
$s\in (0,\infty)\mapsto [e_t]-[e(s,s+t)]$ is weakly
measurable; i.e., that
$$
s\in (0,\infty)\mapsto\<[e_t]-[e(s,s+t)],\eta\>
$$
is a Borel function for every $\eta\in L^2((0,\infty);\Cal C)$.
Since $L^2((0,T);\Cal C)$ is spanned by $\{[y_1]-[y_2]:y_i\in P(T)\}$
for every $T>0$, this reduces to showing that for fixed $T$ and
$y_1,y_2\in P(T)$,
$$
s\in (0,\infty)\mapsto\<[e_t]-[e(s,s+t)],[y_1]-[y_2]\>
$$
is a Borel function.  To see that, pick $u\in P(T-t)$.  Then
we have $[e_t]-[e(s,s+t)]= [e_tu]-[e(s,s+t)u]$, and if we set
$x_1=e_tu$ and $x_2=e(s,s+t)u$ then $x_i\in P(T)$ and
$$
\align
\<[e_t]-[e(s,s+t)],[y_1]-[y_2]\>=&
\<[x_1]-[x_2],[y_1]-[y_2]\>\\
=&\alpha_{11}-\alpha_{12}-\alpha_{21}+\alpha_{22},
\endalign
$$
where $\alpha_{ij}=g(x_i,y_j)$.  We have
$$
\align
\alpha_{1k} =& g(e_tu,y_k),\\
\alpha_{2k}=&g(e(s,s+t)u,y_k)
\endalign
$$
for $k=1,2$.  By the additivity property 2.3 we have
$$
\align
\alpha_{2k}=&g(e(s,s+t)u,y_k)=g(e(s,s+t)u,y(0,t)y_k(t,T))\\
=&g(e(s,s+t),y_k(0,t))+g(u,y_k(t,T))+\psi(e(s,s+t),u)\\
&+\overline{\psi(y_k(0,t),y_k(t,T))}.
\endalign
$$
Noting that neither $\alpha_{11}$ nor $\alpha_{12}$ involves
$s$ and that the terms $\psi(e(s,s+t),u)$ cancel out
of the difference $\alpha_{22}-\alpha_{21}$, we find that
$$
\alpha_{11}-\alpha_{12}-\alpha_{21}+\alpha_{22}=
-g(e(s,s+t),y_1(0,t))+g(e(s,s+t),y_2(0,t))+K
$$
where $K$ does not depend on $s$.  The right side is
a Borel function of $s$ because of the measurability
hypothesis on $g$\qed
\enddemo

We define a family of functions $u_s:(0,\infty)\to \Cal C$
as follows.

\proclaim{Lemma 4.22}
For every $s>0$, the limit
$$
u_s(\lambda) = -\lim_{n\to\infty}\Gamma(s,n)(\lambda) \tag{4.22.1}
$$
exists almost everywhere on $0<\lambda<\infty$ and satisfies
$$
\int_0^T \|u_s(\lambda)\|^2\, d\lambda <\infty \tag{4.22.2}
$$
for every $T>0$.  $\{u_s: s>0\}$ is measurable in the sense
that for every compactly supported $\xi\in L^2((0,\infty);\Cal C)$,
the function
$$
s\in (0,\infty)\mapsto \<u_s,\xi\>
$$
is Borel measurable.  Putting $u_s(\lambda)=0$ for $\lambda\leq 0$ we have
$$
\Gamma(s,t)(\lambda) =
u_{s+t}(\lambda) -u_s(\lambda)-u_t(\lambda-s) \tag{4.22.3}
$$
almost everywhere on $0<\lambda<\infty$, for every $s,t>0$.
\endproclaim

\remark{Remark}
Actually, the limit in (4.22.1) exists in a very strong
sense.  We will show that as $t$ increases with $s$ fixed,
the restrictions
$$
\Gamma(s,t)\restriction_{(0,T]}
$$
stabilize as soon as $t$ is larger than $T$.  Once one knows
this, the assertion (4.22.2) is an obvious consequence
of the fact that each function $\Gamma(s,t)$ belongs to
$L^2(0,\infty);\Cal C)$.
\endremark

\demo{proof of Lemma 4.22}
We first establish the coherence property described in the
preceding remark.  More precisely, we claim that for fixed
$0<s<T<t_1<t_2$ one has
$$
\Gamma(s,t_2)\restriction_{(0,T]}=\Gamma(s,t_1)\restriction_{(0,T]}.
\tag{4.23}
$$
To see that, consider the difference
$$
\Gamma(s,t_2)-\Gamma(s,t_1)=
[e_se_{t_2}]-[e_{s+t_2}]-[e_se_{t_1}]+[e_{s+t_1}].
$$
Writing
$$
\align
[e_se_{t_2}]-[e_se_{t_1}]&=[e_se_{T-s}e(T-s,t_2)]-[e_se_{T-s}e(T-s,t_1)]\\
&= U_T([e(T-s,t_2)]-[e(T-s,t_1)])
\endalign
$$
and noting that
$$
\align
-[e_{s+t_2}]+[e_{s+t_1}]&=-[e_Te(T,t_2)]+[e_Te(T,t_1)]\\
&=-U_T([e(T,t_2)]-[e(T,t_1)]
\endalign
$$
we find that $\Gamma(s,t_2)-\Gamma(s,t_1)$ has the form
$U_T\zeta$ for $\zeta\in L^2$ given by
$$
\zeta=[e(T-s,t_2)]-[e(T-s,t_1)]- [e(T,t_2)]+[e(T,t_1)].
$$
(4.23) follows because every function in the range of
$U_T$ vanishes a.e. on the interval $(0,T]$.

Thus (4.22.1) follows and by the preceding remark we also
have (4.22.2).  It is also clear that for every compactly
supported function $\xi\in L^2$,
$$
\<u_s,\xi\> = \<\Gamma(s,n),\xi\>
$$
for sufficiently large $n=1,2,\dots$.  Thus the measurability
of $\{u_s\}$ follows as well.

Finally, the formula (4.22.3) follows after restricting
all terms in the cocycle equation (4.16.2) to a finite
unterval $0<\lambda\leq T$ and taking the
formal $\lim_{t\to\infty}$ to obtain
$$
-u_{r+s} + u_r +U_ru_s +\Gamma(r,s+ = 0,
$$
for every $r,s>0$ \qed
\enddemo

We must now modify the family $\{u_s: s>0\}$ in order to
obtain a new family $\phi_s = u_s -w$ which has the additional
property that $\phi_s(s)$ vanishes a.e. outside the interval
$0<\lambda\leq s$.  This is accomplished as follows.

Notice that for $s,t>0$,
$$
u_{s+t}(\lambda) = u_s(\lambda) + u_t(\lambda-s)
$$
almost everywhere on the interval $\lambda\geq s+t$.  Indeed,
this is immediate from the fact that
$$
\Gamma(s,t)(\lambda) = u_{s+t}(\lambda)-u_s(\lambda)-u_t(\lambda-s)
$$
and the fact that $\Gamma(s,t)$  vanishes outside the
interval $0<\lambda\leq s+t$.  Thus if we define
$v_t:(0,\infty)\to \Cal C$ by
$$
v_t(\lambda) = u_t(\lambda+t),
$$
then $\{v_t:t>0\}$ is a measurable
family of $\Cal C$-valued functions satisfying
$$
\int_0^T\|v_t(\lambda)\|^2\, d\lambda <\infty
$$
for every $T,t>0$, for which
$$
v_{s+t}(\lambda)=v_s(\lambda)+v_t(\lambda+s)
$$
almost everywhere ($d\lambda$), for every $s,t>0$.  By Theorem
3.2, there is a Borel function
$$
w:(0,\infty)\to \Cal C
$$
which is locally in $L^2$, such that for every $t>0$ we
have
$$
v_t(\lambda) = w(\lambda+t)-w(\lambda)
$$
almost everywhere on the interval $0<\lambda<\infty$.  Set
$w(\lambda)=0$ for $\lambda\leq 0$.  It follows that
$$
u_t(\lambda) - w(\lambda) +w(\lambda-t)
$$
vanishes almost everywhere on the interval $t<\lambda<\infty$.
Hence
$$
\phi_t(\lambda) = u_t(\lambda) -w(\lambda) +w(\lambda-t)
$$
satisfies both conditions
$$
\phi_{s+t}(\lambda) -\phi_s(\lambda) -\phi_t(\lambda-s) = \Gamma(s,t)(\lambda)
$$
a.e. on $0<\lambda<\infty$ and
$$
\phi_t(\lambda) = 0
$$
almost everywhere on the interval $\lambda>t$.  Notice that we
can also define $\phi_t$ as follows
$$
\phi_t(\lambda) =
\cases
u_t(\lambda)-w(\lambda),& 0<\lambda\leq t\\
0,& \lambda>t.
\endcases
 \tag{4.24}
$$

As in the remarks following the statement of Theorem 4.16,
we can now define a fiber map $\log:P\to \Cal P_{\Cal C}$
by
$$
\log(x) = [x]-[e_t]-\phi_t
$$
for every $x\in P(t)$, and every $t>0$, and this function
satisfies $\log(P(t))\subseteq\Cal P_{\Cal C}(t)$ and
$\log(xy)=\log(x)\boxplus\log(y)$, for $x,y\in P$.
\enddemo
\endexample

It remains to establish (4.3.1), and to exhibit a function
$\rho:P\to \Bbb C$ which satisfies (4.3.2).  $\rho$ is
defined as follows.  If $x\in P(t)$ we put
$$
\rho(x) = \<[x]-[e_t],\phi_t\>+g(x,e_t)-
\frac{1}{2}(g(e_t,e_t)+\|\phi_t\|^2).
$$

To see that (4.3.2) is satisfied we choose $x_1,x_2\in P(t)$
and use the definition of $\log$ to write
$$
\align
\<\log(x_1),\log(x_2)\> =& \<[x_1]-[e_t]-\phi_t,[x_2]-[e_t]-\phi_t\>\\
=&\<[x_1]-[e_t],[x_2]-[e_t]\> -\<[x_1]-[e_t],\phi_t\> \\
&-\<\phi_t,[x_2]-[e_t]\> -\|\phi_t\|^2.
\endalign
$$
Noting that $\<[x_1]-[e_t],[x_2]-[e_t]\>$ expands to
$$
g(x_1,x_2)-g(x_1,e_t)-g(e_t,x_2)+g(e_t,e_t),
$$
we obtain the required formula
$$
\<\log(x_1),\log(x_2)\> = g(x_1,x_2) + \rho(x_1) +\overline{\rho(x_2)}.
$$

It remains to show that $\log(P(t))$ is a strongly spanning set
in $\Cal P_{\Cal C}(t)$.  To see that, fix $t>0$, and let us
write $L = \log(P(t))$.  Notice that $L$ carries no linear
structure {\it a priori}, since the only algebraic property
of the $\log$ function is its additivity
$$
\log(xy) = \log(x)\boxplus\log(y).
$$
Nevertheless, we will show that $L$ is ``almost convex".

For every $r>0$, let
$$
B_r = \{\xi\in \Cal P_{\Cal C}(t): \|\xi\|\leq r\}
$$
be the ball of radius $r$ and let
$$
K=\bigcup_{r>0}\overline{L\cap B_r}^w
$$
where $\overline{L\cap B_r}^w$ denotes the closure of
$L\cap B_r$ in the weak topology of the Hilbert space
$\Cal P_{\Cal C}(t)$.

\proclaim{Lemma 4.26}
If $\xi,\eta\in K$ and $\theta$ is a dyadic rational in the
unit interval, then $\theta\xi+(1-\theta)\eta\in K$.
\endproclaim
\demo{proof}
It clearly suffices to prove that
$$
\xi,\eta\in K\implies \frac{1}{2}(\xi+\eta)\in K.
$$
We will show first that this is true in the special
case where $\xi,\eta$ have the form
$$
\xi=\log(x), \quad \eta=\log(y)
$$
with $x,y\in P(t)$.  To this end we claim that there
is a sequence $z_n\in P(t)$ with the properties
$$
\|\log(z_n)\| \leq \|\log(x)\| + \|\log(y)\|\tag{4.27}
$$
and which, in addition, satisfies
$$
\lim_{n\to\infty}\<\log(z_n),\zeta\> = \frac{1}{2}\<\log(x),\zeta\>
+\frac{1}{2}\<\log(y),\zeta\>,
$$
for all $\zeta\in \Cal P_{\Cal C}(t)$. Indeed, for every
$n=1,2,\dots$, consider a dyadic partition of the
interval $[0,t]$ as follows,
$$
\{0=t_0<t_1<\dots <t_{2^n}=t\}
$$
where $t_k=kt/{2^n}$, $0\leq k\leq 2^n$.  Using the propagators
$\{x(r,s):0\leq r<s\leq t\}$ and $\{y(r,s):0\leq r<s\leq t\}$
for $x$ and $y$ we can define $z_n$ as a product
$$
z_n= x_1y_2x_3y_4\dots x_{2^n-1}y_{2^n}
$$
where $x_k = x(t_{k-1},t_k)$ and $y_k=y(t_{k-1},t_k)$.
Then because of the additivity property of $\log$ we have
$$
\log(z_n) = \log(x_1)\boxplus\log(y_2)\boxplus\log(x_3)\boxplus\log(y_4)
\boxplus\dots\boxplus\log(x_{2^n-1})\boxplus\log(y_{2^n}).
$$
Letting $\Cal O_n$ and $\Cal E_n$ be the respective unions of the
odd and even intervals,
$$
\align
\Cal O_n &= \bigcup_{k \text{ odd}} (t_{k-1},t_k]\\
\Cal E_n &= \bigcup_{k \text{ even}} (t_{k-1},t_k]
\endalign
$$
we can rewrite the previous formula for $\log(z_n)$ as follows
$$
\log(z_n) = \log(x)\chi_{{\Cal O}_n} + \log(y)\chi_{{\Cal E}_n},
$$
$\chi_S$ denoting the characteristic function of the set $S\subseteq [0,t]$.
It follows that
$$
\|\log(z_n)\| \leq \|\log(x)\| + \|\log(y)\|.
$$

Moreover, the arument of \cite{2, p. 47} implies that for
any function $w$ in $L^1[0,t]$ we have
$$
\align
\lim_{n\to \infty}\int_{\Cal O_n}w(x)\,dx &= \frac{1}{2}\int_0^tw(x)\,dx
\qquad\text{and}\\
\lim_{n\to \infty}\int_{\Cal E_n}w(x)\,dx &= \frac{1}{2}\int_0^tw(x)\,dx
\endalign
$$
(An equivalent assertion is that the sequence $\chi_{{\Cal E}_n}
\in L^\infty[0,t]$ converges to the constant function $1/2$ in the
weak$^*$-topology of $L^\infty[0,t]$).  So if we fix $\zeta\in
\Cal P_{\Cal C}(t)=L^2((0,t);\Cal C)$ then we have
$$
\<\log(z_n),\zeta\>=\<\log(x)\cdot\chi_{\Cal O_n}+\log(y)\chi_{\Cal
E_n},\zeta\>
= \int_\Cal{O_n}w_1(\lambda)\,d\lambda + \int_{\Cal E_n}w_2(x)\,dx,
$$
where $w_1(\lambda) =\<\log(x)(\lambda),\zeta(\lambda)\>$,
$w_2(\lambda) = \<\log(y)(\lambda),\zeta(\lambda)\>$.  By the preceding
remarks the right side tends in the limit to
$$
\frac{1}{2}\int_0^t \<\log(x)(\lambda),\zeta(\lambda)\>\,d\lambda
+\frac{1}{2}\int_0^t \<\log(y)(\lambda),\zeta(\lambda)\>\,d\lambda
$$
as asserted.

Now let $\xi, \eta$ be arbitrary elements of $K$.  Since the
sets $\overline{L\cap B_r}^w$ increase with $r$, we can assume
that
$$
\xi,\eta\in \overline{L\cap B_r}^w
$$
for some $r>0$.  It follows that there are sequences $x_n,y_n\in P(t)$
satisfying $\|\log(x_n)\|\leq r$, $\|\log(y_n)\|\leq r$ and
$$
\align
\xi&=\lim_{n\to\infty}\log(x_n),\\
\eta&=\lim_{n\to\infty}\log(y_n)
\endalign
$$
weakly.  Now for each fixed $n=1,2,\dots$
the preceding argument implies that
$$
\frac{1}{2}\log(x_n) + \frac{1}{2}\log(y_n)\in\overline{L\cap B_{2r}}^w.
$$
Since the set on the right is weakly closed, we may take the limit on
$n$ to obtain
$$
\frac{1}{2}\xi + \frac{1}{2}\eta \in \overline{L\cap B_{2r}}^w,
$$
as required\qed
\enddemo

\remark{Remark 4.29}
{}From Lemma 4.26 we immediately deduce: {\it the norm closure of}
$$
K=\bigcup_{r>0}\overline{\log(P(t))\cap B_r}^w
$$
{\it is a convex subset of }$\Cal P_{\Cal C}(t)$.
\endremark

\proclaim{Lemma 4.30}
Set $S$ be a convex subset of a Hilbert space $H$ which
spans $H$.  Then $S$ is a strongly spanning set.
\endproclaim
\demo{proof}
If $S$ contains $0$ then the conclusion follows from
\cite{2 Proposition 6.12}.  Thus we will obtain the
more general result if we show that for every strongly
spanning set $S_0\subseteq H$ and every $\xi\in H$,
$S_0+\xi$ is also a strongly spanning set.

To see that, let $\zeta\in H\mapsto W_\zeta\in \Cal U(\Cal B(e^H))$
be the standard representation of the canonical commutation
relations on the symmetric Fock space
$e^H$; $W_\zeta$ is defined by requiring
$$
W_\zeta:\exp(\eta)\mapsto
e^{-\frac{1}{2}\|\eta\|^2-\<\eta,\zeta\>}\exp(\zeta+\eta),
$$
for every $\eta\in H$.  Each $W_\zeta$ is a unitary operator
on $e^H$.
Now let $v$ be a vector in $e^H$ such that $\<v,\exp(\eta)\>=0$
for every $\eta\in S_0+\xi$.  We have to show that $v=0$.  But
for every $\eta_0\in S_0$ we have
$$
\<v,W_\xi(\exp(\eta_0))\> =
e^{-\frac{1}{2}\|\eta\|^2-\<\zeta,\eta\>}\<v,\exp(\eta_0+\xi)\> = 0
$$
for every $\eta_0\in S_0$, hence $W_\xi^*v$ is orthogonal to
$\exp(S_0)$.  Since $\exp(S_0)$ spans $e^H$ we conclude
that $W_\xi^*v = 0$, hence $v=0$\qed
\enddemo

We can now show that $L=\log(P(t))$ is a strongly spanning
subset of $\Cal P_{\Cal C}(t)$.  To see that,
note first that the exponential
map $\exp:\xi\mapsto \exp(\xi)$ is weakly
continuous on {\it bounded} subsets
of $\Cal P_{\Cal C}(t)$.  Indeed, if $\{\xi_\alpha\}$ is a bounded
net in $\Cal P_{\Cal C}(t)$ which converges weakly to $\xi_\infty$,
then for every $\eta\in \Cal P_{\Cal C}(t)$ we have
$$
\align
\lim_{\alpha}\<\exp(\xi_\alpha),\exp(\eta)\> &=
\lim_\alpha e^{\<\xi_\alpha,\eta\>}
=e^{\<\xi_\infty,\eta\>}\\
&=\<\exp(\xi_\infty),\exp(\eta)\>.
\endalign
$$
Since $\exp(\xi_\alpha)$ is a bounded net in
the symmetric Fock space over $\Cal P_{\Cal C}(t)$
and since the set of vectors
$\{\exp(\eta):\eta\in \Cal P_{\Cal C}(t)\}$ span
this space, it follows that $\exp(\xi_\alpha)$ converges
weakly to $\exp(\xi_\infty)$.

Now choose a vector $v$ in the symmetric Fock space
over $\Cal P_{\Cal C}(t)$ with the property that $\<v,\exp(L)\>=\{0\}$.
The preceding paragraph implies that $v$ is orthogonal
to the set of vectors $\exp(\overline{L\cap B_r}^w)$ for
every $r>0$, and taking the union over $r>0$ we obtain
$\<v,\exp(K)\> = \{0\}$.  Since the exponential map
$\exp$ is {\it metrically} continuous on its entire
domain $\Cal P_{\Cal C}(t)$,  it follows that
$v$ is orthogonal to the set of vectors
$\exp(\overline{K})$, $\overline{K}$ denoting the
closure of $K$ in the norm topology of $\Cal P_{\Cal C}(t)$.
By remark 4.29 and Lemma 4.30, we conclude that $v=0$.

That completes the proof of Theorem 4.3\qed


\example{Classification of multiplicative forms}
We conclude this section with a discussion of how Theorem 4.3
gives a classification of multiplicative structures on path
spaces.  Let $(P,g)$ be a metric path space and let
$e^g:P^2\to \Bbb C$ be its associated multiplicative form,
defined on $P(t)\times P(t)$ for $t>0$ by
$$
x,y\in P(t)\mapsto e^{g(x,y)}.
$$
$e^g$ is a positive definite function on $P(t)\times P(t)$
and hence there is a Hilbert space $E(t)$ and a function
$F_t:P(t)\to E(t)$ satisfying
$$
\align
E(t) =& \overline{\text{span}}\{F_t(x): x\in P(t)\}\\
\<F_t(x),F_t(y)\>=&e^{g(x,y)},\qquad x,y\in P(t).
\endalign
$$
Thus we have a family of Hilbert spaces $p:E\to (0,\infty)$
$$
E=\{(t,\xi):t>0,\xi\in E(t)\}
$$
with projection $p(t,\xi) = t$.

We define a binary operation
$\xi,\eta\in E\mapsto \xi\cdot\eta\in E$ as
follows.  Since $g$ is additive there is a function
$\psi:P\times P\to \Bbb C$ such that for all
$x_1,x_2\in P(s)$, $y_1,y_2\in P(t)$
$$
g(x_1y_1,x_2y_2) = g(x_1,x_2)+g(y_1,y_2)+
\psi(x_1,y_1)+\overline{\psi(x_2,y_2)}.
$$
For $x\in P(s)$, $y\in P(t)$ we try to define the product
$F_s(x)\cdot F_t(y)$ by
$$
F_s(x)\cdot F_t(y) = e^{-\psi(x,y)}F_{s+t}(xy).  \tag{4.31}
$$
It follows that for $x_1,x_2\in P(s)$,
$y_1,y_2\in P(t)$ we have
$$
\align
\<F_s(x_1)\cdot F_t(y_1),F_s(x_2)\cdot F_t(y_2)\> &=
e^{g(x_1y_1,x_2y_2)-\psi(x_1,y_1)-\overline{\psi(x_2,y_2)}}
=e^{g(x_1,x_2)+g(y_1,y_2)}\\
&=\<F_s(x_1),F_s(x_2)\>\<F_t(y_1),F_t(y_2)\>.
\endalign
$$
The latter formula implies that there is a unique unitary
operator
$$
W_{s,t}:E(s)\otimes E(t)\to E(s+t)
$$
satisfying
$$
W_{s,t}(F_s(x)\otimes F_t(y))=F_s(x)\cdot F_t(y).
$$
Thus we can define a bounded bilinear map
$(\xi,\eta)\in E(s)\times E(t)\mapsto \xi\cdot\eta\in E(s+t)$
by way of
$$
\xi\cdot\eta = W_{s,t}(\xi\otimes \eta)
$$
and this mapping extends the operation (4.31).

To see that this operation on $E$ is associative, it suffices to
show that it is associative on generators, i.e.,
$$
F_r(x)\cdot(F_s(y)\cdot F_t(z))= (F_r(x)\cdot F_s(y))\cdot F_t(z)
$$
for all $x\in P(r), y\in P(s), z\in P(t)$.  Using the definition
(4.31), one observes that this will follow provided that
$\psi$ satisfies
$$
\psi(x,y)+\psi(xy,z) = \psi(x,yz) + \psi(y,z).  \tag{4.32}
$$
In fact, the equation (4.32) can be arranged {\it a priori}
from the definition of $\psi$ (2.4).  But it is easier at this
point to invoke Theorem 4.3.  The latter asserts that there is
a complex-valued function $\rho$ defined on $P$ such that
$$
\psi(x,y) = \rho(xy)-\rho(x)-\rho(y)
$$
for all $x,y\in P$.  Substituting this into (4.32) one finds
(using associativity of the multiplication in $P$)
that both sides of (4.32) reduce to
$$
\rho(xyz)-\rho(x)-\rho(y)-\rho(z).
$$
This proves associativity of the multiplication in $E$.

The preceding discussion implies that this multiplication acts
like tensoring.  Thus we have a product structure satisfying
all the axioms of a product system except measurability requirements.
Using Theorem 4.3, we can describe this product structure as follows.

\proclaim{Corollary 4.33}
Let $(P,g)$ be a metric path space and let $E$ be the product
structure obtained from the positive definite functions
$$
(x,y)\in P(t)\times P(t) \mapsto e^{g(x,y)},
$$
$t>0$.  Assume that $E(t)$ is not one-dimensional for every $t$.
Then $E$ is isomorphic to the product structure of one of the standard
product systems $E_1,E_2,\dots, E_\infty$.
\endproclaim
\demo{proof}
By Theorem 4.3, there is a separable Hilbert space $\Cal C$, a
complex-valued function $\rho:P\to \Bbb C$ and a fiber map
$\log:P\to \Cal P_{\Cal C}$ satisfying the conditions (4.3.1)
and (4.3.2).  In view of the remarks following Theorem 4.3,
 we may take the defect to be of
the form
$$
\psi(x,y) = \rho(xy)-\rho(x)-\rho(y).
$$
Using the formula (4.31) we find that if we rescale
$F_t:P(t)\to E(t)$ according to
$$
g_t(x)=e^{-\rho(x)}F_t(x),
$$
then the definition of multiplication in $E$ simplifies to
$$
G_s(x)\cdot G_t(y) = G_{s+t}(xy),
$$
for $x\in P(s), y\in P(t)$, and all $s,t>0$.

Moreover, formula (4.3.2) implies that for $x_1,x_2\in P(t)$
we have
$$
\<G_t(x_1),G_t(x_2)\>=e^{g(x_1,x_2)-\rho(x_1)-\overline{\rho(x_2)}}
=e^{\<\log(x_1),\log(x_2)\>}.
$$

Now consider the exponential map
$$
\exp:L^2((0,\infty);\Cal C)\to e^{L^2((0,\infty);\Cal C)}.
$$
The latter formula asserts that
$$
\<G_t(x_1),G_t(x_2)\> = \<\exp(\log(x_1)),\exp(\log(x_2))\>
$$
for all $x_1,x_2\in P(t)$.  This implies that we can define
an isometry
$$
W_t:E(t)\to e^{\Cal P_{\Cal C}(t)}
$$
by way of
$$
W_t(G_t(x)) = \exp(\log(x)),
$$
for $x\in P(t)$.  By (4.3.1), each $W_t$ is a unitary operator.
The total map
$$
W:E\to e^{\Cal P_{\Cal C}}
$$
is an isomorphism of families of Hilbert spaces.

It remains to verify that $W$ preserves
multiplication, i.e., that
$$
W_{s+t}(\xi\cdot\eta)= W_s(\xi)W_t(\eta)
$$
for every $\xi\in E(s)$, $\eta\in E(t)$.  Recalling that the
multiplication in $e^{\Cal P_{\Cal C}}$ is defined by
$$
\exp(f)\exp(g) = \exp(f\boxplus g)
$$
for $f\in \Cal P_{\Cal C}(s)$, $g\in \Cal P_{\Cal C}(t)$ we find that
for all $x\in P(s)$, $y\in P(t)$,
$$
\align
W_{s+t}(G_s(x)\cdot G_t(y))&=W_{s+t}(G_{s+t}(xy))=\exp(\log(xy))\\
&=\exp(\log(x)\boxplus\log(y))=\exp(\log(x))\exp(\log(y)),
\endalign
$$
and hence
$$
W_{s+t}(G_s(x)\cdot G_t(y))=W_s(G_s(x))W_t(G_t(y)).
$$
The assertion follows from the bilinearity of multiplication and
the fact that $E(r)$ is spanned by $G_r(P(r))$ for every $r>0$

Finally, note that if the space $\Cal C$ of ``coordinates" is
the trivial Hilbert space $\{0\}$
then $\Cal P_{\Cal C}(t)=L^2((0,t);\Cal C)$ is trivial as
well and hence
$$
e^{\Cal P_{\Cal C}(t)}
$$
is one-dimensional for every $t>0$.  By virtue of the isomorphism
$W:E\to e^{\Cal P_{\Cal C}}$, this has been ruled out in the hypothesis
of Corollary 4.33.  Thus $n=\dim(\Cal C)$ is a positive integer
or $\aleph_0$.  In this case, $W$ implements an isomorphism
of the product structure $E$ onto the product structure of $E_n$
\qed
\enddemo
\endexample

\heading{Part II.  Continuous Tensor Products}
\endheading

\subheading{Introduction to Part II}

Let $p:E\to (0,\infty)$ be a product system.  Thus
each fiber $E(t)=p^{-1}(t)$ is a separable Hilbert space
and we are given an associative multiplication
$(x,y)\in E\times E\mapsto xy\in E$ which acts
like tensoring in the sense that for fixed $s,t>0$,
$$
(x,y)\in E(s)\times E(t)\mapsto xy\in E(s+t)
$$
is a bilinear mapping with the properties
$$
\align
E(s+t) &= \overline{\text{span}}E(s)E(t)  \tag{II.1}\\
\<x_1y_1,x_2y_2\>&=\<x_1,x_2\>\<y_1,y_2\>, \tag{II.2}
\endalign
$$
for $x_i\in E(s)$, $y_i\in E(t)$.  In addition, there are
natural measurability axioms which we will not repeat here
\cite{2}.  We will write $E=\{E(t):  t>0\}$ instead
of $p:E\to (0,\infty)$ when it is convenient.

A nonzero vector $x\in E(t)$ is called {\it decomposable}
if for every $0<s<t$ there are vectors $y\in E(s)$,
$z\in E(t-s)$ for which
$$
x = y z. \tag{II.3}
$$
The set of decomposable vectors in $E(t)$ will be written
$D(t)$.

There are product systems which contain
no decomposable vectors at all.  But if
there is a $t_0>0$ for which $D(t_0)\neq \emptyset$, then
$D(t)\neq \emptyset$ for every $t>0$ and we clearly have
$$
D(s+t) = D(s)D(t)
$$
for every $s,t>0$.  This multiplicatively structured family of
sets $\{D(t): t>0\}$ comes close to defining a path space, except
that the factorizations of (II.3) are not unique.  However,
if $y_1,y_2\in E(s)$ and $z_1,z_2\in E(t-s)$
satisfy $y_1z_1=y_2z_2$, then because of the
identification of $E(t)$ with the tensor product
$E(s)\otimes E(t-s)$ described in (II.1) and (II.2), we see
that there must be a nonzero complex number $\lambda$ such that
$$
y_2=\lambda y_1,\qquad z_2=\lambda^{-1}z_1.
$$
Thus we may obtain a path space structure by passing from
each $D(t)$ to its associated projective space $\Delta(t)$.

More explictly, $\Delta(t)$ is obtained by identifying
two vectors $x_1,x_2$ in $D(t)$ which are
nonzero scalar multiples of each other.
We consider $\Delta(t)$ to be a set with no additional structure.
There is a natural projection $x\in D(t)\mapsto \dot x\in \Delta(t)$.
Any complex function $f:D(t)\to \Bbb C$ which is homogeneous
of degree zero in the sense
that $f(\lambda x)=f(x)$ for all nonzero scalars $\lambda$
can be promoted to a function $\dot f:\Delta(t)\to \Bbb C$ by way of
$\dot f(\dot x) = f(x)$, $x\in D(t)$.  In fact, it will be convenient
to abuse notation slightly and identify functions on $\Delta(t)$ with
homogenous functions defined on $D(t)$. The path space
$p:\Delta \to (0,\infty)$ is defined by
$$
\Delta = \{(t,\dot x): t>0, x\in D(t)\},
$$
with projection $p(t,\dot x) = t$ and multiplication
$(s,\dot x)(t,\dot y) = (s+t,\dot{xy})$.
We remind the reader that $\Delta$,
like any path space, is to be considered
a fibered set with no additional structure beyond
the multiplication it carries.
$\Delta^2$ will denote the fiber product
$$
\Delta^2 = \{(t,\dot x,\dot y): t>0, \dot x,\dot y\in \Delta(t)\}.
$$
For example, if for each $t>0$ we are given a function $f_t:D(t)\to \Bbb C$
which satisfies $f_t(\lambda x) = f_t(x)$ for $x\in D(t)$ and $\lambda\neq 0$,
then according to the abuses that have been agreed to we can define a
function $\phi : \Delta\to \Bbb C$ by way of
$$
\phi(t,x) = f_t(x),\qquad x\in D(t).
$$

Of course, the inner product restricts to a
positive definite function on every $D(t)$
$$
(x,y)\in D(t)\times D(t)\mapsto \<x,y\>\in \Bbb C,
$$
but this function cannot be promoted to one defined on
$\Delta(t)\times\Delta(t)$.  We will see in section 6
below that the inner product of any two
vectors in $D(t)$ must be nonzero.
Thus if we choose a fixed element $e\in D(t)$
then we may form the renormalized inner product
$$
P_t(x,y)=\frac{\<x,y\>}{\<x,e\>\<e,y\>}
$$
and the latter function can be promoted to a
positive definite function on $\Delta(t)\times\Delta(t)$.
If we choose a family $\{e_t\in D(t): t>0\}$ of
decomposable vectors then we obtain a function
$$
P:\Delta^2\to \Bbb C
$$
which restricts to a positive definite function
on $\Delta(t)\times \Delta(t)$ for every $t>0$.
Of course, these renormalized versions of the inner
product depend on the particular choice of
$\{e_t:t>0\}$.

The results of sections 5--9 below combine to show that it is
possible to find additive forms $g: \Delta^2\to \Bbb C$
for this path space which are logarithms of the inner
product in the following sense.

\proclaim{Theorem A}
Let $\{e_t\in D(t): t>0\}$ be a family of decomposable
vectors which is left-coherent in the sense that for
every $s,t$ satisfying $0<s<t$
there is a vector $e(s,t)\in D(t-s)$ such
that
$$
e_t = e_se(s,t).
$$
Then there is an additive form $g:\Delta^2\to \Bbb C$
(which will depend on $e$) such that for every $t>0$ and
every $x_1,x_2\in D(t)$ we have
$$
\<x_1,x_2\> = \<x_1,e_t\>\overline{\<x_2,e_t\>}e^{g(x_1,x_2)}.
$$
\endproclaim

In fact, we will show that $g:\Delta^2\to \Bbb C$ is a
``continuous" logarithm which ``vanishes at zero";
moreover, it is uniquely determined by these requirements
once $\{e_t\}$ is fixed.
We emphasize that for typical choices of $e$ the form
$g$ will have {\it nonzero} defect.

In sections 5--7, we establish certain continuity and
nonvanishing properties of inner products of
decomposable vectors.  In section 8 we construct
$g$ as a ``continuous" logarithm, and in section 9
we show that it has the required positivity properties.

\subheading{5.  Continuity of the modulus}

Let $E= \{E(t): t>0\}$ be a product system.  Suppose that we are
given vectors $x\in E(s)$, $y\in E(t)$ with $0<s<t$.  $x$ is called
a {\it left divisor}  (resp. {\it right divisor}) of $y$ if there is a
vector $z\in E(t-s)$ such that $y = xz$ (resp. $y = zx$).  Notice that
in either case we have $\|x\|\cdot\|z\| = \|y\|$ and hence both
$x$ and $z$ are nonzero whenever $y$ is nonzero.  Notice too that,
while factorizations of the form II.3 are not unique, we do have
both left and right cancellation laws.  That is, if $y\in E(s)$
and $z_1,z_2\in E(t)$, then $yz_1=yz_2\implies z_1=z_2$, and
$z_1y=z_2y\implies z_1=z_2$.

Given $0 < T\leq\infty$, we say that a family of vectors
$\{x_t\in E(t): 0 < t < T\}$ is {\it left-coherent}
(resp. {\it right-coherent})if for every
$0<s_1<s_2<T$, $x_{s_1}$ is a left (resp. right) divisor of
$x_{s_2}$.  Our analysis is based on the following
continuity property of coherent families of vectors.  The
proof makes use of a central technical result from \cite{3}.

\proclaim{Theorem 5.1}Let $\{x_t: 0<t<T\}$ and $\{y_t: 0<t<T\}$
be two left-coherent (resp. right-coherent) families of vectors
satisfying $\|x_t\| = \|y_t\| = 1$ for all $t$.  Then
$$
\lim_{t\to0+} |\<x_t,y_t\>| = 1.
$$
\endproclaim

\remark{Remarks}  Notice that the quantity $|\<x_t,y_t\>|$ does not
exceed $1$ and increases as $t$ decreases to $0$.  Indeed, if
$0<s<t<T$ then in the left-coherent case we can write
$x_t = x_su$, $y_t=y_sv$ where $u,v$ are unit vectors in $E(t-s)$,
hence
$$
|\<x_t,y_t\>| = |\<x_su,y_sv\>| =
|\<x_s,y_s\>|\cdot|\<u,v\>|\leq|\<x_s,y_s\>|\leq 1.
$$
It follows that the essential assertion of Theorem 5.1 is that
$$
\sup_{0<t<T} |\<x_t,y_t\>| = 1.
$$

Notice too that we cannot draw any conclusion about continuity
of the inner product itself.  Indeed, if we start with $\{y_t\}$
as above and define $x_t$ by $x_t = f(t)y_t$ where $f$ is an arbitrary
function from $(0,\infty)$ to the unit circle in the
complex plane, then $<x_t,y_t> = f(t)$ can be pathological.
\endremark

\demo{proof of Theorem 5.1}
We will prove Theorem 5.1 for left-coherent families.  With that in hand,
the corresponding conclusion for right-coherent families follows from
it by considering the product system $E^o$ opposite to $E$ (see \cite{7,8}).
Moreover, by passing to subfamilies if necessary, we may assume that
$T$ is finite and positive.

For every $s\in(0,T)$ we consider a projection $P_s\in\Cal B(E(T))$
which is defined as follows.  Since $\|x_s\| = 1$, the operator
$$
L_{x_s}:z\in E(T-s) \to x_sz\in E_T
$$
is an isometry whose range projection is given by
$$
Ps = L_{x_s}L_{x_s}^* \in\Cal B(E(T)).
$$

The action of $P_s$ on a product vector of the form
$$
y = y_1y_2,\qquad y_1\in E(s), y_2\in E(T-s)
$$
is given by
$$
P_s y = \<y_1,x_s\>x_sy_2.
$$
Note that if $0<s_1<s_2<T$ then $P_{s_1}\geq P_{s_2}$, or equivalently
that
$$
x_{s_2}E(T-s_2) \subseteq x_{s_1}E(T-s_1).   \tag{5.2}
$$
Indeed, because of left-coherence $x_{s_2}$ can be factored into a product
$x_{s_2} = x_{s_1}z$ where $z\in E(s_2-s_1)$.  Hence
$zE(T-s_2)\subseteq E(T-s_1)$, from which (5.2) is evident.

Therefore, as $s$ decreases to $0$ the projections $P_s$ increase
strongly to a limit projection $P_0$:
$$
P_0 = \operatornamewithlimits{strong\,limit}_{s\to 0+} P_s.
$$
We claim that $P_0 = {\bold 1}$.  Since $P_0$ is clearly a
nonzero projection in $\Cal B(E(T))$, it suffices to
exhibit an irreducible $*$-algebra $\Cal A$ of operators on
$E(T)$ which commutes with $P_0$.  We can exhibit such an algebra
as follows.

For every $0<s<T$ we can identify $E(T)$ with the tensor product
$E(s)\otimes E(T-s)$ by associating a product $xy\in E(T)$,
$x\in E(s)$, $y\in E(T-s)$ with $x\otimes y\in E(s)\otimes E(T-s)$.
In this identification, $P_s$ is associated with the projection
$$
[x_s]\otimes{\bold 1}_{E(T-s)},
$$
$[x_s]$ denoting the rank-one projection
$\xi\in E(s)\mapsto \<\xi,x_s\>x_s\in E(s)$.  Thus $P_s$ commutes with
the von Neumann algebra
$$
\Cal A_s = {\bold 1}_s\otimes\Cal B(E(T-s)).
$$
As $s$ decreases to $0$ the operator algebras $\Cal A_s$ increase,
and hence $P_0$ commutes with the union
$$
\Cal A = \bigcup_{0<s<T}\Cal A_s.
$$
By Proposition 5.5 of \cite{3}, the union $\Cal A$ is irreducible.
Hence $P_0={\bold 1}$.

It follows that
$$
\lim_{s\to 0+}\|P_s\xi\| = \|\xi\|,
$$
for every vector $\xi\in E(T)$.  Taking $\xi = y_T$ and
noting that for every $0<s<T$, $y$ factors into a product
$y_T = y_sz_{T-s}$ where $z_r$ is a unit vector in $E(r)$
for every $0<r<T$, we have
$$
P_sy_T = P_s(y_sz_{T-s}) = \<y_s,x_s\>x_sz_{T-s}
$$
and hence
$$
1 = \|y_T\| = \lim_{s\to 0+}\|P_sy_T\| =
\lim_{s\to 0+}|\<y_s,x_s\>|\cdot\|x_s\|\cdot\|z_{T-s}\| =
\lim_{s\to 0+}|\<y_s,x_s\>|,
$$
as required\qed
\enddemo

With 5.1 in hand, we can now establish the following basic
result on continuity of inner products.

\proclaim{Theorem 5.3} Let $\{x_t: 0<t<T\}$ and $\{y_t: 0<t<T\}$
be two families of vectors which satisfy the hypotheses of Theorem
5.1.  Then the function
$$
\phi(t) = |\<x_t,y_t\>|
$$
is continuous on the interval $0<t<T$ and tends to $1$ as
$t\to 0+$.
\endproclaim
\demo{proof}
As in the proof of 5.1, it suffices to consider the case where
both families are left-coherent, and where $\|x_t\| = \|y_t\| = 1$
for all $t$.  Choose $t_0\in(0,\infty)$.  We will show that
$\phi$ is both left continuous and right continuous at $t_0$.

\demo{right continuity}  For every $\lambda\in(0,T-t_0)$ we may
find $u_\lambda$, $v_\lambda \in E(\lambda)$ such that
$$
\align
x_{t_0+\lambda} &= x_{t_0}u_\lambda,\\
y_{t_0+\lambda} &= y_{t_0}v_\lambda.
\endalign
$$
Notice that $\|u_\lambda\| = \|v_\lambda\| = 1$.  we claim
that $\{u_\lambda: 0<\lambda<T-t_0\}$ and
$\{v_\lambda: 0<\lambda<T-t_0\}$ are left-coherent families.
Indeed, if $0<\lambda_1<\lambda_2<T-t_0$ then since
$\{x_t\}$ is left-coherent we may find $z\in E(\lambda_2-\lambda_1)$
such that $x_{t_0+\lambda_2} = x_{t_0+\lambda_1}z$.  Hence
$$
x_{t_0}u_{\lambda_2} = x_{t_0+\lambda_2} = x_{t_0+\lambda_1} =
x_{t_0}u_{\lambda_1}z.
$$
{}From the left cancellation property we conclude
that $u_{\lambda_2}=u_{\lambda_1}z$, proving that
$\{u_\lambda\}$ is left-coherent.  The same is true of $\{v_\lambda\}$.

Thus by Theorem 5.1 we have
$$
\align
\lim_{\lambda\to0+}&|\<x_{t_0+\lambda},y_{t_0+\lambda}\>| =
\lim_{\lambda\to0+}|\<x_{t_0}u_\lambda,y_{t_0}v_\lambda\>| = \\
&|\<x_{t_0},y_{t_0}\>|\lim_{\lambda\to0+}|\<u_\lambda,v_\lambda\>|
= |\<x_{t_0},y_{t_0}\>|,
\endalign
$$
proving right continuity at $\lambda$.

\enddemo

\demo{left continuity}  Using left-coherence of
$\{x_t\}$ and $\{y_t\}$, we may find vectors $u_\lambda$,
$v_\lambda\in E_\lambda$ for every $\lambda\in(0,t_0)$
such that
$$
\align
x_{t_0} &= x_{t_0-\lambda}u_{\lambda}^\prime \tag{5.4}\\
y_{t_0} &= y_{t_0-\lambda}v_{\lambda}^\prime.
\endalign
$$
Notice that $\{u_\lambda^\prime\}$ and $\{v_\lambda^\prime\}$ are
right-coherent families.  The proof is similar to what was done above.
For example,  for $0<\lambda_1<\lambda_2<t_0$ we may find
$z^\prime\in E(\lambda_2-\lambda_1)$ such that
$$
x_{t_0-\lambda_1} = x_{t_0-\lambda_2}z^\prime.
$$
Hence
$$
x_{t_0-\lambda_2}u_{\lambda_2}^\prime = x_{t_0-\lambda_1}u_{\lambda_1}^\prime
= x_{t_0-\lambda_2}z^\prime u_{\lambda_1}^\prime,
$$
and we obtain
$$
u_{\lambda_2}^\prime = z^\prime u_{\lambda_1}^\prime
$$
after cancelling $x_{t_0-\lambda_2}$ from the left.
The proof that $\{v_\lambda^\prime\}$ is right-coherent is of
course the same.

Using (5.4) we have
$$
|\<x_{t_0},y_{t_0}\>| = |\<x_{t_0-\lambda},y_{t_0-\lambda}\>|
\cdot |\<u_\lambda^\prime, v_\lambda^\prime\>|
$$
for every $0<\lambda<t_0$.  Because of Theorem 5.1,
$$
\lim_{\lambda\to0+}|\<u_\lambda^\prime,v_\lambda^\prime\>| = 1,
$$
hence
$$
\lim_{\lambda\to0+}|\<x_{t_0-\lambda},y_{t_0-\lambda}\>| =
|\<x_{t_0},y_{t_0}\>|,
$$
proving left continuity at $t_0$.  The continuity of $\phi$ at
$t=0$ was established in Theorem 5.1\qed
\enddemo
\enddemo

\subheading{6.  Inner products of decomposable vectors}

Let $t>0$ and let $x\in D(t)$ be a decomposable vector.
Then for $0<s<t$ there are vectors $a_s$, $b_s\in E(s)$ such
that
$$
x = a_s b_{t-s}. \tag{6.1}
$$

\remark{Remark}  Note that if $\|x\| = 1$,
then $\|a_s\|\cdot\|b_{t-s}\|=1$
and we may perform an obvious renormalization to achieve
$\|a_s\| = \|b_s\| = 1$ for every $s$.

It is not obvious that each $a_s$ and each $b_s$ is a
decomposable vector.  The following lemma establishes
this fact, and moreover it shows that $\{a_s: 0<s<t\}$
(resp. $\{b_s: 0<s<t\}$) is a
left-decomposable (resp. right-decomposable) family.

\proclaim{Lemma 6.2}  Choose $s_1,s_2$ with $0<s_1<s_2<t$, and
suppose that $a_i\in E(s_i)$ and $b\in E(t-s_i)$ satisfy
$$
a_1b_1=a_2b_2.
$$
Then $a_1$ is a left-divisor of $a_2$ and $b_2$ is a
right-divisor of $b_1$; i.e., there are vectors
$c_1,c_2\in E(s_2-s_1)$ such that
$$
\align
a_2&=a_1c_1\\
b_1&=c_2b_2.
\endalign
$$
\endproclaim
\demo{proof}
We may assume that $x=a_1b_1=a_2b_2$ satisfies
$\|x\|=1$, and thus after an obvious renormalization
we may also assume $\|a_i\|=\|b_i\|=1$.


We require an enhanced version of the right cancellation law.  Note
that if $z$ is any element of $E$, say $z\in E(\lambda)$ for $\lambda >0$,
then the right multiplication map
$R_z:u\in E\mapsto uz\in E$ restricts to a bounded linear map on fiber
spaces, carrying $E(\mu)$ to $E(\mu + \lambda)$, and thus has a fiber
adjoint.  Let
$$
R_z^* : \{E(\mu): \mu>\lambda\} \to E
$$
be the total map defined by these adjoints.  We claim that that $R_z^*$
commutes with any left multiplication operator $L_a: u\mapsto au$
in the sense that $L_a$ commutes
with the restriction of $R_z^*$ to any fiber space
in any case in which the formulas make sense.
That is, if $a\in E(\lambda)$ and $u\in E(\mu)$ with
$\mu > \lambda$, then we have
$$
aR_z^*(u) = R_z^*(au).  \tag{6.3}
$$
To see this, simply note that since $E(\mu)$ is spanned by
$E(\mu-\lambda)E(\lambda)$, it suffices to verify that (6.3) is valid
for vectors $u$ of the form $u_1u_2$ with $u_1\in E(\mu-\lambda)$ and
$u_2\in E(\lambda)$.  In that case the left side of 6.3 is
$$
aR_z^*(u_1u_2) = a(u_1\<u_2,z\>)
$$
while the right side is
$$
R_z^*(a(u_1u_2)) = R_z^*((au_1)u_2) = au_1\<u_2,z\>,
$$
as asserted.

We apply these remarks to the proof of Lemma 6.2 as follows.
Write
$$
a_2=a_2\|b_2\|^2=R_{b_2}^*(a_2b_2)
=R_{b_2}^*(a_1b_1)=a_1R_{b_2}^*(b_1),
$$
and therefore we can take $c_1=R_{b_2}^*(b_1)\in E(s_2-s_1)$.

The other decomposition $b_1=c_2b_2$ follows from this by
considering the product system $E^o$ opposite to $E$.  Indeed,
if we interpret the equation $a_1b_1=a_2b_2$ in $E^o$, it
becomes $b_1a_1=b_2a_2$.  By what was proved above, there
is an element $c_2\in E^o(s_2-s_1)$ for which $b_1=b_2c_2$,
and if we interpret the latter in $E$ we obtain $b_1=c_2b_2$.
\qed
\enddemo

\remark{Remark 6.4}
We may conclude that every vector $x\in D(t)$ can be
associated with a propagator $\{x(r,s)\in D(s-r):0\leq r<s\leq t\}$
which satisfies $x(0,t)=x$.  Indeed, for each $0<s<t$ we
can find a nonzero left-divisor $x_s\in E(s)$ for $x$.
Set $x_t=x$.  By Lemma 6.2, $\{x_s:0<s\leq t\}$ is a
left-coherent family with $x_s\in D(s)$ for every $s$.
Because of the left cancellation law, we can therefore
define a propagator $\{x(r,s):0\leq r<s\leq t\}$ by setting
$$
x_s=x_rx(r,s)
$$
for $0<r<s\leq t$, and by setting $x(0,s)=x_s$.
\endremark

\proclaim{Theorem 6.5}  For any $t>0$ and any two vectors
$x,y\in D(t)$ we have $\<x,y\>\neq 0$.
\endproclaim
\demo{proof}
There is clearly no loss if we assume that
$\|x\| = \|y\| = 1$.  By the preceding remarks, we can find
propagators $\{x(r,s):0\leq r<s\leq t\}$ and
$\{y(r,s):0\leq r<s\leq t\}$ for $x$ and $y$ with the property
that $x(0,t)=x$ and $y(0,t)=y$.  By an obvious renormalization,
we can also assume that $\|x(r,s)\|=\|y(r,s)\|=1$ for every
$0\leq r<s\leq t$.

Notice that
for every $0<r\leq t$ we have
$$
\lim_{\lambda\to 0+}|\<x(r-\lambda,r),y(r-\lambda,r)\>|=1.  \tag{6.6}
$$
Indeed, this follows by applying 5.3 to the right-coherent
normalized sections $a_\lambda=x(r-\lambda,r)$ and
$b_\lambda=y(r-\lambda,r)$, $0<\lambda<r$.

Now the function $f:[0,t]\to \Bbb R$ defined by
$$
f(s) =
\cases
|\<x(0,s),y(0,s)\>|,& 0<s\leq t,\\
1,& s=0
\endcases
$$
is continuous, by Theorem 5.3.  We claim that $f$ is
never zero.  For if there were an $r\in [0,t]$ for which
$f(r)=0$, then there is a smallest one $r_0$, and we must
have $0<r_0\leq t$.  But for every $s\in (0,r_0)$ we
can write
$$
\align
f(r_0) =& |\<x(0,r_0),y(0,r_0)\>|=|\<x(0,s)x(s,r_0),y(0,s)y(s,r_0)\>|\\
=&f(s)|\<x(s,r_0),y(s,r_0)\>|.
\endalign
$$
$f(s)$ is nonzero for every
$s<r_0$, and because of (6.6) the term
$|\<x(s,r_0),y(s,r_0)\>|$ is nonzero when
$s$ is sufficiently close to $r_0$.  This contradicts the
fact that $f(r_0)$ was supposed to be zero.  Taking
$s=t$ we find that $|\<x,y\>| = f(t)\neq 0$\qed
\enddemo

\subheading{7.  Continuity and normalization of the inner product}

We will write $D$ for the set of all left-coherent decomposable
sections.  Thus $D$ consists of all sections
$$
t\in (0,\infty) \mapsto x_t\in E(t)
$$
which are left-coherent and for which $x_t$ is never zero.  It is possible,
of course, that $D=\emptyset$.  But if $D$ is not empty then we are
interested in establishing continuity of the inner product function
$$
t\in(0,\infty) \mapsto \<x_t,y_t\>,\tag{7.1}
$$
defined by two elements $x,y\in D$.
In this section we will show that if one normalizes the elements of
$D$ appropriately then inner products of the form (7.1) are
continuous.

\remark{Remarks}
Note that in general, nothing can be said about continuity of the
inner products (7.1) (see the remarks following Theorem 5.1).

Notice too that, even though unique factorization fails in the
multiplicative family $\{D(t):t>0\}$ we still have a
left cancellation law, and this implies that there is a
bijective correspondence between left-coherent sections and
propagators.  For example, if $x\in D(t)$ and we find
a left coherent family $\{x_s: 0<s\leq t\}$ for which
$x_t=x$ then we may define a unique propagator
$\{x(r,s)\in D(t-s): 0\leq r<s\leq t\}$ by
$$
x_s=x_rx(r,s)
$$
for $0<r<s$ (by the left cancellation property), and
where for $r=0$ we put $x(0,s)=x_s$.  The left-coherent family
is recovered from its propagator via $x_s=x(0,s)$, $0<s\leq t$.
\endremark

The normalization in $D$ is done as follows.  Choose an arbitrary
element $e\in D$ satisfying $\|e_t\|=1$ for every $t>0$; $e$
will be fixed throughout the remainder of this section.

\proclaim{Definition 7.3}  $D^e$ is the set of all $x\in D$
satisfying $\<x_t,e_t\> = 1$ for every $t>0$.
\endproclaim

\proclaim{Lemma 7.4}  For every $x\in D^e$, the norm
$\|x_t\|$ is a continuous nondecreasing function of
$t$ with
$$
\lim_{t\to 0+}\|x_t\| = 1.
$$
In particular, we have $\|x_t\|\geq 1$ for every $t>0$.
\endproclaim

\demo{proof}
Let $u_t = \|x_t\|^{-1}x_t$.  Then both $\{e_t\}$ and $\{u_t\}$ are
left-coherent families of {\it unit} vectors.  So by Theorems 5.1 and
5.3 we may conclude that $|\<u_t,e_t\>|$ is continuous in $t$ over the
interval $(0,\infty)$ and satisfies
$$
\lim_{t\to 0} |\<u_t,e_t\>| = 1.
$$
Since
$$
|\<u_t,e_t\>| = \|x_t\|^{-1}|\<x_t,e_t\>| = \|x_t\|^{-1},
$$
the continuity assertion follows.

To see that $\|x_t\|$ increases with $t$, choose $0<s<t$.  By left-coherence
of $\{x_t\}$ and $\{e_t\}$ we can write
$$
\align
e_t &= e_su,\\
x_t &=x_sv,
\endalign
$$
where $u=e(s,t)$, and $v=x(s,t)$ belong to $E(t-s)$.  Note that
$u$ must be a unit vector because $\|e(s,t)\| = \|e_t\|/\|e_s\| = 1$.
Notice too that $\|v\|\geq 1$.  Indeed,
since $\<e_r,x_r\> = 1$ for all $r$ we can write
$$
\<u,v\> = \<e_s,x_s\>\<u,v\> = \<e_su,x_sv\> = \<e_t,x_t\> = 1,
$$
so by the Schwarz inequality
$$
1 = |\<u,v\>|\leq \|u\|\cdot\|v\| = \|v\|.
$$
It follows that
$$
\|x_s\| \leq \|x_s\|\cdot\|v\| = \|x_sv\|=\|x_t\|
$$
as asserted \qed
\enddemo

Our principal result on the continuity of inner products is
the following.

\proclaim{Theorem 7.5}  Let $x,y\in D^e$.  Then the inner product
$\<x_t,y_t\>$ is continuous and nonzero on $0<t<\infty$, and satisfies
$$
\lim_{t\to 0+}\<x_t,y_t\> = 1.
$$
\endproclaim
\demo{proof}
We will deduce Theorem 7.5 from the following inequality.  For
every $s,t,T$ satisfying $0<s<t\leq T<\infty$ we claim that
$$
|\<x_s,y_s\> - \<x_t,y_t\>| \leq \|x_T\|\cdot \|y_T\|
\sqrt{(\|x_t\|^2 - \|x_s\|^2)(\|y_t\|^2-\|y_s\|^2)}. \tag{7.6}
$$
To prove (7.6), we write
$$
\align
x_t &=x_su\\
y_t &=y_sv\\
e_t &=e_sf
\endalign
$$
where $u=x(s,t),v=y(s,t),f=e(s,t)\in E(t-s)$ and $\|f\|=1$.  Notice that
$$
\<u,f\> = \<v,f\> = 1. \tag{7.7}
$$
Indeed, since $\<x_s,e_s\> = 1$ we have
$$
\<u,f\> = \<x_s,e_s\>\<u,f\> = \<x_su,e_sf\> = \<x_t,e_t\> = 1,
$$
and similarly $\<v,f\> = 1$.  We can therefore estimate the quantity
$\<u,v\> - 1$ as follows:
$$
|\<u,v\>-1| = |\<u-f,v-f\>| \leq \|u-f\|\cdot\|v-f\|.
$$
By Lemma 7.3 we have $\|x_s\|\geq 1$.  Hence we may use
$\<u,f\> = 1$ again to obtain
$$
\|u-f\|^2 = \|u\|^2 -1 = \frac{\|x_t\|^2}{\|x_s\|^2} -1
=\|x_s\|^{-2}(\|x_t\|^2-\|x_s\|^2) \leq \|x_t\|^2-\|x_s\|^2.
$$
Similarly,
$$
\|v-f\|\leq \sqrt{\|y_s\|^2 - \|y_s\|^2}.
$$
Thus
$$
|\<u,v\>-1| \leq \sqrt{(\|x_t\|^2-\|x_s\|^2)(\|y_t\|^2-\|y_s\|^2)}.
$$
The inequality (7.6) follows after multiplying the preceding
inequality through by
$|\<x_s,y_s\>|$, noting that
$$
\align
|\<x_s,y_s\>|\cdot|\<u,v\>-1| &= |\<x_s,y_s\>\<u,v\> - \<x_s,y_s\>| \\
&= |\<x_t,y_t\>-\<x_s,y_s\>|,
\endalign
$$
and using
the Schwarzenegger inequality and Lemma 7.4 to estimate
the factor $|\<x_s,y_s\>|$ on the right by way of
$$
|\<x_s,y_s\>|\leq\|x_s\|\cdot\|y_s\|\leq \|x_T\|\cdot\|y_T\|.
$$
This establishes (7.6).

Now Lemma 7.4 implies that $\|x_t\|^2$ and $\|y_t\|^2$
are continuous increasing
functions tending to $1$ as $t\to 0+$, and from
(7.6) we immediately conclude that $\<x_t,y_t\>$ is continuous on $(0,\infty)$.
If we allow $s$ to tend to $0$ in (7.6) and use $\lim_{t\to0+}\|x_s\|=1$
from Lemma 7.4 the result is
$$
|\<x_t,y_t\>-1| \leq \|x_T\|\cdot\|y_T\|
\sqrt{(\|x_t\|^2-1)(\|y_t\|^2-1)},
$$
from which we  deduce
$$
\lim_{t\to 0+}\<x_t,y_t\> = 1.
$$
That establishes continuity on the closed interval $[0,\infty)$ \qed
\enddemo

\subheading{8. Continuous logarithms}

Fix $t>0$.  We have seen above that if $x,y$ are two decomposable vectors
in $E(t)$ then the inner product $\<x,y\>$ is not zero.  Thus one
might attempt to define a logarithm  function
$(x,y)\in D(t)\mapsto L(t;x,y)\in \Bbb C$ with the property
that
$$
e^{L(t;x,y)} = \<x,y\>,
$$
in such a way that the logarithms
fit together consistently for different values of $t$.
We will show that this is in fact possible, provided that one is careful
to define the logarithm so as to remove as much ambiguity as
possible.  In section 9 we will show that $L$ is conditionally
positive definite.

Let $\Delta= \{\Delta(t): t>0\}$ be the path space obtained
from $\{D(t): t>0\}$ as in the introduction to Part II above.
Despite the fact that $\Delta$ and $\Delta^2$ are lifeless
sets, there is a useful notion of continuity for complex functions
defined on them.  We will say that $\phi:\Delta\to \Bbb C$ is
continuous if, for every left-coherent section
$t\in (0,\infty)\mapsto u_t\in D(t)$, the function $f(t) = \phi(t;u_t)$
is continuous on the interval $(0,\infty)$, and the limit
$$
f(0) = \lim_{t\to 0+}f(t)
$$
exists.  For functions $\psi:\Delta^2\to \Bbb C$, continuity
means that for any pair of sections $u,v\in D$, the function
$g(t)=\psi(t;u_t,v_t)$ is continuous for positive $t$ and extends
continuously to $[0,\infty)$.  We will say that $\phi$ (resp.\ $\psi$)
{\it vanishes at the origin} if $f(0)=0$ (resp.\ $g(t)=0$) for all choices
of $u$ (resp.\ $u,v$).

\remark{Remark 8.1}  Notice that certain normalized inner products
give rise to continuous functions $F:\Delta^2\to\Bbb C$.
For example, with $e$ as above put
$$
F(t;x,y) = \frac{\<x,y\>}{\<x,e_t\>\<e_t,y\>},
$$
for $t>0$, $x,y\in D(t)$.  To see that $F$ is continuous
choose $u,v\in D$ and put $u^\prime_t = \<u_t,e_t\>^{-1}u_t$
and $v^\prime_t = \<v_t,e_t\>^{-1}v_t$.  Then we have
$$
F(t;u_t,v_t) = \<u^\prime_t,v^\prime_t\>
$$
because of the homogeniety of $F$.  Theorem 7.5 implies that
the right side is continuous in $t$ and tends to $1$ as
$t\to 0+$.  It follows that $F:\Delta^2\to\Bbb C$
is a continuous function.  Needless to say, $F$ depends on $e$.
\endremark

\proclaim{Theorem 8.2}  Let $e\in D$ satisfy $\|e_t\|=1$ for
every $t>0$.  Then there is a unique continuous function
$L^e:\Delta^2\to\Bbb C$ which vanishes at the origin and
satisfies
$$
e^{L^e(t;x,y)}= \frac{\<x,y\>}{\<x,e_t\>\<e_t,y\>},\tag{8.2a}
$$
for every $t>0$, $x,y\in D(t)$.

If $f\in D$ satisfies $\|f_t\|=1$ for every $t>0$ and
$L^f:\Delta^2\to \Bbb C$ is the corresponding logarithm,
then there is a continuous function $\phi:\Delta\to\Bbb C$
which vanishes at $0$ and satisfies
$$
L^f(t;x,y)=L^e(t;x,y)+ \phi(t;x) + \overline{\phi(t;y)}\tag{8.2b}
$$
for all $t>0$, $x,y\in D(t)$.
\endproclaim

\remark{Remark 8.3}  The function $L^e:\Delta^2\to \Bbb C$ is
called the {\it $e$-logarithm} of the inner product on $E$.
Notice that we can use $L^e$ to define a logarithm of
the non-normalized inner product in the following way.  For every
$t>0$ and every $x\in D(t)$ let $f(t;x)$ be a
complex number such that
$$
e^{f(t;x)} = \<x,e_t\>.
$$
The function $f$ need have no regularity properties whatsoever, and
may even be non-measurable.  Nevertheless, once we settle on
$f$ then we can define
$$
L(t;x,y) = L^e(t;x,y) + f(t;x) + \overline f(t;y)
$$
and this new $L$ will satisfy
$$
e^{L(t;x,y)} = \<x,y\>,\qquad x,y\in D(t), t>0.
$$
We will see below that $L^e(t;\cdot,\cdot)$ is a positive definite
function on $D(t)\times D(t)$, and hence $L(t;\cdot,\cdot)$ is
a {\it conditionally positive definite logarithm of the inner product}
$\<\cdot,\cdot\>:D(t)\times D(t)\to \Bbb C$.  Now the construction
of a Hilbert space from a conditionally positive definite function
has the property that $L(t;\cdot,\cdot)$ and $L^e(t;\cdot,\cdot)$
determine the same Hilbert space (see section 7).  It follows that
the function $f$ has no effect on the invariantly defined Hilbert spaces
that concern us.  Moreover, for the same reason {\it (8.2b) implies that
these Hilbert spaces will also be independent of the particular choice of
normalized section $e\in D$}.
\endremark

\demo{proof of Theorem 8.2}  For uniqueness, notice that if
$L^e, M^e:\Delta^2\to \Bbb C$ both satisfy the conditions associated
with (8.2a), then
$$
\phi(t;x,y) = L^e(t;x,y) - M^e(t;x,y)
$$
is a continuous complex-valued function on $\Delta^2$ which vanishes
at the origin and satisfies
$$
e^{\phi(t;x,y)} = 1
$$
identically.  To see that $\phi = 0$ choose $t>0$ and $x,y\in D(t)$.
Let $u,v\in D$ be left-coherent sections such that $u_t$ and $v_t$
are, respectively, scalar multiples of $x$ and $y$
(see Theorem 10.1).  The function
$s\in(0.\infty)\to\phi(t;u_s,v_s)\in\Bbb C$ is continuous, vanishes
as $s\to 0+$ and satisfies
$$
e^{\phi(s;u_s,v_s)} = 1
$$
for all $s>0$.  Hence $\phi(s;u_s,v_s)=0$ for all $s$.  By homogeneity,
$\phi(t;x,y) = \phi(t;u_t,v_t) = 0$.

For existence, fix $t>0$, $x,y\in D(t)$.  We define $L^e(t;x,y)$ as
follows.  Again, we find sections $u,v\in D$ such that
$u_t=\lambda x$, $v_t=\mu y$, with $\lambda \mu\neq 0$.  Now the
function
$$
s\in(0,\infty)\mapsto \frac{\<u_s,v_s\>}{\<u_s,e_s\>\<e_s,v_s\>}
$$
is continuous, never $0$, and tends to $1$ as $s\to 0+$ (see Remark 8.1).
Thus there is a unique continuous function $l\in C[0,\infty)$ satisfying
the conditions $l(0)=0$ and
$$
e^{l(s)} = \frac{\<u_s,v_s\>}{\<u_s,e_s\>\<e_s,v_s\>},\qquad s>0.\tag{8.4}
$$
We define $L^e(t;x,y) = l(t)$.

To see that $L^e(t;x,y)$ is well-defined, choose another pair
$u^\prime, v^\prime\in D$ so that $u^\prime_t=\lambda^\prime x$ and
$v^\prime_t=\mu^\prime y$ with $\lambda^\prime \mu^\prime\neq 0$.
Choose $l^\prime\in C[0,\infty)$ with $l^\prime(0)=0$ so that
(8.4) is satisfied with $u^\prime$, $v^\prime$ replacing $u$, $v$.
We have to show that $l^\prime(t) = l(t)$.  But for $0<s\leq t$ the
uniqueness of factorizations of the two vectors $x,y$ implies that
there are nonzero complex numbers $\alpha_s,\beta_s$ so that
$u_s^\prime = \alpha_s u_s$, $v_s^\prime=\beta_s v_s$ for
$0<s\leq t$.  It follows that the right side of (8.4) is unaffected
by passing from $u,v$ to $u^\prime,v^\prime$.  Hence $l^\prime=l$
and finally $l^\prime(t) = l(t)$.

Notice that the continuity of $L^e$ follows from its definition.  Indeed,
for $u,v$ and $l$ related by (8.4) in the definition of $L^e(t;x,y)$,
we must also have
$$
L^e(s;u_s,v_s) = l(s)
$$
for every $0<s\leq t$.  In particular, the function
$s\in (0,t]\mapsto L^e(s;u_s,v_s)$ is continuous and tends to
$0$ as $s\to 0+$.  Since $t,x,y$ can be chosen arbitrarily, the sections
$u,v\in D$ are also arbitrary.  It follows that $L^e:\Delta^2\to \Bbb C$
is continuous and vanishes at $0$.

To prove (8.2b), pick $f\in D$ so that $\|f_t\|=1$ for every $t$, and
consider the function $\Phi:\Delta\to \Bbb C$ defined by
$$
\Phi(t;x)=\frac{|\<e_t,f_t\>| \<x,e_t\>}{\<f_t,e_t\>\<x,f_t\>},
\qquad t>0, \quad x,y\in D(t).
$$
We claim that for every $u\in D$, $\Phi(t,u_t)$ is continuous in
$t\in (0,\infty)$ and tends to $1$ as $t\to0+$.  Indeed, putting
$u^\prime=\<u_t,f_t\>^{-1}u_t$ and $e^\prime_t=\<e_t,f_t\>^{-1}e_t$,
then $e^\prime$ and $u^\prime$ are elements of $D$ satisfying
$\<e^\prime_t,f_t\>=\<u^\prime_t,f_t\>=1$, and
$$
\Phi(t;u_t)=|\<e_t,f_t\>| \<u^\prime_t,e^\prime_t\>
$$
for $t>0$.  The claim follows because both $|\<e_t,f_t\>|$ and
$\<u^\prime_t,e^\prime_t\>$ are continuous in $t$ and tend to
$1$ as $t\to 0+$ by Theorem 7.5.

To define $\phi$ we proceed as we did in defining
the function $L^e$ above.  Fix
$t>0$, $x\in D(t)$, and choose $u\in D$ so that $u_t=\lambda x$ for
some complex number $\lambda\neq 0$.  By the preceding paragraph
there is a unique continuous function $l\in C[0,t]$ so that
$l(0)=0$ and
$$
e^{l(s)} = \Phi(s;u_s),\qquad 0<s\leq t.
$$
Put $\phi(t;x) = l(t)$.  One shows that $\phi$ is well-defined and
continuous as one did for $L^e$.

Finally, since both $L^f(t;x,y)-L^e(t;x,y)$ and
$\phi(t;x) + \overline{\phi(t;y)}$ define continuous functions
on $\Delta^2$ which vanish at $0$, (8.2b) will follow if we show that
$$
e^{L^f(t;x,y)-L^e(t;x,y)} = e^{\phi(t;x)+\overline{\phi(t;y)}}.\tag{8.5}
$$
The left side of (8.5) is the quotient of
$$
e^{L^f(t;x,y)} = \frac{\<x,y\>}{\<x,f_t\>\<f_t,y\>}
$$
by the quantity
$$
e^{L^e(t;x,y)} = \frac{\<x,y\>}{\<x,e_t\>\<e_t,y\>}.
$$
Thus the $\<x,y\>$ terms cancel out and the left side of (8.5) reduces to
$$
\frac{\<x,e_t\>\<e_t,y\>}{\<x,f_t\>\<f_t,y\>}.  \tag{8.6}
$$
Similarly, the right side is the product of
$$
e^{\phi(t;x)} = \frac{|\<e_t,f_t\>| \<x,e_t\>}{\<f_t,e_t\>\<x_t,f_t\>}
$$
with the quantity
$$
\overline{e^{\phi(t;y)}} = \frac{|\<e_t,f_t\>|
\<e_t,y\>}{\<e_t,f_t\>\<f_t,y\>}.
$$
After performing the indicated multiplication all terms involving
$\<e_t,f_t\>$ cancel and the result agrees with (8.6) \qed
\enddemo

\remark{Remark 8.7} We remark that for each $e$, the $e$-logarithm is
self-adjoint in the sense that
$$
\overline{L^e(t;x,y)} = L^e(t;y,x), \qquad t>0,\quad x,y\in D(t).
$$
To see that, simply note that the function $F:\Delta^2\to \Bbb C$
defined by
$$
F(t;x,y) = \overline{L^e(t;y,x)}
$$
has all of the defining properties of an $e$-logarithm, and hence
$F=L^e$ by uniqueness.
\endremark

In the applications of part III, we will need to know that
the function $L^e$ defines an addtitive form on $\Delta^2$.
The following establishes this fact.

\proclaim{Proposition 8.8}
Fix $s>0$.  Then there is a continuous function
$\psi_s:\Delta\to \Bbb C$, vanishing at $0$, such that
for all $x_1,x_2\in \Delta(s)$, all $t>0$ and
all $y_1,y_2\in \Delta(t)$
we have
$$
L^e(s+t;x_1y_1,x_2y_2)-L^e(s;x_1,x_2)-L^e(t;y_1,y_2)
= \psi_s(t;y_1) + \overline{\psi_s(t;y_2)}.  \tag{8.9}
$$
\endproclaim

\demo{proof}
Fix $s>0$.  We claim that there is a continuous function
$\psi_s:\Delta \to \Bbb C$ which vanishes at $0$ and satisfies
$$
e^{\psi_s(t;y)} =
\frac{|\<e(s,s+t),e_t\>|\<y,e_t\>}{\<y,e(s,s+t)\>\<e(s,s+t),y\>}
$$
for all $y\in D(t)$, $t>0$.  In order to define $\psi_s$, fix
$t>0$ and $y_0\in D(t)$.  Choose a left-decomposable section
$y\in D$ such that $y_t$ is  a scalar multiple of $y_0$.  The
function
$$
t\in (0,\infty)\mapsto |\<e(s,s+t),e_t\>|
$$
is continuous and tends to $1$ as $t\to 0+$ by Theorem 2.3.
Similarly, by Remark 6.1, the function
$$
t\in (0,\infty) \to \frac{\<y_t,e_t\>}{\<y_t,e(s,s+t)\>\<e(s,s+t),y_t\>}
$$
has the same properties.  Thus
$$
t\in (0,\infty)\to
\frac{|\<e(s,s+t),e_t\>|\<y_t,e_t\>}{\<y_t,e(s,s+t)\>\<e(s,s+t),y_t\>}
$$
is a continuous function which tends to $1$ as $t\to 0+$.  It follows
that there is a unique continuous function $l:[0,\infty)\to \Bbb C$
such that $l(0)=0$ and
$$
e^{l(t)} =
\frac{|\<e(s,s+t),e_t\>|\<y_t,e_t\>}{\<y_t,e(s,s+t)\>\<e(s,s+t),y_t\>}
$$
for $t>0$.  We define $\psi_s(t;y_0)=l(t)$.

$\psi_s(t;\cdot)$ is a homogeneous function of degree $0$ on
$D(t)$ and hence we may consider $\psi_s$ to be a function
defined on $\Delta$.  $\psi_s$ is continuous because
of the way it was defined.

It remains to show that $\psi_s$ satisfies (8.9).  For that,
it suffices to show that for any pair of left-decomposable
sections $t\to y_t, y^\prime_t\in D(t)$ we have
$$
L^e(s+t;x_1y_t,x_2y^\prime_t)-L^e(s;x_1,x_2)-L^e(t;y_t,y^\prime_t)
= \psi_s(t;y_t) + \overline{\psi_s(t;y^\prime_t)}
$$
for every $t>0$.  To see this, let $L(t)$ and $R(t)$
be the left and right sides of the preceding formula.
$R(t)$ is continuous on $(0,\infty)$
and tend to $0$ as $t\to 0+$ by definition of $\psi_s$,
and we claim that $L(\cdot)$ has these two properties as well.
Indeed, $L(t)$ is continuous for positive $t$ because of
the continuity of $L^e$.
To see that $L(t)\to 0$ as $t\to0+$,
consider the left-coherent sections $u,u^\prime$ defined by
$$
u_r=
\cases
(x_1)_r,& 0<r\leq s\\
(x_1)_sy_{r-s},& r>s
\endcases
,\qquad
u^\prime_r=
\cases
(x_2)_r,& 0<r\leq s\\
(x_2)_sy^\prime_{r-s},& r>s
\endcases
..
$$
By continuity of $L^e(r;u_r,u^\prime_r)$
at $r=s$ we obtain
$$
\align
\lim_{t\to 0+}L(t) &= \lim_{t\to 0+}(L^e(s+t;u_{s+t},u^\prime_{s+t})
-L^e(s;x_1,x_2) - L^e(t;y_t,y^\prime_t)) \\
&= L^e(s;u_s,u^\prime_s) -L^e(s;x_1,x_2)=0,
\endalign
$$
as asserted.

Since both $L(t)$ and $R(t)$ are continuous on $(0,\infty)$
and tend to $0$ as $t\to 0+$, it suffices to show that
$$
e^{L(t)} = e^{R(t)},\qquad \text{for } t>0.
$$
But
$$
\align
e^{R(t)}=
&\frac{|\<e(s,s+t),e_t\>|\<y_t,e_t\>}{\<y_t,e(s,s+t)\>\<e(s,s+t),e_t\>}
\cdot
\left(\overline{\frac{|\<e(s,s+t),e_t\>|\<y^\prime_t,e_t\>}
{\<y^\prime_t,e(s,s+t)\>\<e(s,s+t),e_t\>}}\right)\\
=&\frac{\<y_t,e_t\>\<e_t,y^\prime_t\>}
{\<y_t,e(s,s+t)\>\<e(s,s+t),y^\prime_t\>}
\endalign
$$
while
$$
e^{L(t)}=\frac{\<x_1,y_t,x_2y^\prime_t\>}
{\<x_1y_t,e_{s+t}\>\<e_{s+t},x_2y^\prime_t\>}\cdot
\frac{\<x_1,e_s\>\<e_s,y_2\>}{\<x_1,x_2\>}\cdot
\frac{\<y_1,e_t\>\<e_t,y_2\>}{\<y_1,y_2\>}.
$$
\enddemo

Using the formulas
$$
\align
\<x_1y_t,x_2y^\prime_t\>&=\<x_1,x_2\>\<y_t,y^\prime_t\>\\
\<x_1y_t,e_{s+t}\>&=\<x_1y_t,e_se(s,s+t)\>=\<x_1,e_s\>\<y_t,e(s,s+t)\>\\
\<e_{s+t},x_2y^\prime_t\>&=\<e_s,x_2\>\<e(s,s+t),y^\prime_t\>
\endalign
$$
and performing the obvious cancellations, we obtain
$$
e^{L(t)}=\frac{\<y_t,e_t\>\<e_t,y^\prime_t\>}
{\<y_t,e(s,s+t)\>\<e(s,s+t,y^\prime_t\>}
=e^{R(t)},
$$
as required\qed

\subheading{9. Infinite divisibility of the inner product}

We have indicated in remark 8.3 how to find functions $L$
of the form
$$
L(t;x,y) = L^e(t;x,y) + f(t,x) + \overline{f(t;y)}
$$
that are logarithms of the inner product restricted to
decomposable vectors:
$$
e^{L(t;x,y)}= \<x,y\>,\qquad t>0, \quad x,y\in D(t).
$$
It is essential for the constructions of Part III that
such an $L$ should have the property that for fixed $t>0$,
it defines a conditionally positive
definite function on $D(t)\times D(t)$.  According to the
remarks at the beginning of section 2, it would
be enough to exhibit a sequence of positive definite functions
$\Phi_n:D(t)\times D(t)\to \Bbb C$ (depending on $t$) such that
$$
\Phi_n(x,y)^n = \<x,y\>,\qquad x,y\in D(t)
$$
for every $n=1,2,\dots$.
Unfortunately, there are no natural candidates for the positive
definite functions $\Phi_n$.  Thus we will have to establish the
conditional positive definiteness of $L(t;\cdot,\cdot)$ directly,
by making use of the structure of the product system itself.
Actually, we will prove somewhat more than we require.

\proclaim{Theorem 9.2}  Let $e\in D$ satisfy $\|e_t\|=1, t>0$.
Then for every $t>0$ the function
$$
(x,y)\in D(t)\times D(t)\mapsto L^e(t;x,y)
$$
is positive definite.
\endproclaim

The proof of Theorem 9.2 will occupy the remainder of this section.
Note that once 9.2 has been proved, one can immediately deduce

\proclaim{Corollary 9.3}  For every $t>0$, the inner product of
$E$ restricts to an infinitely divisible positive definite function
on $D(t)\times D(t)$.
\endproclaim

\demo{proof of Theorem 9.2}  Let $e$ be an element of $D$ satisfying
$\|e_t\|=1$ for all $t>0$, which will be fixed throughout the remainder
of this section.  We will define a function
$$
P^e:\Delta^2\to \Bbb C
$$
with the property that each function
$$
(x,y)\in D(t)\times D(t)\mapsto P^e(t;x,y)
$$
is {\it obviously} positive definite, and which is also a continuous
$e$-logarithm that vanishes $0$.  The conclusion $L^e=P^e$
will then follow by the uniqueness assertion of Theorem 8.2 and hence
we obtain 9.2.

Turning now to the proof, fix $t>0$ and choose $x,y\in D(t)$ with
$\<x,e_t\>=\<y,e_t\>=1$.  Because of the normalization of $x$ and $y$
there are left-coherent families $\{x_s: 0<s\leq t\}$
and $\{y_s: 0<s\leq t\}$ satisfying $\<x_s,e_s\> = \<y_s,e_s\>=1$
for all $s$ with
the property that $x_t=x$ and $y_t=y$.  Moreover, the two families are
uniquely determined by these conditions because of the
uniqueness of factorizations.  If $I=(a,b]$ is a subinterval of
$(0,t]$ having positive length then we will write $x_I$ (resp. $y_I$)
for the value of the propagator $x(a,b)$ (resp. $y(a,b)$).  Finally,
if $\Cal P=\{0=s_0<s_1<\dots<s_m=t\}$ and
$\Cal Q=\{0=t_0<t_1<\dots<t_n=t\}$ are two partitions of $(0,t]$, we
will write $\Cal P\leq \Cal Q$ for the usual ordering $\Cal P\subseteq\Cal Q$.
Thus we can define a net of complex numbers on the increasing directed
set of partitions by
$$
B_{\Cal P}(t;x,y) = \sum_{I\in \Cal P}(\<x_I,y_I\>-1).  \tag{9.4}
$$

\proclaim{Proposition 9.5}  For each $t>0$ and every partition $\Cal P$
of $(0,t]$, $B_{\Cal P}(t;\cdot,\cdot)$ is a positive definite function
whose associated Hilbert space is separable.  The net is decreasing in
the sense that for
$\Cal P\leq \Cal Q$, $B_{\Cal P}(t;\cdot,\cdot)-B_{\Cal Q}(t;\cdot,\cdot)$
is a positive definite function on $D(t)\times D(t)$.
\endproclaim

\demo{proof}Let $\Cal P = \{0=t_0<t_1<\dots<t_n=t\}$.  To see that
$B_{\Cal P}(t;\cdot,\cdot)$ is positive definite, choose $x,y\in D(t)$
and let $\{x_s:0<s\leq t\}$ and $\{y_s:0<s\leq t\}$ be the unique families
defined above.  It will be convenient to write $x_k=x(t_{k-1},t_k)$,
$y_k = y(t_{k-1},t_k)$ and $e_k=e(t_{k-1},t_k)$.  Noting that
$\<x_k,e_k\>=\<y_k,e_k\>=1$ we have
$\<x_k,y_k\>-1 = \<x_k-e_k,y_k-e_k\>$ and thus we can write
$$
B_{\Cal P}(t;x,y) = \sum_{k=1}^n (\<x_k,y_k\>-1) =
\sum_{k=1}^n\<x_k-e_k,y_k-e_k\>.
$$
Notice that $x_k-e_k$ and $y_k-e_k$ belong to the Hilbert space
$E(t_k-t_{k-1})$ for every $k=1,2,\dots, n$.  Thus we can define
a vector-valued function
$$
F_{\Cal P}:D(t)\to E(t_1)\oplus E(t_2-t_1)\oplus\dots\oplus E(t_n-t_{n-1})
$$
by way of
$$
F_{\Cal P}(x) = (x_1-e_1,x_2-e_2,\dots,x_n-e_n).
$$
The preceding formula for $B_{\Cal P}(t;x,y)$ now becomes
$$
B_{\Cal P}(t;x,y) = \<F_{\Cal P}(x),F_{\Cal P}(y)\>.
$$
This formula shows that $B_{\Cal P}(t;\cdot,\cdot)$ is a positive
definite function.  Moreover, since the range of the function $F_{\Cal P}$
is contained in a separable Hilbert space, it follows that
the Hilbert space associated with $B_{\Cal P}(t;\cdot, \cdot)$ is separable.

It remains to show that if $\Cal P_1$ and $\Cal P_2$
are two partitions satisfying $\Cal P_1\leq \Cal P_2$ then
$B_{\Cal P_1}(t;\cdot,\cdot)-B_{\Cal P_2}(t;\cdot,\cdot)$
is a positive definite function.
Now since the partial order of positive definite functions
defined by ($B_1\leq B_2 \iff B_2-B_1$ is positive definite) is
transitive and since $\Cal P_2$ is obtained from $\Cal P_1$ by
a sequence of steps in which one refines a single interval at
every step, we can reduce to the case in which
$\Cal P_1 = \{0=s_0<s_1<\dots<s_m=t\}$ and
$\Cal P_2$ is
obtained by adding a single point $c$ to $\Cal P_1$, where
$s_{k-1} < c < s_k$ for some $k=1,2,\dots,m$.  In this case the
difference $\Delta = B_{\Cal P_1} - B_{\Cal P_2}$ is given by
$$
\align
\Delta(x,y) =& \<x(s_{k-1},s_k),y(s_{k-1},s_k)\> - 1 \\
&-(\<x(s_{k-1},c),y(s_{k-1},c)\> +\<x(c,s_k),y(c,s_k)\>-2)\\
=&\<x(s_{k-1},s_k),y(s_{k-1},s_k)\>-\<x(s_{k-1},c),y(s_{k-1},c)\>
-\<x(c,s_k),y(c,s_k)\>\tag{9.6} \\
&+ 1.
\endalign
$$
If we write $x_1=x(s_{k-1},c)$, $x_2=x(c,s_k)$, $y_1=y(s_{k-1},c)$,
$y_2=y(c,s_k)$
then the right side of (9.6) can be rewritten as follows
$$
\align
\<x_1x_2,y_1y_2\>& - \<x_1,y_1\> - \<x_2,y_2\> +1\\
&= \<x_1,y_1\>\<x_2,y_2\>-\<x_1,y_1\>-\<x_2,y_2\>+1\\
&= (\<x_1,y_1\>-1)(\<x_2,y_2\>-1) = \<x_1-e_1,y_1-e_1\>\<x_2-e_2,y_2-e_2\>\\
&= \<(x_1-e_1)(x_2-e_2),(y_1-e_1)(y_2-e_2)\>
\endalign
$$
where $e_1=e(s_{k-1},c)$, $e_2=e(c,s_k)$, and where the inner product
in the last term on the right
is taken in the Hilbert space $E((s_k-c)+(c-s_{k-1}) = E(s_k-s_{k-1})$.
The last term clearly defines a positive definite function
of $x$ and $y$\qed
\enddemo

The conditions of Proposition 9.5 imply that the pointwise
limit $\lim_{\Cal P}B_{\Cal P}(t;x,y)$ exists.  This is a consequence of the
following elementary result.

\proclaim{Lemma 9.7}  Let $I$ be a directed set and let
$\{P_\alpha: \alpha\in I\}$ be a net of positive definite functions
on a set $X$
which is decreasing in the sense that $\alpha\leq\beta$ implies
that $P_\alpha-P_\beta$ is positive definite.  Then
$$
P_\infty(x,y) = \lim_\alpha P_\alpha(x,y)
$$
exists for every $x,y\in X$ and $P_\infty$ is a positive definite
function.  If the Hilbert space associated with some $P_\alpha$
is separable then so is the Hilbert space associated with $P_\infty$.
\endproclaim

\demo{proof}
Since a pointwise limit of positive definite functions is obviously
positive definite, we merely show that the above limit exists and defines
a separable Hilbert space.

Fix two elements $x,y\in X$.
For every $\alpha\in I$ consider the $2\times 2$
complex matrix
$$
A_\alpha =
\pmatrix
P_\alpha(x,x) & P_\alpha(x,y)\\
P_\alpha(y,x) & P_\alpha(y,y)
\endpmatrix.
$$
We may consider $\{A_\alpha: \alpha\in D\}$ as a net of
self adjoint operators on the
two dimensional Hilbert space $\Bbb C^2$.  We have $A_\alpha\geq 0$ for
every $\alpha$ because $P_\alpha$ is positive definite, and
$\alpha \leq \beta\implies A_\beta \leq A_\alpha$ because the net
$P_\alpha$ is decreasing.  Hence the net of operators $A_\alpha$ must
converge in the weak operator topology to a positive operator
$$
A = \lim_\alpha A_\alpha.
$$
Considering $A$ as a $2\times 2$ matrix, the element $\lambda$ in the
12 position satisfies
$$
\lambda = \lim_\alpha P_\alpha(x,y),
$$
establishing the existence of the required limit.

For separability, notice that
there are Hilbert spaces $H_\alpha$, $H_\infty$ and functions
$F_\alpha:X\to H_\alpha$, $F_\infty:X \to H_\infty$ with the
property
$$
\align
P_\alpha(x,y) &= \<F_\alpha(x),F_\alpha(y)\>,\\
P_\infty(x,y) &= \<F_\infty(x),F_\infty(y)\>,
\endalign
$$
and where we may also assume $H_\alpha$ (resp. $H_\infty$) is spanned by
$F_\alpha(X)$ (resp. $F_\infty(X)$).  By hypothesis,
we can find $\alpha$ so that $H_\alpha$ is separable.  Since
$P_\alpha - P_\infty$ is positive definite it follows that there is
a unique contraction $T:H_\alpha\to H_\infty$ having the property
$T(F_\alpha(x)) = F_\infty(x)$ for every $x\in X$.  Thus $T$ maps
$H_\alpha$ onto a dense subspace of $H_\infty$.  Since $H_\alpha$ is
separable we conclude that $H_\infty$ is separable as well.\qed
\enddemo

By 9.5 and 9.7, we may define a positive definite function
$P_\infty(t;\cdot,\cdot)$ on $D(t)\times D(t)$ by
$$
B_\infty(t;x,y)=\lim_{\Cal P}B_{\Cal P}(t;x,y).
$$
Finally, we define $P^e:\Delta^2\to \Bbb C$ by
$$
P^e(t;x,y) = B_\infty(t;\<x,e_t\>^{-1}x,\<y,e_t\>^{-1}y).
$$
It remains to show that $P^e$ is an $e$-logarithm, i.e., that it
is continuous, vanishes at the origin, and exponentiates correctly.

We deal first with continuity and vanishing at $0$.  Choose
$u,v\in D$ so that $\<u_t,e_t\>=\<v_t,e_t\>=1$ for every $t$.  We
have to show that the function
$$
t\in (0,\infty)\mapsto P^e(t;u_t,v_t) = B_\infty(t;u_t,v_t)
$$
is continuous on $(0,\infty)$ and tends to $0$ as $t\to 0+$.  This
will follow from the following two estimates.

\proclaim{Proposition 9.8}  If $0<s<t$ and $u,v$ are as above, then
$$
\align
|P^e(s;u_s,v_s)|\leq& (\|u_s\|^2 -1)(\|v_s\|^2-1),\\
|P^e(t;u_t,v_t)-P^e(s;u_s,v_s)|\leq&(\|u_t\|^2-\|u_s\|^2)(\|v_t\|^2-\|v_s\|^2).
\endalign
$$
\endproclaim
\demo{proof}Consider the first of the two inequalities.  Because
of the fact that
$$
P^e(s;u_s,v_s) = \lim_{\Cal P}B_{\Cal P}(s;u_s,v_s),
$$
it suffices
to show that for every partition $\Cal P = \{0=s_0<s_1<\dots<s_m=s\}$ of
the interval $(0,s]$ we have
$$
|B_{\Cal P}(s;u_s,v_s)|^2 \leq (\|u_s\|^2-1)(\|v_s\|^2-1).  \tag{9.9}
$$
For that,
let us write $u^k=u(s_{k-1},s_k)$, $v^k=v(s_{k-1},s_k)$,
$e^k = e(s_{k-1},s_k)$ for $k=1,2,\dots,m$.  Because of the normalizations
$\<u^k,e^k\> = \<v^k,e^k\>=1$ we have
$$
\<u^k,v^k\>-1 = \<u^k-e^k,v^k-e^k\>
$$
and hence
$$
\align
|B_{\Cal P}(s;u_s,v_s)|^2 &= |\sum_{k=1}^m\<u^k-e^k,v^k-e^k\>|^2
\leq (\sum_{k-1}^m\|u^k-e^k\|\cdot\|v^k-e^k\|)^2 \\
&\leq \sum_{k=1}^m\|u^k-e^k\|^2\sum_{k=1}^m\|v^k-e^k\|^2.
\endalign
$$
Now we can write
$$
\|u^k-e^k\|^2 = \|u^k\|^2 -1 = \|u(s_{k-1},s_k)\|^2 -1.
$$
If $k=1$ this is just $\|u_{t_1}\|^2-1$, and if $k>1$ it becomes
$$
\align
\frac{\|u_{s_k}\|^2}{\|u_{s_{k-1}}\|^2} - 1
&= \|u_{s_{k-1}}\|^{-2}(\|u_{s_k}\|^2-\|u_{s_{k-1}}\|^2)
\\&\leq \|u_{s_k}\|^2-\|u_{s_{k-1}}\|^2.
\endalign
$$
Thus we can estimate $\sum\|u^k-e^k\|^2$ using a telescoping series to
obtain
$$
\sum_{k=1}^m\|u^k-e^k\|^2 \leq \|u_{s_m}\|^2-1=\|u_s\|^2-1.
$$
Similarly,
$$
\sum_{k=1}^m\|v^k-e^k\|^2 \leq \|v_{s_m}\|^2-1=\|v_s\|^2-1,
$$
and the first of the two inequalities follows.

The proof of the second is similar, and we merely indicate the
changes.  It suffices to show that for any partition $\Cal P$ of
$(0,t]$ which contains $s$, we have
$$
|B_{\Cal P}(t;u_t,v_t)-B_{\Cal P\cap(0,s]}(s;u_s,v_s)|
\leq (\|u_t\|^2-\|u_s\|^2)(\|v_t\|^2-\|v_s\|^2).
$$
The desired inequality will follow by taking the limit on $\Cal P$.
Suppose that
$$
\Cal P = \{0=s_0<s_1<\dots<s_m=s=t_0<t_2<\dots<t_n=t\}.
$$
Then if we write out the formula for $B_{\Cal P}(t;u_t,v_t)$ we
find that
$$
B_{\Cal P}(t;u_t,v_t)= B_{\Cal P\cap(0,s]}(s;u_s,v_s)
+\sum_{l=1}^n (\<u(t_{l-1},t_l),v(t_{l-1},t_l)\>-1).
$$
Thus we have to show that
$$
|\sum_{l=1}^n (\<u(t_{l-1},t_l),v(t_{l-1},t_l)\>-1)|^2\leq
(\|u_t\|^2-\|u_s\|^2)(\|v_t\|^2-\|v_s\|^2).
$$
But if we write $u^l=u(t_{l-1},t_l)$, $v^l=v(t_{l-1},t_l)$,
$e^l=e(t_{l-1},t_l)$, then we have
$$
\<u^l,v^l\>-1 = \<u^l-e^l,v^l-e^l\>
$$
and as in the proof of the first inequality it suffices to show
that
$$
\sum_{l=1}^n\|u^l-e^l\|^2\leq \|u_t\|^2-\|u_s\|^2,
$$
with a similar estimate for $u$ replaced with $v$.  But noting that
$$
\{s=t_0<t_1<\dots<t_n=t\}
$$
is a partition of the interval $(s,t]$, we can make similar
estimates as those made in the preceding argument to estimate the
sum
$$
\sum_{l=1}^n\|u^l-e^l\|^2=\sum_{l=1}^n(\|u(t_{l-1},t_l)\|^2-1)
$$
with a telescoping series whose sum is
$\|u_t\|^2-\|u_s\|^2$\qed
\enddemo

{}From the inequalities of Proposition 9.8 and Lemma 7.4,
we immediately conclude that $P^e:\Delta^2\to \Bbb C$ is continuous
and vanishes at the origin.

Thus, to show that $P^e$ is an $e$-logarithm (and therefore coincides with
$L^e$) it remains only to show that for $t>0$ and $x,y\in D(t)$ we have
$$
e^{P^e(t;x,y)} = \frac{\<x,y\>}{\<x,e_t\>\<e_t,y\>}.
$$
Since both sides are homogeneous functions of degree zero in $x$ and $y$,
it suffices to prove the formula for $x,y$ normalized so that
$\<x,e_t\>=\<y,e_t\>=1$.  That is, we must prove that
$$
e^{B_\infty(t;x,y)} = \<x,y\>,\tag{9.10}
$$
for all $x,y\in D(t)$ satisfying $\<x,e_t\>=\<y,e_t\>=1$.  We will
deduce (9.10) from the following lemma, which may be considered
a generalization of the familiar formula
$$
\lim_{n\to\infty}(1+z/n)^n = e^z.
$$
$l^1$ (resp. $l^2$) will denote the space of all sequences of complex
numbers $z=(z(1), z(2),\dots)$ which are summable (resp. square summable).
The norm of $z\in l^2$ is denoted $\|z\|_2$.

\proclaim{Lemma 9.11}  Let $I$ be a directed set and let
$\{z_\alpha: \alpha\in I\}$ be a net of sequences in $l^1\cap l^2$
satisfying
$$
\align
\lim_\alpha\|z_\alpha\|_2 &= 0, {\text and}  \tag{9.11a}\\
\lim_\alpha\sum_{k=1}^\infty z_\alpha(k) &= \zeta\in \Bbb C.  \tag{9.11b}
\endalign
$$
Then for every $\alpha\in I$ the infinite product
$\prod_{k=1}^\infty(1+z_\alpha(k))$ converges absolutely and we have
$$
\lim_\alpha \prod_{k=1}^\infty(1+z_\alpha(k)) = e^\zeta.
$$
\endproclaim

\remark{Remarks}  Notice that since $\sum_{k=1}^\infty|z_\alpha(k)|$ converges
for every $\alpha\in I$, every infinite product
$p_\alpha=\prod_{k=1}^\infty(1+z_\alpha(k))$ converges absolutely as well.

Note too that that Lemma 9.11 reduces to the familiar formula
$$
\lim_{n\to\infty}(1+\zeta/n)^n = e^\zeta,
$$
by taking $I = \{1,2,\dots\}$ and
$$
z_n = (\underbrace{\zeta/n, \zeta/n,\dots,\zeta/n,}_{n{\text\ times}}0, \dots).
$$
In this case (9.11a) follows from the fact that
$\|z_n\|_2 = |\zeta|/{n^{1/2}}\to 0$ as $n\to \infty$,
and (9.11b) is an identity for every $n$.
\endremark

\demo{proof of Lemma 9.11}  Let $\log$ denote the principal branch
of the complex logarithm on the region $\{1+z: |z|<1\}$.  Then for
sufficiently large $\alpha$ we have
$$
\sup_{k\geq 1}|z_\alpha(k)|\leq \|z_\alpha\|_2<1,
$$
and for such an $\alpha$ $\log(1+z_\alpha(k))$ is defined for every
$k=1,2,\dots$.
We will show that the series $\sum_k\log(1+z_\alpha(k))$ is absolutely
convergent for large $\alpha$ and, in fact,
$$
\lim_\alpha \sum_{k=1}^\infty \log(1+z_\alpha(k)) = \zeta.  \tag{9.12}
$$
The required conclusion follows after exponentiating (9.12).

Now since
$$
\lim_{z\to 0}\frac{\log(1+z)-z}{z^2} = 1/2 < 1,
$$
we can find $\epsilon>0$ so that
$$
|\log(1+z)-z|<|z|^2
$$
for all $z$ with $|z|<\epsilon$.  If $\alpha$ is large enough that
$\|z_\alpha\|_2<\epsilon$ then we have
$$
\sum_{k=1}^\infty |\log(1+z_\alpha(k)) - z_\alpha|
\leq \sum_{k=1}^\infty |z_\alpha(k)|^2 < \epsilon^2.
$$
In particular, $\sum_k|\log(1+z_\alpha(k))|<\infty$ because
$z_\alpha\in l^1$.  Now we obtain
(9.12) by noting that for large $\alpha$,
$$
\align
|\sum_k\log(1+z_\alpha(k)) -\zeta| \leq
&\sum_k|\log(1+z_\alpha(k)) - z_\alpha(k)| + |\sum_k z_\alpha(k) - \zeta| \\
\leq &\epsilon^2 + |\sum_k z_\alpha(k) - \zeta|,
\endalign
$$
and using the hypothesis (9.11b) to estimate the second term on
the right.\qed
\enddemo

We apply Lemma 9.11 to the directed set $I$ of all partitions $\Cal P$
of the interval $(0,t]$ as follows.  For every
$$
\Cal P = \{0=t_0 < t_1 < \dots < t_n=t\}
$$
we have from (9.8)
$$
B_{\Cal P}(t;x,y) = \sum_{k=1}^n z_{\Cal P}(k)
$$
where $z_{\Cal P}(k) = \<x(t_{k-1},t_k),y(t_{k-1},t_k)\>-1$.  Now for
$0\leq a<b\leq t$ we have
$$
\align
|\<x(a,b),y(a,b)\>-1| &= |\<x(a,b)-e(a,b),y(a,b)-e(a,b)\>|\\
&\leq \|x(a,b)-e(a,b)\|\cdot\|y(a,b)-e(a,b)\|.
\endalign
$$
Let $\phi,\psi:[0,t]\to \Bbb R$ be the functions $\phi(0)=\psi(0)=1$,
$\phi(s)=\|x_s\|^2$, $\psi(s) = \|y_s\|^2$ for $0<s\leq 1$.
Then $\phi$ and $\psi$
are continuous monotone increasing, and an estimate like those in the proof
of Theorem 9.8 shows that
$$
\|x(a,b)-e(a,b)\|^2 = \|x(a,b)\|^2-1 \leq \phi(b) - \phi(a)
$$
and
$$
\|y(a,b)-e(a,b)\|^2 = \|y(a,b)\|^2-1 \leq \psi(b) - \psi(a).
$$
Thus we can estimate $|z_{\Cal P}(k)|$ as follows.  Setting
$$
\epsilon(\alpha) = \sup_{|t-s|\leq \alpha} |\phi(t) - \phi(s)|
$$
for positive $\alpha$, we have
$$
|z_{\Cal P}(k)|^2 \leq (\phi(t_k)-\phi(t_{k-1}))(\psi(t_k)-\psi(t_{k-1}))
\leq \epsilon(|\Cal P|) (\psi(t_k)-\psi(t_{k-1})),
$$
 $|\Cal P|$ denotinging the norm of $\Cal P$.  It follows that
$$
\sum_{k=1}^n|z_{\Cal P}(k)|^2 \leq \epsilon(|\Cal P|)(\phi(t) - 1).
$$
Since $\lim_{\Cal P}\epsilon(|\Cal P|) = 0$ because
$\lim_{\Cal P}|\Cal P| = 0$ and $\phi$ is uniformly continuous,
we conclude that $\lim_{\Cal P}\|z_{\Cal P}\|_2 = 0$.  This establishes
condition (9.11a). In this case (9.11b) is the formula
$$
\lim_{\Cal P}B_{\Cal P}(t;x,y) = B_\infty(t;x,y).
$$
Since $1+z_{\Cal P}(k)=\<x(t_{k-1},t_k),y(t_{k-1},t_k)\>$
we may conclude from Lemma 9.11 that
$$
\lim_{\Cal P}\prod_{I\in \Cal P}\<x_I,y_I\> = e^{B_\infty(t;x,y)}.
$$
But for every partition $\Cal P$ we have
$$
\<x,y\> = \prod_{I\in \Cal P}\<x_I,y_I\>,
$$
and thus the preceding discussion implies that
$$
\<x,y\> = e^{B_\infty(t;x,y)}
$$
which is the required formula (9.10). \qed
\enddemo

\remark{Remark 9.13}
The argument we have given yields a somewhat stronger
result involving {\it sequential} convergence.  For each
$t>0$, let $\Cal P_{1,t}\leq \Cal P_{2,t} \leq\dots$ be
any increasing sequence of partitions of the interval
$[0,t]$ with the property
$$
\lim_{n\to\infty}|\Cal P_{n,t}| = 0,
$$
for every $t$.  Then for every $x,y\in D(t)$ satisfying
$\<x, e_t\> = \<y,e_t\> = 1$, we claim that
$$
L^e(t;x,y) = \lim_{n\to \infty}\sum_{I\in \Cal P_{n,t}}
(\<x_I,y_I\> - 1).\tag{9.14}
$$
To see that, fix $x,y$ and define $B_n^\prime(t;x,y)$ by
$$
B_n^\prime(t;x,y) = \sum_{I\in\Cal P_{n,t}}(\<x_I,y_I\> -1).
$$
(9.5) and (9.7) imply that the limit
$$
B_\infty^\prime(t;x,y) = \lim_{n\to\infty}B_n^\prime(t;x,y)
$$
exists.  Moreover, we may apply Lemma 9.11 as in the preceding
argument to conclude that
$$
\<x,y\> = e^{B_\infty^\prime(t;x,y)}.
$$

Finally, the estimates of Proposition 9.8 are valid for
$B_\infty^\prime(t;x,y)$ as well as for $B_\infty(t;x,y)$.
Thus for any pair of left-decomposable sections
$u,v$ satisfying $\<u_s,e_s\> = \<v_s,e_s\> = 1$, the function
$$
t\in (0,\infty)\mapsto B_\infty^\prime(t;u_t,v_t)
$$
is continuous and tends to $0$ as $t\to 0+$.  It follows that
$$
B_\infty^\prime(t;u_t,v_t) = B_\infty(t;u_t,v_t) = L^e(t;u_t,v_t)
$$
for all $t$.  Fixing $t$ and choosing $u, v$ appropriately, we
obtain (9.14).
\endremark

We will make use of formula (9.14) in the sequel.

\subheading{10. Existence of measurable propagators}

The reader may have noticed that none of the results of
sections 5 through 9 made any reference to measurability
or to measurable sections, even though there is a natural
Borel structure on the product system $E$.  Measurability
was simply not an issue in those matters.  On the other hand,
in the applications that will be discussed in part III it
will be necessary do deal with measurable elements
of $D$, and with a measurable reference section $e$; this is
necessary in order to satisfy the measurability hypothesis
of Definition 2.2.  The existence of sufficiently many
measurable elements of $D$ is established in the following
result.

\proclaim{Theorem 10.1}  Let $t_0>0$ and let $u\in E(t_0)$ be a nonzero
decomposable vector.  Then there is a left-coherent
decomposable Borel section $t\in (0,\infty)\to x_t\in D(t)$
such that $x_{t_0}$ is a scalar multiple of $u$.
\endproclaim
\demo{proof} We may assume without loss of generality that $\|u\|=1$.
Since $u$ is decomposable, we can find elements
$a_s$, $u_s\in E(s)$ such that
$$
u = a_su_{t_0-s},\qquad 0<s<t_0; \tag{10.2}
$$
and by renormalizing again if necessary we can arrange that
$\|a_s\|=\|u_s\|=1$.  Set $a_{t_0}=u$.  Then $\{a_s: 0<s\leq t_0\}$
is a left-decomposable family of unit vectors.
The following result implies that we can
choose $a$ in a measurable way.

\proclaim{Proposition 10.3}Let $t\in(0,t_0]\mapsto a_t\in E(t)$ be
any left-coherent family of unit vectors.  Then there is a function
$t\in(0,t_0]\mapsto\lambda_t\in\Bbb C$ such that
$|\lambda_t|=1$ for every $t$ and $t\mapsto\lambda_ta_t$
is a measurable section of $E\restriction_{(0,t_0]}$.
\endproclaim

Let us assume, for the moment, that the technical result 10.3 has
been established.  If we replace $a_s$ with $\lambda_sa_s$ in
(10.2) then we may assume that $a_s$ is measurable in $s$.  We
can now define a measurable section $t\in(0,\infty)\mapsto x_t\in E(t)$ by
$$
x_t =
\cases
\<a_t,e_t\>^{-1}a_t,&\text{for $0<t\leq t_0$}\\
a_{t_0}e(t_0,t),&\text{for $t>t_0$},
\endcases
$$
where $\{e(s,t): 0<s<t<\infty\}$ is the propagator associated
with $e$.  It is clear that $x\in D^e$ and that $x_{t_0}$ is a
scalar multiple of $u$, completing the proof of Theorem 10.1.
\enddemo

\demo{proof of 10.3}  Consider the family of operators
$\{P_s: 0<s<t_0\}\subseteq\Cal B(E(t_0))$ defined by
$$
P_s(\xi) = a_s\cdot L_{a_s}^*(\xi),\qquad \xi\in E(t_0)
$$
(see the proof of Theorem 2.1 for a discussion of operators of
the form $L^*_z, z\in E$).
Notice that since $E(t_0)$ is spanned by $E(s)E(t_0-s)$, $P_s$ is
uniquely determined by its action on decomposable vectors
$\xi_1\xi_2$, $\xi_1\in E(s)$, $\xi_2\in E(t_0-s)$:
$$
P_s(\xi_1\xi_2) = \<\xi_1,a_s\>a_s\xi_2.
$$
This also shows that if we identify $E(t_0)$ with $E(s)\otimes E(t_0-s)$
then $P_s$ becomes
$$
[a_s]\otimes {\bold 1}_{t_0-s},
$$
$[a_s]$ denoting the projection onto the one-dimensional
subspace $\Bbb C\cdot a_s$ of $E(s)$ and
${\bold 1}_r$ denoting the identity operator in $\Cal B(E(r))$.
Now since $\{a_s: 0<s<t_0\}$ is left-coherent, the family of projections
$\{P_s: 0<s<t_0\}$ satisfies
$$
s_1 < s_2 \implies P_{s_2} \leq P_{s_1}
$$
(see the proof of Theorem 2.1).  Thus for any fixed $\xi\in E(t_0)$,
$$
s\mapsto \<P_s\xi,\xi\>
$$
is a monotone decreasing function of $(0,t_0)$, and hence measurable.
By polarization it follows that $\<P_s\xi,\eta\>$ is measurable in
$s$ for any $\xi,\eta\in E(t_0)$.  Hence
$$
s\in (0,t_0) \mapsto P_s\in \Cal B(E(t_0))
$$
is a measurable projection-valued operator function.

For each $0<s<t_0$, let $Q_s$ denote the rank-on projection
$Q_s = [a_s] \in \Cal B(E(s))$.  We have just seen that
$$
s\mapsto Q_s\otimes{\bold 1}_{t_0-s} \tag{10.5}
$$
is measurable and we claim now that $Q$ itself is measurable.
Equivalently, we claim that for any pair of measurable sections
$\xi_s, \eta_s\in E(s)$, the complex valued function
$s\in(0,t_0)\mapsto \<Q_s\xi_s,\eta_s\>$ is measurable.  To see that we
choose any measurable section $t\mapsto u_t\in E(t)$ of unit vectors
and write
$$
\<Q_s\xi_s,\eta_s\> = \<Q_s\xi_s,\eta_s\>\<u_{t_0-s},u_{t_0-s}\> =
\<Q_s\otimes{\bold 1}_{t_0-s}(\xi_su_{t_0-s}),\eta_su_{t_0-s}\>,
$$
for $0<s<t_0$.  The right side is measurable in $s$ by (10.6) and
the fact that both $\xi_su_{t_0-s}$ and $\eta_su_{t_0-s}$ are
measurable functions of $s$.

Finally, since $s\in(0,t_0)\mapsto Q_s\in \Cal B(E(s))$ is
a measurable family we claim that there is a measurable section
$s\in(0,t_0)\mapsto b_s\in E(s)$ satisfying $\|b_s\| = 1$ for
every $s$ and
$$
Q_sb_s = b_s, \qquad 0<s<t_0.  \tag{10.6}
$$
To prove (10.6) we choose a measurable basis for the family of
Hilbert spaces $p:E\to (0,\infty)$.  That is, we find a sequence of
measurable sections $e_n:t\in(0,\infty)\mapsto e_n(t)\in E(t)$,
$n=1,2,\dots$ such that $\{e_1(t), e_2(t),\dots \}$ is an orthonormal
basis for $E(t)$ for every $t>0$.  This is possible because of the
last axiom for product systems \cite{2, (1.8) {\it iii}}.  For
every $t>0$ we must have $Q_te_n(t)\neq 0$ for some $n$, and we
define $n(t)$ to be the smallest such $n$.  Now for every positive
integer $k$ we have
$$
\{t\in (0,\infty): n(t)> k\} = \bigcap_{i=1}^{k}\{t\in (0,\infty):
\<Q_te_i(t),e_i(t)\> = 0\},
$$
and the right side is a Borel set since each function
$t\mapsto \<Q_te_i(t),e_i(t)\>$ is measurable.  It follows that
the function $t\in (0,\infty)\mapsto n(t)\in \Bbb R$ is measurable.
Hence $\xi_t = e_{n(t)}$ defines a measurable section of $E$
having the property that $Q_t\xi_t\neq 0$ for every $t>0$.  Thus we
obtain a section $b$ as required by (10.6) by setting
$$
b_t = \|Q_t\xi_t\|^{-1}Q_t\xi_t, \qquad 0<t<t_0.
$$
Now since $Q_t$ is the projection onto the one-dimensional space
$\Bbb Ca_t$ it follows that there is a complex number $\lambda_t$
such that $b_t = \lambda_ta_t$.  Since $\|a_t\| = \|b_t\| = 1$ we
have $|\lambda_t|=1$, completing the proof of Proposition 10.3. \qed
\enddemo

\heading{Part III.  Applications}
\endheading
In the following two sections we apply the preceding results
to classify certain product systems and certain
\esg s.

\subheading{11. Decomposable continuous tensor products}
A product system $p:E\to (0,\infty)$ is said to be
{\it decomposable} if for every $t>0$, $E(t)$ is the
closed linear span of the set $D(t)$ of all its decomposable
vectors.  It is easy to see that if this condition is satisfied
for a single $t_0>0$ then it is satisfied for every
$t>0$.

\proclaim{Theorem 11.1}
A decomposable product system is either isomorphic
to the trivial product system
with one-dimensional spaces $E(t)$,
$t>0$, or it is isomorphic to one of the standard
product systems $E_n$, $n=1,2,\dots,\infty$.
\endproclaim

\remark{Remarks}
We recall that any product system with one-dimensional
fibers $E(t)$ for every $t>0$ is isomorphic to the
trivial product system $p:Z\to (0,\infty)$, where
$Z=(0,\infty)\times\Bbb C$ with multiplication
$(s,z)(t,w) = (s+t,zw)$, with the usual inner product
on $\Bbb C$, and with projection $p(t,z)=t$ \cite{6,
Corollary of Prop. 2.3} and \cite{7}.
\endremark

\demo{proof of Theorem 11.1}

By Theorem 10.3, we can find a Borel section
$t\in (0,\infty)\mapsto e(t)\in D(t)$.  By replacing $e_t$
with $e_t/\|e_t\|$ if necessary we may assume that
$\|e_t\| = 1$ for every $t$.  Let $L^e:\Delta^2\to \Bbb C$
be the function provided by Theorem 8.2 which satisfies
$$
e^{L^e(x,y)} = \frac{\<x,y\>}{\<x,e_t\>\<e_t,y\>}\tag{11.2}
$$
for every $x,y\in D(t)$, $t>0$.  We will show first that
$(\Delta,L^e)$ is a metric path space; that is, $L^e$
is an additive form on $\Delta^2$.

Indeed, most of that assertion follows immediately from
the results of Parts I and II.  Theorem 9.2 implies that
$L^e$ restricts to a positive definite function on
$\Delta(t)\times \Delta(t)$ for every $t>0$, and because of
Propositions 9.5 and 9.7 taken together with with formula
9.14, $L^e$ must satisfy the separability condition of
Definition 2.6.  Proposition 8.8 shows that $L^e$ is
additive with defect function $\psi:\Delta\times\Delta\to \Bbb C$
of the form $\psi(x,y)=\psi_s(t;y)$ for every $x\in D(s)$,
$y\in D(t)$, $s,t>0$.  Thus we need only
establish the measurability criterion of Definition 2.2.

Notice that Definition 2.2 makes reference to the propagator
$\{\dot z(r,s): 0\leq r<s\leq t\}$ associated with an element
$\dot z\in \Delta(t)$.  While this propagator in the path
space $\Delta$ is uniquely determined by $\dot z$, an element
$z\in D(t)$ does not determine a unique propagator
$\{z(r,s)\in D(s-r): 0\leq r<s\leq t\}$ because of the
failure of unique factorization in $\{D(t): t>0\}$.
Nevertheless, any propagator in $\Delta$ can be lifted
to a {\it measurable} propagator in $\{D(t): t>0\}$.
More precisely, given any element $z\in D(t)$ then
Theorem 10.1 provides a left-coherent {\it measurable} family
$\{z_s\in D(s): 0<s\leq t\}$ for which $z_t$ is a scalar
multiple of $z$.  Recalling that a left-coherent family
gives rise to a unique propagator because of the left
cancellation law in $\{D(t): t>0\}$, we conclude that
the propagator $\{z(r,s)\in D(s-r):0\leq r<s\leq t\}$
associated with $\{z_s: 0<s\leq t\}$ projects to the required
propagator in $\Delta$, viz
$$
\dot z_s=\dot z_r\dot z(r,s) ,\qquad 0<r<s\leq t,
$$
and of course $\dot z(0,s)=\dot z_s$ for $0<s\leq t$.

In order to establish the measurability criterion of
Definition 2.2, choose $T_1,T_2$ satisfying $0<T_1<T_2$  and
choose $x\in D(T_1)$, $y\in D(T_2)$.  By the preceding
remarks, we may find a measurable left-coherent section
$y_s\in D(s)$, $0<s\leq T_2$ whose propagator projects to
the propagator of $\dot y\in \Delta(T_2)$.
It suffices to show that
the function
$$
\lambda\in (0,T_2-T_1)\mapsto L^e(T_1;x,y(\lambda,\lambda+T_1))\tag{11.3}
$$
is a complex-valued Borel function.

To prove (11.3) we will make use of
the fact that $L^e$ is the limit of a {\it sequence} of functions
for which the fact of measurability is obvious.
Let $\Cal P_1, \Cal P_2,\dots$
be a sequence of finite partitions of the interval $[0,T_1]$
such that $\Cal P_{n+1}$ is a refinement of $\Cal P_n$ and such
that the norms $|\Cal P_n|$ tend to zero as $n\to \infty$.
By (9.14) we may conclude that for every $u,v\in D(T_1)$
satisfying $\<u,e_{T_1}\>=\<v,e_{T_1}\>=1$, we have
$$
L^e(T_1;u,v) = \lim_{n\to \infty}\sum_{I\in \Cal P_n}(\<u_I,v_I\>-1).
$$
Taking
$$
\align
u=&\<x,e_{T_1}\>^{-1}x \qquad\text{and} \\
v=&v_\lambda=\<y(\lambda,\lambda+T_1),e_{T_1}\>^{-1}y(\lambda,\lambda+T_1),
\endalign
$$
we see that
the left side of (11.3) is exhibited as the limit of a convergent
sequence of functions
$$
f_n(\lambda) = \sum_{I\in \Cal P_n}(\<u_I,(v_\lambda)_I\>-1).
$$
Thus it suffices to show
that each $f_n$ is Borel-measurable on $0<\lambda<T_2-T_1$.
In order to see that, choose an interval $I=(a,b]\subseteq (0,T_1]$
for which $a<b$ and look at the inner product $\<u_I,(v_\lambda)_I\>$.
Noting that
$$
\align
y(\lambda,\lambda+T_1)=&
y(\lambda,\lambda+a)y(\lambda+a,\lambda+b)y(\lambda+b,\lambda+T_1)\\
e_{T_1} =& e_ae(a,b)e(b,T_1),
\endalign
$$
we can write down an obvious propagator for $v_\lambda$ and we find
that
$$
v_\lambda(a,b) = \<y(\lambda+a,\lambda+b),e(a,b)\>^{-1}y(\lambda+a,\lambda+b).
$$
It follows that $f_n(\lambda)$ is a finite linear combination of functions
of the form
$$
\lambda\in (0,T_2-T_1)\mapsto
\frac{\<u(a,b),y(\lambda+a,\lambda+b)\>}
{\<e(a,b),y(\lambda+a,\lambda+b))\>}-1.  \tag{11.4}
$$
Now for any element $w\in D(b-a)$ we have
$$
\align
\<y_{\lambda+a}w,y_{\lambda+b}\>&=
\<y_{\lambda+a}w,y_{\lambda+a}y(\lambda+a,\lambda+b)\>\\
&=\|y_{\lambda+a}\|^2\<w,y(\lambda+a,\lambda+b)\>.
\endalign
$$
Thus the right side of (11.4) can be written
$$
\frac{\<y_{\lambda+a}u(a,b),y_{\lambda+b}\>}
{\<y_{\lambda+a}e(a,b),y_{\lambda+b}\>} - 1.
$$
This is obviously a measurable function of $\lambda$
because $s\mapsto y_s$ is a measurable section, left
multiplication by a fixed element of $E$ is a measurable
mapping of $E$ into itself, and because the inner product
$\<\cdot,\cdot\>:E^2\to \Bbb C$ is measurable

Thus we have established the fact that $(\Delta,L^e)$
is a metric path space.  By Theorem 4.3 there is a
separable Hilbert space $\Cal C$, a function
$\rho:\Delta\to \Bbb C$, and a logarithm mapping
$$
\log:\Delta\to \Cal P_{\Cal C}
$$
such that $\log(xy) = \log(x)\boxplus\log(y)$ for every
$x,y\in \Delta$ and
$$
L^e(t;x_1,x_2)=
\<\log(x_1),\log(x_2)\>+\rho(x_1)+\overline{\rho(x_2)}\tag{11.5}
$$
for all $x_1,x_2\in \Delta(t)$, $t>0$.  In fact, in order
to obtain the necessary measurability properties, we must
use the specific function $\rho$ defined in the proof
of Theorem 4.3.

Now when the path space $\Cal P_{\Cal C}$ is exponentiated,
it gives rise to the standard product system $E_\Cal C$.
In more detail, consider the symmetric Fock space $H_{\Cal C}$ over the
one-particle space $L^2((0,\infty);\Cal C)$, and consider
the exponential map
$$
\exp:L^2((0,\infty);\Cal C)\to H_{\Cal C}
$$
defined by
$$
\exp(f) = \sum_{n=0}^\infty \frac{1}{\sqrt{n!}}f^{\otimes n}.
$$
For every $t>0$ we define
$$
E_{\Cal C}(t) = \overline{\text{span}}\{\exp(f): f\in \Cal P_{\Cal C}(t)\}.
$$
$E_{\Cal C}$ is the total space of this family of Hilbert spaces,
with multiplication
$\xi\in E_\Cal C(s),\eta\in E_{\Cal C}(t)\mapsto
\xi\eta\in E_{\Cal C}(s+t)$
defined uniquely by requiring that the generating vectors should
multiply thus:
$$
\exp(f)\exp(g) = \exp(f\boxplus g),
\qquad f\in \Cal P_{\Cal C}(s), g\in \Cal P_{\Cal C}(t).
$$

We will use (11.5) to construct an isomorphism of product
systems $W:E\to E_{\Cal C}$.  But in order to define
$W$ it is necessary to solve another cohomological problem.
The result is summarized as follows.  We write
$$
p:D = \{(t,x): x\in D(t), t>0\}\to (0,\infty)
$$
for the fiber space determined by the family of sets
$D(t),t>0$ with projection $p(t,x)=t$.  Notice that $D$
carries a natural Borel structure as a subspace of $E$.

\proclaim{Theorem 11.6}
There is a Borel-measurable function $u:(0,\infty)\to \Bbb C$
satisfying $|u(t)|=1$ for every $t>0$, such that the function
$f:D\to \Bbb C$ defined by
$$
f(x) = u(t)\<x,e_t\>e^{\rho(x)},\qquad x\in D(t),\quad t>0
$$
is multiplicative: $f(xy)=f(x)f(y)$, $x\in D(s), y\in D(t)$.
\endproclaim
\demo{proof of Theorem 11.6}
Let $f_0(x) = \<x,e_t\>e^{\rho(x)}$.  Notice first that
$f_0$ is Borel-measurable.  Indeed, recalling the formula
for $\rho(x)$ given in the proof of Theorem 4.3, we have
for every $x\in D(s)$ and $s>0$,
$$
\rho(x) = \<[x]-[e_s],\phi_s\>+L^e(s;x,e_s)-
\frac{1}{2}(L^e(s;e_s,e_s)+\|\phi_s\|^2).  \tag{11.7}
$$
Because of the representation of $L^e$ as a sequential
limit in (9.14) we see that $L^e:D^2\to \Bbb C$ is a Borel
function, and it follows that
$x\in D(s)\mapsto [x]-[e_s]$ defines a Borel
map of $D$ into the Hilbert space $H_\infty$.  Finally,
since $s\in (0,\infty)\mapsto \phi_s$ is a Borel function
taking values in the space of all locally square integrable
$\Cal C$-valued functions on $(0,\infty)$,
we see that the right side of (11.7)
defines a Borel function on the subspace $D$ of $E$.  Since
the inner product is a Borel function on $E$, it follows
that
$$
f_0(x) = \<x,e_s\>e^{\rho(x)},\qquad x\in D(s)
$$
defines a complex-valued Borel function on $D$.

We consider the associated ``coboundary"
$c:D\times D\to \Bbb C$:
$$
c(x,y)= \frac{f_0(xy)}{f_0(x)f_0(y)}.
$$
Note that $c$ is a Borel function as well, since the
multiplication operation of $E$ is Borel measurable.
We will show that for every $s,t>0$, $x_i\in D(s)$,
$y_i\in D(t)$,
$$
\align
|c(x_1,y_1)|=&1,\tag{11.8.1}\\
c(x_1,y_1)=&c(x_2,y_2).\tag{11.8.2}
\endalign
$$
Assume for the moment that the equations (11.8) have been
established.  It follows that there is a function
$$
c_0:(0,\infty)\times(0,\infty)\to \{z\in \Bbb C: |z|=1\}
$$
such that for all $x\in D(s)$, $y\in D(t)$ we have
$$
c(x,y) = c_0(s,t).
$$
$c_0$ is clearly measurable because $c_0(s,t) = c(e_s,e_t)$,
and $r\mapsto e_r\in D(r)$ is a measurable section.
We will then show that $c_0$ satisfies the multiplier
equation,
$$
c_0(r,s+t)c_0(s,t)=c_0(r+s,t)c_0(r,s),\tag{11.9}
$$
for every $r,s,t>0$.
By \cite{6, Corollary of Propostion 2.3}, there is
a Borel function $u:(0,\infty)\to \{z\in \Bbb C: |z|=1\}$
such that
$$
c_0(s,t) = \frac{u(s)u(t)}{u(s+t)},\qquad s,t>0.
$$
Once we have $u$ it is clear that the function
$f(x)=u(t)f_0(x)$, $x\in D(t)$, $t>0$ satisfies
$$
f(xy)=f(x)f(y).
$$

Thus, we must prove the formulas (11.8) and (11.9).
For $t>0$, consider the function $W_t:D(t)\to E_{\Cal C}(t)$
defined by
$$
W_t(x) = f_0(x)\exp\log(x) = \<x,e_t\>e^{\rho(x)}\exp\log(x).
$$
We claim that for $x_1,x_2\in D(t)$ we have
$$
\<W_t(x_1),W_t(x_2)\>=\<x_1,x_2\>.  \tag{11.10}
$$
Indeed, using Theorem 4.3 we have
$$
\align
\<\exp\log(x_1),\exp\log(x_2)\> =& e^{\<\log(x_1),\log(x_2)\>}
=\exp(L^e(t;x_1,x_2)-\rho(x_1)-\overline{\rho(x_2)})\\
=&\frac{\<x_1,x_2\>}{\<x_1,e_t\>\<e_t,x_2\>}
e^{-(\rho(x_1)+\overline{\rho(x_2)})}.
\endalign
$$
Thus the left side of (11.10) is
$$
f_0(x_1)\overline{f_0(x_2)}\<\exp\log(x_1),\exp\log(x_2)\>=
\<x_1,e_t\>\overline{\<x_2,e_t\>}\frac{\<x_1,x_2\>}{\<x_1,e_t\>\<e_t,x_2\>}=
\<x_1,x_2\>,
$$
as asserted.

The family $\{W_t: t>0\}$ obeys the following
multiplicative rule.  For $x\in D(s)$, $y\in D(t)$, we claim
$$
W_{s+t}(xy) = c(x,y)W_s(x)W_t(y).  \tag{11.11}
$$
Notice that the multiplication on the left side of (11.11) is
performed in $E$ and on the right side it is
performed in $E_{\Cal C}$.  To check (11.11), recall that the
multiplication in $E_{\Cal C}$ is related to the
operation $\boxplus$ by the following formula: if
$f\in E_{\Cal C}(s)$ and $g\in E_{\Cal C}(t)$ then
$$
\exp(x)\exp(g) = \exp(f\boxplus g).
$$
Thus, using the additivity property of the logarithm
mapping the left side of (11.11) can be rewritten
$$
\align
f_0(xy) \exp\log(xy) =& f_0(xy)\exp(\log(x)\boxplus\log(y))
= f_0(xy)\exp(\log(x))\exp(\log(y)) \\
=&\frac{f_0(xy)}{f_0(x)f_0(y)}W_s(x)W_t(y)=c(x,y)W_s(x)W_t(y).
\endalign
$$

We claim next that for $x_i\in D(s)$, $y_i\in D(t)$ we have
$$
c(x_1,y_1)\overline{c(x_2,y_2)}=1.  \tag{11.12}
$$
To see this, choose $x_i\in D(s)$, $y_i\in D(t)$.  Noting
that by (11.10)
$$
\align
\<x_1y_1,x_2y_2\>=\<x_1,x_2\>\<y_1,y_2\>=&
\<W_s(x_1),W_s(x_2)\>\<W_t(y_1),W_t(y_2)\>\\
=&\<W_s(x_1)W_t(y_1),W_s(x_2)W_t(y_2)\>,
\endalign
$$
we have
$$
\align
&c(x_1,y_1)\overline{c(x_2,y_2)}\<x_1y_1,x_2y_2\>=\\
&c(x_1,y_1)\overline{c(x_2,y_2)}\<W_s(x_1)W_t(y_1),W_s(x_2)W_t(y_2)\>=\\
&\<c(x_1,y_1)W_s(x_1)W_t(y_1),c(x_2,y_2)W_s(x_2)W_t(y_2)\>=\\
&\<W_{s+t}(x_1y_1),W_{s+t}(x_2y_2)\> =
\<x_1y_1,x_2y_2\>.
\endalign
$$
The claim follows after cancelling $\<x_1y_1,x_2y_2\>\neq 0$.

Set $x_2=x_1,y_2=y_1$ in (11.12) to obtain $|c(x_1,y_1)|=1$.
Thus if we multiply through in (11.12) by $c(x_2,y_2)$ we
obtain
$$
c(x_1,y_1)=c(x_2,y_2),
$$
hence (11.8.1) and (11.8.2) are established.

Now we can define a function
$c_0:(0,\infty)\times(0,\infty)\to \Bbb C$ by
$c_0(s,t)=c(e_s,e_t)$.  $c_0$ is a Borel function which, because
of (11.11), obeys
$$
W_{s+t}(xy) = c_0(s,t)W_s(x)W_t(y),\qquad x\in D(s), y\in D(t).
$$
The latter formula implies that $c_0$ must satisfy the multiplier
equation (11.9).  Indeed, for $r,s,t>0$ and
$x\in D(r),y\in D(s),z\in D(t)$ we have
$$
W_{r+s+t}(x(yz))=c_0(r,s+t)W_r(x)W_{s+t}(yz)=
c_0(r,s+t)c_0(s,t)W_r(x)W_s(y)W_t(z),
$$
while
$$
W_{r+s+t}((xy)z)=c_0(r+s,t)W_{r+s}(xy)W_{t}(z)=
c_0(r+s,t)c_0(r,s)W_r(x)W_s(y)W_t(z).
$$
Since $W_r(x)W_s(y)W_t(z)\neq 0$, equation (11.9) follows.

The argument given above can now be applied to complete
the proof of Theorem 11.6\qed
\enddemo

We are now in position to write down an isomorphism
$E\cong E_{\Cal C}$.  Choose $u$ as in Theorem 11.6.
For every $t>0$, $x\in D(t)$ define
$W_t(x)\in E_{\Cal C}(t)$ by
$$
W_t(x)= u(t)\<x,e_t\>e^{\rho(x)}\exp(\log(x)).
$$
Equation (11.10) implies that
$$
\<W_t(x_1),W_t(x_2)\>=\<x_1,x_2\>.
$$
Since $E(t)$ is spanned by $D(t)$, $W_t$ can be extended
uniquely to a linear isometry of $E(t)$ into $E_{\Cal C}(t)$,
and we will denote the extended mapping by the same letter $W_t$.
The range of $W_t$ is given by
$$
W_t(E(t)) = \overline{\text{span}}[\exp\log(D(t))]=E_{\Cal C}(t),
$$
because the set $\log(D(t))$ is a strongly spanning subset
of $\Cal P_{\Cal C}(t)$.  Thus the total map
$$
W: E\to E_{\Cal C}
$$
is an isomorphism of families of Hilbert spaces.

Because of the way we chose the function $u$, the multiplication
formula (11.11) simplifies to
$$
W_{s+t}(xy)=W_s(x)W_t(y),\qquad x\in D(s), y\in D(t).
$$
Using bilinearity and the fact that $D(r)$ spans $E(r)$ for
every $r>0$, the latter implies that $W$ is a homomorphism of
product structures in that
$$
W_{s+t}(\xi\eta)=W_s(\xi)W_t(\eta),\qquad \xi\in E(s),y\in E(t).
$$

In particular, $W$ is a bijection of the standard Borel space
$E$ onto the standard Borel space $E_{\Cal C}$.  Thus to see that
$W$ is a Borel isomorphism (and therefore an isomorphism of product
systems), it suffices to show that it is measurable.  The proof
of that is a routine variation on the argument
presented in detail in \cite{2, pp 55--57}, and we omit it.

Finally, notice that it is possible that $\Cal C$ is the trivial
Hilbert space $\{0\}$.  However, in this case $E_\Cal C$ is the
trivial product system with one-dimensional fibers.  If
$\Cal C\neq \{0\}$ and we let $n$ be the dimension of $\Cal C$,
then $n=1,2,\dots,\aleph_0$ and $E_\Cal C$ is the standard
product system $E_n$\qed
\enddemo

\subheading{12. Decomposable $E_0$-semigroups}
Let $\alpha=\{\alpha_t; t\geq 0\}$ be an \esg, and
for every $t>0$ let $\Cal E(t)$ be the operator space
$$
\Cal E(t) = \{T\in \Cal B(H): \alpha_t(A)T=TA, \forall A\in \Cal B(H)\}.
$$
$\Cal E(t)$ is a Hilbert space relative to the inner
product defined on it by
$$
T^*S=\<S,T\>{\bold 1}.
$$
The family of Hilbert spaces
$$
p: \Cal E=\{(t,T): t>0, T\in \Cal E(t)\}\to (0,\infty)
$$
with projection $p(t,T)=t$
is actually a product system with respect to operator multiplication
$(s,S)(t,T)=(s+t,ST)$ \cite{2}.  In particular,
$\Cal E(s+t)$ is the norm closed linear span of the
set of all products $\{ST: S\in \Cal E(s), T\in \Cal E(t)\}$.

As with any product system, it makes sense to speak of
decomposable elements of $\Cal E(t)$; thus,
an operator $T\in \Cal E(t)$ is decomposable if,
for every $0<s<t$, T admits a factorization
$T=AB$ where $A\in \Cal E(s)$ and $B\in \Cal E(t-s)$.  The
set of all decomposable operators in $\Cal E(t)$ is denoted
$\Cal D(t)$.  Let $H(t) = [\Cal D(t)H]$
be the closed linear span of the ranges of all operators in
$\Cal D(t)$.  The spaces $H(t)$ are obviously decreasing with
$t$.

\proclaim{Definition 12.1}
$\alpha=\{\alpha_t: t\geq 0\}$ is called decomposable if
$[\Cal D(t)H]=H$ for some (and therefore every) $t>0$.
\endproclaim

The following simple result shows that this terminology
is consistent with the notion of decomposability for
product systems.

\proclaim{Proposition 12.2}
An \esg\ $\alpha$ is decomposable iff its associated product
system is decomposable in the sense of section 11.
\endproclaim
\demo{proof}
Suppose first that $\Cal E$ is a decomposable product system.
Then for every $t>0$ $\Cal E(t)$ is the norm-closed linear span
of $\Cal D(t)$.  Since we have $[\Cal E(t)H]=H$ for the product
system of an arbitrary \esg, it follows that $H=[\Cal D(t)H]$ as
well.

Conversely, assuming that $\alpha$ satisfies Definition 12.1 we
pick $t>0$ and an operator $T\in \Cal E(t)$ such that $T$ is
orthogonal to $\Cal D(t)$.  Because of the definition of the
inner product in $E(t)$ it follows that for every $S\in \Cal D(t)$
we have
$$
T^*S = \<S,T\>\bold{1} = 0,
$$
hence $TH$ is orthogonal to $[\Cal D(t)H]=H$, hence $T=0$.
Thus $\Cal D(t)$ spans $\Cal E(t)$\qed
\enddemo

We can immediately deduce from Theorem 11.1 that the product
system of a decomposable \esg\ is either the trivial product
system or it is isomorphic to one of the standard product
systems $E_n$, $n=1,2,\dots,\infty$.  Since the product system
of an \esg\ is a complete invariant for cocycle conjugacy and
since the standard product systems $E_n$ are associated with $CCR$
flows (or $CAR$ flows) \cite{2}, we can infer
the following classification of \esg s as an immediate
consequense of Theorem 11.1.

\proclaim{Theorem 12.3}
Let $\alpha$ be a decomposable \esg\ acting on $\Cal B(H)$
which is nontrivial in the sense that
it cannot be extended to a group of automorphisms of
$\Cal B(H)$.  Then $\alpha$ is cocycle conjugate to a $CCR$
flow.
\endproclaim

\remark{Concluding remarks}
Note that Theorem 12.3 implies that every decomposable
\esg\ $\alpha$ has plenty of intertwining semigroups...that
is, semigroups of isometries $U=\{U_t: t\geq0\}$ acting on
$H$ for which
$$
\alpha_t(A)U_t=U_tA,\qquad A\in \Cal B(H), t>0.
$$
Indeed, a decomposable \esg\ must be completely spatial.

In a more philosophical vein, we conclude
that {\it any construction of
\esg s which starts with a path space cannot
produce anything other than $CCR$ flows and
their cocycle perturbations}.  For example, in
\cite{2 pp 14--16} we gave examples of product systems
using Gaussian random processes and Poisson
processes.  The latter examples did not appear to
contain enough units to be standard ones.  However,
a closer analysis showed that there were ``hidden"
units, and in fact there were enough of them so that
these product systems were indeed standard.  Theorems 11.1
and 12.4 serve to clarify this phenomenon because these
examples constructed from random processes were clearly
decomposable.

We believe that Theorem 12.4 is analogous to
the familiar description of representations
of the compact operators (i.e., every representation
is unitarily equivalent to a multiple of the identity
representation), or to the Stone-von Neumann theorem.
Certainly, the conclusion of 12.4 implies that
decomposable \esg s exhibit ``type I" behavior.
What is interesting here is that the \esg s whose
product systems are isomorphic to a given one $E$ correspond
bijectively with the ``essential" representations of the
spectral \cstar\ $C^*(E)$ \cite{3, section 3}, \cite{4}, and \cite{5}.
This correspondence has
the feature that two representations of $C^*(E)$
are unitarily equivalent iff the corresponding
\esg s are conjugate.  In the case of decomposable
\esg s we obtain one of the
standard spectral \cstar s $C^*(E_n)$,
$n=1,2,\dots,\infty$.  Those \cstar s are continuous
analogues of the Cuntz algebra $\Cal O_\infty$ \cite{9},
and are far from being type I.
Nevertheless, if we agree to identify
two representations of $C^*(E_n)$ up to unitary equivalence
modulo cocycle perturbations (that is to say, up to
cocycle perturbations of the associated \esg s),
then the resulting set of equivalence classes
of representations is smooth: it
is parameterized by a single integer
$n=1,2,\dots,\infty$.

On the other hand, we remind the
reader that Powers has constructed
examples of \esg s that are of type II
(they have some intertwining semigroups but not
enough of them) \cite{27}, and others that are of type III
(they have no intertwining semigroups whatsoever) \cite{25}.
None of these more exotic \esg s can be decomposable.  The
structure of the product systems associated with such
\esg s remains quite mysterious.

\endremark

\vfill
\pagebreak

\end